\documentclass[12pt]{article}
\usepackage{a4wide,amsmath,amssymb}
\usepackage{color}

\usepackage{tikz}
\usetikzlibrary{calc}
\usepackage[ruled,vlined]{algorithm2e}
\usepackage{multirow}
\usepackage{comment}

\newtheorem{theorem}{Theorem}[section]

\newtheorem{lemma}[theorem]{Lemma}
\newtheorem{proposition}{Proposition}
\newtheorem{assumption}{Assumption}
\newtheorem{remark}{Remark}

\def\calR{\mathcal{R}}
\def\calT{\mathcal{T}}
\def\ul#1{\underline{#1}}
\def\ol#1{\overline{#1}}
\def\ProjMdel{P_{\tilde{M}^\delta}}
\newcommand{\dd}{{~\mathrm{d}}}
\newcommand{\ad}{{\mathrm{ad}}}
\newcommand{\ele}{{\mathrm{element}}}
\newcommand{\nod}{{\mathrm{node}}}
\newcommand{\meas}{{\mathrm{meas}}}
\newcommand{\ex}{{\mathrm{ex}}}

\DeclareMathOperator*{\argmin}{argmin}

\begin{document}
\title{Iterative regularization for constrained minimization formulations of nonlinear inverse problems}
\author{Barbara Kaltenbacher \and Kha Van Huynh\\
Department of Mathematics, Alpen-Adria-Universität Klagenfurt, Austria\\
              {barbara.kaltenbacher@aau.at, van.huynh@aau.at} 
}

\maketitle

\begin{abstract}
In this paper we the formulation of inverse problems as constrained minimization problems and their iterative solution by gradient or Newton type. We carry out a convergence analysis in the sense of regularization methods and discuss applicability to the problem of identifying the spatially varying diffusivity in an elliptic PDE from different sets of observations. Among these is a novel hybrid imaging techology known as impedance acoustic tomography, for which we provide numerical experiments.  
\end{abstract}

\noindent
{\bf key words:} inverse problems \and iterative regularization \and coefficient identification in elliptic PDEs \and impedance acoustic tomography

\section{Introduction}\label{sec:intro}
Inverse problems usually consist of a model
\begin{equation}\label{Axu0}
A(x,u)=0
\end{equation}
where the operator $A$ acts on the state $u$ of a system and contains unknown parameters $x$, 
and an observation equation
\begin{equation}\label{Cxuy}
C(x,u)=y
\end{equation}
quantifying the additionally available information that is supposed to allow for identifying the parameters $x$; by a slight notation overload, we will often summarize $(x,u)$ into a single element, which we again call $x$.

The classical formulation of an inverse problem is as an operator equation 
\begin{equation}\label{Fxy}
F(x)=y
\end{equation}
where usually $x$ is the searched for parameter (some coefficient, initial or boundary conditions in a PDE or ODE model) but -- in an all-at-once formulation -- might as well include the state, i.e., the PDE solution. 
In a conventional reduced setting $F=C\circ S$ is the concatenation of an observation operator $C$ with a parameter-to-state map $S$ satisfying $A(x,S(x))=0$, whereas an all-at-once setting considers the inverse problem as a system $\left\{\begin{array}{c} A(x,u)=0\\C(x,u)=y\end{array}\right.$, which by the above mentioned replacement $x:=(x,u)$ takes the form \eqref{Fxy}, see, e.g. \cite{all-at-once,aao_time}.

We here follow the idea of generalizing this to a formulation of an inverse problem as a constrained minimization problem 
\begin{equation}\label{minJ}
\min J(x) \mbox{ s.t. } x\in M\,,
\end{equation}
where in a reduced type setting, $x$ is the parameter and in an all-at-once-setting $x=(x,u)$ contains both parameter and state. In what follows, it will not be necessary to distinuish between these two cases notationally.

Straightforward instances for equivalent minimization based formulations of \eqref{Axu0}, \eqref{Cxuy} are, e.g.,
\[
\min \tfrac12 \|C(x,u)-y\|^2\mbox{ s.t. }A(x,u)=0\,,
\]
\begin{equation}\label{minAstC}
\min \tfrac12 \|A(x,u)\|^2\mbox{ s.t. }C(x,u)=y\,,
\end{equation}
or in the context of \eqref{Fxy} comprising both the reduced $F(x)=C(S(x))$ and the all-at-once $F(x,u)=\left(\begin{array}{c} A(x,u)\\C(x,u)\end{array}\right)$ setting simply 
\begin{equation}\label{nllsq}
\min \tfrac12 \|F(x)-y\|^2\,.
\end{equation}
For further examples of such formulations, see., e.g., \cite{minIP,examples_minIP}. In particular we point to the variational formulation according to Kohn and Vogelius, see, e.g, \cite{KohnVogelius87}.

Here $J$ is a proper functional acting on a Banach space $X$, and we make the normalization assumption 
\begin{equation}\label{normalize}
J\geq0 \mbox{ on } M \mbox{ and }J(x^\dagger)=\min_{x\in M}J(x)=0
\end{equation}
for $x^\dagger$ solving the inverse problem, i.e., we assume to know the minimal value of $J$ (but of course not the minimizer, which is what we intend to retrieve).

The first order optimality condition for a minimizer of \eqref{minJ} is
\begin{equation}\label{stat0}
\langle \nabla J(x^\dagger),x-x^\dagger\rangle\geq0
\quad \mbox{ for all } x\in M\,,
\end{equation}

Typically inverse problems also in the form \eqref{minJ} are ill-posed in the sense that solutions to \eqref{minJ} do not depend continuously on the data $y$ that enters the definition of the cost function $J$ and/or of the feasible set $M$. Since in a realistic setting the data is contamianted with measurement noise, i.e., only $y^\delta\approx y$ is given, regularization needs to be employed.
We first of all do so by possibly adding some regularizing constraints -- in particular we think of bound constraints in the sense of Ivanov regularization --  and/or by relaxing constraints like fit to the data in the sense of Morozov regularization. In the context of \eqref{minAstC}, this, e.g., means that we replace $M=\{x\in X\,:\, Cx=y\}$ by  
$\tilde{M}^\delta=\{x\in X \ : \|Cx-y^\delta\|\leq\tau \delta\mbox{ and } \tilde{\calR}(x)\leq\rho\}$, for the noise level $\delta\geq\|y-y^\delta\|$, some constants $\tau>1$, $\rho>0$ and some functional $\tilde{\calR}$ satisfying $\tilde{\calR}(x^\dagger)\leq\rho$. 

Thus we consider the partly regularized problem
\begin{equation}\label{minJdelta}
\min J^\delta(x) \mbox{ s.t. } x\in \tilde{M}^\delta
\end{equation}
which we intend to solve iteratively, where further regularization is incorporated by early stopping and potentially also by adding regularizing terms during the iteration. 
As in the above example of $\tilde{M}^\delta$, we will generally assume $x^\dagger$ to be feasible also for this modified problem, and also approximately minimal 
\begin{equation}\label{etadelta}
\begin{aligned}
&x^\dagger\in \tilde{M}^\delta \mbox{ and } J^\delta(x^\dagger)\leq \eta(\delta) 
\quad \mbox{ for all }\delta\in(0,\bar{\delta})\,, \\
&\mbox{ where } \eta(\delta)>0 \mbox{ and }\eta(\delta)\to0\mbox{ as }\delta\to0\,.
\end{aligned}
\end{equation}
With \eqref{minJdelta} we formally stay in the same setting as in \eqref{minJ} and, like in \eqref{normalize}, assume 
\begin{equation}\label{normalize_delta}
J^\delta\geq0 \mbox{ on } \tilde{M}^\delta \,.
\end{equation}
The key difference to \eqref{minJ} lies in the fact that $J^\delta$ and $\tilde{M}^\delta$ might depend on the noise level and this will in fact be crucial since we will study convergence as $\delta$ tends to zero. 

Since we consider formulations of inverse problems as constrained minimization problems, an essential step is to consider extensions of iterative methods such as gradient or Newton type methods, to constrained minimization problems.
Along with these two paradigms concerning the search direction, we will consider two approaches for guaranteeing feasibility of the sequence, namely projection onto the admissible set in the context of gradient methods in Section \ref{sec:grad} and sequential quadratic programming SQP type constrained minimization in Section \ref{sec:Newton}. 

Some key reference for gradient, i.e., Landweber type iterative methods are \cite{Eicke92} on projected Landweber iteration for linear inverse problems, \cite{HNS95} on (unconstrained) nonlinear Landweber iteration Landweber, and more recently \cite{Kindermann17} on gradient type methods under very general conditions on the cost function or the forward operator, respectively.
Extensions with a penalty term (also allowing for the incorporation of constraints) for linear inverse problems can be found in \cite{BotHein12}; For nonlinear problems we also point to \cite{JinWang13,WangWangHan19}, however, they do not seem to be applicable to constrained problems, since the penalty term is assumed to be $p$-convex and thus cannot be an indicator function.

Newton type methods for the solution of nonlinear ill-posed problems have been extensively studied in Hilbert spaces (see, e.g., \cite{BakKok04,KalNeuSch08} and the references therein) and more recently also in a in Banach space setting. In particular, the iteratively regularized Gauss-Newton method \cite{Baku92} or the Levenberg-Marquardt method \cite{HankeLevMar} easily allow to incorporate constraints in their variational form.
Projected Gauss-Newton type methods for constrained ill-posed problems have been considered in, e.g., \cite{muge12}.

The remainder of this paper is organized as follows. 
In Section \ref{sec:grad} we will study a projected version of Landweber iteration, thus a gradient type method in a Hilbert space setting and prove its convergence under certain convexity assumptions on the cost function. 
Section \ref{sec:Newton} turns to a general Banach space setting and discusses Newton SQP methods as well as their convergence.
Finally, in Section we investigate applicability to the identification of a spatially varying diffusion coefficient in an elliptic PDE from diefferent sets of boundary condtions which leads to three different inverse probems: Inverse groundwater filtration (often also used as a model problem and denoted by $a$-problem) impedance acoustic tomography and electrical impedance tomography. Numerical experiments in Section \ref{sec:numexp_EIT-IAT} illustrate our theoretical findings. 

\section{A projected gradient method}\label{sec:grad}
In this section, we consider the projected gradient method for \eqref{minJdelta}  
\begin{equation}\label{projgrad}
\tilde{x}_{k+1}=x_k-\mu_k\nabla J^\delta(x_k)\,, \quad x_{k+1}=\ProjMdel(\tilde{x}_{k+1})
\end{equation}
and extend some of the results from \cite{Kindermann17} to the constrained setting, or from a different viewpoint, extend some of the results from \cite{Eicke92} to the nonlinear setting. 
In \eqref{projgrad}, $\mu_k>0$ is a stepsize parameter and $\nabla J^\delta(x_k)\in X$ is the Riesz representation of ${J^\delta}'(x_k)\in X^*$ as in this section we restrict ourselves to a Hilbert space setting. The reason for this is the fact that in general Banach spaces, ${J^\delta}'(x_k)$ would have to be transported back into $X$ by some duality mapping, which adds nonlinearity and therefore, among others, complicates the choice of the step size, see e.g. \cite{itnlBanach} for the unconstrained least squares case \eqref{nllsq}.
Moreover, throughout this section we will assume $\tilde{M}^\delta$ to be closed and convex and denote by $\ProjMdel$ the metric (in the Hilbert space setting condsidered in this section also orthogonal) projection onto $\tilde{M}^\delta$, which is characterized by the variational inequality
\begin{equation}\label{proj}
x=\ProjMdel(\tilde{x})\ \Leftrightarrow \ \left(x\in \tilde{M}^\delta \mbox{ and }
\forall z\in \tilde{M}^\delta \, : \  \langle \tilde{x}-x,z-x\rangle \leq0 \right)
\end{equation}
With $z:=x_k\in \tilde{M}^\delta$, this immediately implies 
\[
0\geq \langle \tilde{x}_{k+1}-x_{k+1},x_k-x_{k+1}\rangle 
= \langle x_k-x_{k+1}-\mu_k \nabla J^\delta(x_k),x_k-x_{k+1}\rangle
\]
hence 
\begin{equation}\label{estgrad1}
\|x_{k+1}-x_k\|^2\leq -\mu_k \langle \nabla J^\delta(x_k),x_{k+1}-x_k\rangle
\end{equation}
and thus, using the Cauchy-Schwarz inequality, the estimate
\begin{equation}\label{Deltak_nablaJ}
\|x_{k+1}-x_k\| \leq \mu_k \| \nabla J^\delta(x_k)\|\,.
\end{equation}
Moreover, as well known for (projected) gradient methods, under the Lipschitz type condition on the gradient 
\begin{equation}\label{Lipschitz}
J^\delta(x)-J^\delta(x_+)-\langle \nabla J^\delta(x)(x-x_+)\geq -\tfrac{L}{2}\|x-x_+\|^2 
\quad \mbox{ for all } x,x_+\in \tilde{M}^\delta
\end{equation}
for $\mu_k\leq\ol{\mu}<\frac{2}{L}$, from \eqref{estgrad1} we get monotonicity of the cost function values
\[
J^\delta(x_k)-J^\delta(x_{k+1})\geq (\tfrac{1}{\mu_k}-\tfrac{L}{2}) \|x_{k+1}-x_k\|^2
\]
and square summability of the steps
\[
\sum_{k=0}^\infty \|x_{k+1}-x_k\|^2 \leq\frac{1}{\tfrac{1}{\ol{\mu}}-\tfrac{L}{2}} J^\delta(x_0)\,.
\]
Monotonicity of the error under additional convexity assumptions easily follows from nonexpansivity of the projection, which yields
\begin{equation}\label{estgrad2}
\begin{aligned} 
&\|x_{k+1}-x^\dagger\|^2-\|x_k-x^\dagger\|^2
=\|\ProjMdel(\tilde{x}_{k+1})-\ProjMdel(x^\dagger)\|^2-\|x_k-x^\dagger\|^2 \\
&\leq\|\tilde{x}_{k+1}-x^\dagger\|^2-\|x_k-x^\dagger\|^2 
=\|\tilde{x}_{k+1}-x_k\|^2+2\langle \tilde{x}_{k+1}-x_k,x_k-x^\dagger\rangle\\
&=\mu_k^2\|\nabla J^\delta(x_k)\|^2-2 \mu_k\langle \nabla J^\delta(x_k),x_k-x^\dagger\rangle\,.
\end{aligned}
\end{equation}
Under the monotonicity condition on $\nabla J^\delta$ (i.e., convexity condition on $J^\delta$) 
\begin{equation}\label{convex1}
\langle \nabla J^\delta(x)- \nabla J^\delta(x^\dagger),x-x^\dagger\rangle
\geq\gamma \|\nabla J^\delta(x)\|^2
\quad \mbox{ for all } x\in \tilde{M}^\delta
\end{equation}
(which for $\gamma=0$ follows from convexity of $J^\delta$, i.e., monotonicity of $\nabla J^\delta$) 
and assuming approximate stationarity 
\begin{equation}\label{stat}
\langle \nabla J^\delta(x^\dagger),x-x^\dagger\rangle\geq-\eta(\delta)
\quad \mbox{ for all } x\in \tilde{M}^\delta\,,
\end{equation}
we get from \eqref{estgrad2}, that for all $k\leq k_*-1$ with $k_*$ defined by 
\begin{equation}\label{discrprinc_grad}
k_*=k_*(\delta)=\min\{k\, : \, \|\nabla J^\delta(x_k)\|^2\leq\tau\eta(\delta)\}
\end{equation}
the estimate
\begin{equation}\label{estgrad3}
\begin{aligned}
&\|x_{k+1}-x^\dagger\|^2-\|x_k-x^\dagger\|^2\\
&\leq \mu_k^2\|\nabla J^\delta(x_k)\|^2-2\mu_k\langle \nabla J^\delta(x_k)-\nabla J^\delta(x^\dagger)\rangle +2\mu_k\eta(\delta)\\
&\leq -\mu_k(2-\tfrac{\mu_k}{\gamma}-\tfrac{2}{\tau\gamma})
\langle \nabla J^\delta(x_k)- \nabla J^\delta(x^\dagger),x_k-x^\dagger\rangle\\
&\leq -\mu_k(2\gamma-\mu_k-\tfrac{2}{\tau})\|\nabla J^\delta(x_k)\|^2 \leq 0
\end{aligned}
\end{equation}
for $\tau>\frac{1}{\gamma}$, $0<\underline{\mu}\leq\mu_k\leq\bar{\mu}<2(\gamma-\frac{1}{\tau})$, 
hence summability
\[
\sum_{k=0}^{k_*} \langle \nabla J^\delta(x_k)- \nabla J^\delta(x^\dagger),x_k-x^\dagger\rangle
\leq \frac{1}{\underline{\mu}(2-\frac{\bar{\mu}}{\gamma}-\tfrac{2}{\tau\gamma})}\|x_0-x^\dagger\|^2\,.
\]
\begin{equation}\label{sum1}
\sum_{k=0}^{k_*}\|\nabla J^\delta(x_k)\|^2 
\leq \frac{1}{\underline{\mu}(2\gamma-\bar{\mu}-\tfrac{2}{\tau})} \|x_0-x^\dagger\|^2\,, 
\end{equation}
Alternatively, under a condition following from \eqref{convex1}, \eqref{stat} and comprising both convexity and approximate stationarity
\begin{equation}\label{convex2}
\langle \nabla J^\delta(x),x-x^\dagger\rangle
\geq\gamma \|\nabla J^\delta(x)\|^2-\eta(\delta)
\quad \mbox{ for all } x\in \tilde{M}^\delta
\end{equation}
which for $k\leq k_*-1$ implies 
\[
(\gamma\tau-1)\eta(\delta)\leq \langle\nabla J^\delta(x_k),x_k-x^\dagger\rangle,
\] 
as well as
\begin{equation}\label{estgrad4}
(1+\tfrac{1}{\gamma\tau-1}) \langle\nabla J^\delta(x_k),x_k-x^\dagger\rangle 
\geq \gamma \|\nabla J^\delta(x_k)\|^2\,,
\end{equation}
we get from \eqref{estgrad2} 
\begin{equation}\label{estgrad5}
\begin{aligned}
&\|x_{k+1}-x^\dagger\|^2-\|x_k-x^\dagger\|^2
\leq -\mu_k(2-\tfrac{\mu_k}{\gamma}(1+\tfrac{1}{\gamma\tau-1}))
\langle \nabla J^\delta(x_k),x_k-x^\dagger\rangle\\
&\leq -\mu_k\left(\tfrac{2\gamma}{1+\tfrac{1}{\gamma\tau-1}}-\mu_k\right)
\|\nabla J^\delta(x_k)\|^2 \leq 0
\end{aligned}
\end{equation}
for $\tau>\frac{1}{\gamma}$, $0<\underline{\mu}\leq\mu_k\leq \bar{\mu} < \tfrac{2\gamma}{1+\tfrac{1}{\gamma\tau-1}}$
hence summability
\[
\sum_{k=0}^{k_*} \langle \nabla J^\delta(x_k),x_k-x^\dagger\rangle
\leq \frac{1}{\underline{\mu}(2-\tfrac{\bar{\mu}}{\gamma}(1+\tfrac{1}{\gamma\tau-1}))}
\|x_0-x^\dagger\|^2\,.
\]
which via \eqref{estgrad4} also implies summability of $\|\nabla J^\delta\|^2$.
\begin{equation}\label{sum2}
\sum_{k=0}^{k_*} \|\nabla J^\delta\|^2
\leq (1+\tfrac{1}{\gamma\tau-1})
\frac{1}{\gamma\underline{\mu}(2-\tfrac{\bar{\mu}}{\gamma}(1+\tfrac{1}{\gamma\tau-1}))}
\|x_0-x^\dagger\|^2\,.
\end{equation}

\medskip

The estimates \eqref{sum1}, \eqref{sum2} imply convergence of the gradient to zero as $k\to\infty$ in the noise free case and finiteness of the stopping index $k_*$ in case of noisy data.
In the noiseless case $\delta=0$ Opial's Lemma (Lemma \ref{lem:opialdiscr} in the Appendix) with $S=\{x^*\in X\,:\,\forall x\in M:\,\langle \nabla J(x^*),x-x^*\rangle\geq0\}$, due to monotonicity of $\|x_k-x^*\|$ and the Bolzano-Weierstrass Theorem, implies weak convergence of $x_k$ as $k\to\infty$ to a stationary point.
In case of noisy data, one could think of applying the continuous version of Opial's Lemma (Lemma \ref{lem:opialcont} in the Appendix) with $t:=\frac{1}{\delta}$, $x(t):=x_{k_*(\delta)}$. However, we do not have monotonicity of the final iterates $x_{k_*(\delta)}$ as a function of $\delta$. Still, in case of uniqueness, that is, if $S$ is a singleton $S=\{x^\dagger\}$, then boundedness of the sequence $\|x_{k_*(\delta)}-x^*\|$ by $\|x_0-x^*\|$ together with a subsequence-subsequence argument yields its weak convergence of $x_{k_*(\delta)}$ to $x^\dagger$ as $\delta\to0$.

For this purpose, we have to impose certain continuity assumptions on the cost function and the constrains, namely
\begin{eqnarray}\label{MJcontinuous}
&&\mbox{For any sequence }(z_n)_{n\in\mathbb{N}}\subseteq X\,, \ (\delta_n)_{n\in\mathbb{N}}\in(0,\bar{\delta}]\,, \ \delta_n\to0\mbox{ as }n\to\infty\nonumber\\
&&\Bigl(\forall n\in\mathbb{N}: \ z_n\in \tilde{M}^{\delta_n} \mbox{ and } z_n \rightharpoonup z \mbox{ and } \nabla J^{\delta_n}(z_n)\to0\Bigr) \nonumber\\ 
&&\qquad\qquad \Rightarrow \
\Bigl(z\in M \mbox{ and } \forall x\in M:\,\langle \nabla J(z),x-z\rangle\geq0\Bigr)
\end{eqnarray}
which in the noiseless case becomes
\begin{equation}\label{MJcontinuous0}
\begin{aligned}
&\mbox{For any sequence }(z_n)_{n\in\mathbb{N}}\subseteq X\\
& \Bigl(\forall n\in\mathbb{N}: \ z_n\in M \mbox{ and } z_n \rightharpoonup z \mbox{ and } \nabla J(z_n)\to0\Bigr)
\\ 
&\qquad \qquad\Rightarrow \
\Bigl(z\in M \mbox{ and } \forall x\in M:\,\langle \nabla J(z),x-z\rangle\geq0\Bigr)
\end{aligned}
\end{equation}

\begin{proposition}\label{lem:convLW}
Let \eqref{normalize}, \eqref{normalize_delta}, \eqref{convex2} hold, and let the sequence of iterates $x_k$ be defined by \eqref{projgrad} with $k_*$ defined by \eqref{discrprinc_grad}. 

Then for $\delta=0$ if $M$ and $\nabla J$ satisfy \eqref{MJcontinuous0}, the sequence $(x_k)_{k\in\mathbb{N}}$ converges weakly to a solution $x^*\in M$ of the first order optimality condition \eqref{stat0} as $k\to\infty$.

If $\delta>0$ and additionally \eqref{etadelta}, \eqref{stat}, and \eqref{MJcontinuous} holds, then the family $(x_{k_*(\delta)})_{\delta\in(0,\bar{\delta}]}$ converges weakly subsequentially to a stationary point $x^\dagger$ according to \eqref{stat0} as $\delta\to0$. If this stationary point is unique, then the whole sequence converges weakly to  $x^\dagger$. The same assertion holds with stationarity \eqref{stat0} (with \eqref{MJcontinuous}) replaced by 
\begin{itemize}
\item[(a)] minimality, i.e., $x^\dagger$ (and $z$) $\in \mbox{argmin}\{J(x)\,:\,x\in M\}$
\item[or by]  
\item[(b)] $\|\nabla J(x^\dagger)\|=0$ (and $\|\nabla J(z)\|=0$).       
\end{itemize}
\end{proposition}
Note that case (a) makes uniqueness harder, whereas (b) makes uniqueness easier than \eqref{stat0}.

\begin{remark}\label{rem:strongconvLW} 
Strong convergence can be shown for the modified projected Landeweber method from \cite[Section 3.2]{Eicke92}.
However, this requires a source condition to hold.
\end{remark}

\begin{remark}\label{rem:convcond} 
Let us finally comment on the convexity condition \eqref{convex2}.

In the special case $J^\delta(x)= \tfrac12\|F(x)-y^\delta\|^2$ cf \eqref{nllsq}, condition \eqref{convex2} becomes
\begin{equation}\label{weaktcc_grad}
\langle F(x)-y^\delta,F'(x)(x-x^\dagger)\rangle
\geq\gamma \|F'(x)^*(F(x)-y^\delta)\|^2-\eta(\delta)
\end{equation}
which, follows, e.g., from
\begin{equation}\label{tccFxy_grad}
\begin{aligned}
&\|F'(x)\|\leq 1 \mbox{ and } \\
&\langle F(x)-F(x^\dagger)-F'(x)(x-x^\dagger),F(x)-y^\delta\rangle 
\leq (1-\gamma-\kappa)\|F(x)-y^\delta\|^2
\end{aligned}
\end{equation}
with $\|F(x^\dagger)-y^\delta\|\leq 4\kappa\eta(\delta)$.
The latter is closely related to the usual normalization and tangential cone conditions for Landweber iteration, see, e.g., \cite{HNS95,Kindermann17}.
It is, e.g., satisfied for linear $F$ as well as for some specific coefficient identification problems, see, e.g., \cite{HNS95} for the reduced setting, \cite{all-at-once} for the all-at-once setting, and \cite{TCC_time} for some time dependent problems in both reduced and all-at-once formulation.
\end{remark}

\section{An SQP type constrained Newton method}\label{sec:Newton}
A quadratic approximation of the cost function combined with a Tikhonov type additive regularization term yields the iteration scheme
\begin{equation}\label{Newton}
\begin{aligned}
&x_{k+1}\in X_{k+1}(\alpha):=\mbox{argmin}_{x\in \tilde{M}^\delta} Q_k^\delta(x) + \alpha_k\calR(x)\\
&\mbox{ where }Q_k^\delta(x)=J^\delta(x_k)+G^\delta(x_k)(x-x_k)+\tfrac12 H^\delta(x_k)(x-x_k)^2
\end{aligned}
\end{equation}
with 
\begin{equation}\label{GHR}
\begin{aligned}
&G^\delta(x_k):X\to\mathbb{R} \mbox{ linear }, \quad H^\delta(x_k):X^2\to\mathbb{R} \mbox{ bilinear },\\ 
&\mathcal{R}:X\to[0,\infty] \mbox{ proper with domain }\mbox{dom}(\calR)\supseteq \bigcup_{\delta\in(0,\bar{\delta})} \tilde{M}^\delta \cup M
\end{aligned}
\end{equation}
where $G^\delta$ and $H^\delta$ should be viewed as (approximations to) the gradient and Hessian of $J$, $G^\delta(x_k)\approx {J^\delta}'(x_k)$, $H^\delta(x_k)\approx {J^\delta}''(x_k)$, and $\calR$ is a regularization functional.
Since we do not necessarily neglect ${J^\delta}''(x_k)$, this differs from the iteratively regularized Gauss-Newton method IRGNM studied, e.g., in \cite{Baku92,IRGNMIvanov,NewtonIvanov}.

Here $X$ is a general Banach space.

To guarantee existence of minimizers, besides \eqref{GHR} we will make the following assumption 
\begin{assumption}\label{ass0}
For some topology $\calT_0$ on $X$
\begin{itemize}
\item for all $r\geq\mathcal{R}(x^\dagger)$, the sublevel set $\tilde{M}^\delta_r:=\{x\in M^\delta\, : \, \calR(x)\leq r\}$ is $\calT_0$ compact.
\item the mapping $Q_k^\delta+\alpha_k\calR$ is $\calT_0$ lower semicontinuous
\end{itemize}
\end{assumption}
Uniqueness of a minimizer of \eqref{Newton} will not necessarily hold; the sequence $(x_k)_{k\in\{1,\ldots,k_*\}}$ will therefore be defined by an arbitrary selection of minimizers of \eqref{Newton}.

The overall iteration is stopped according to the discrepancy principle
\begin{equation}\label{discrprinc}
k_*=k_*(\delta)=\min\{k\, : \, J^\delta(x_k)\leq\tau\eta(\delta)\}
\end{equation}
for some constant $\tau>1$.

As far as the sequence of regularization parameters $\alpha_k$ is concerned, we will choose it a priori or a posteriori, see \eqref{alpha_apriori}, \eqref{alpha_aposteriori} below.

A special case of this with 
\begin{equation}\label{JFxy}
\begin{aligned}
J^\delta(x)=\frac12\|F(x)-y^\delta\|^2, \quad 
G^\delta(x)h={J^\delta}'(x)h=\langle F(x)-y^\delta, F'(x)h\rangle, \\ 
{H^\delta}(x)(h,\ell)=\langle F'(x)h, F'(x)\ell\rangle
\end{aligned}
\end{equation}
(note that $H^\delta$ in general does not coincide with the Hessian of $J^\delta$) in Hilbert space is the iteratively regularized Gauss-Newton method for the operator equation formulation \eqref{Fxy} of the inverse problem, see, e.g., \cite{Baku92,IRGNMIvanov,NewtonIvanov}.
   
Another special case we will consider is the quadratic one
\begin{equation}\label{Jquad}
\begin{aligned}
J^\delta(x)=f^\delta+\bar{G}^\delta x + \tfrac12 \bar{H}^\delta x^2, \quad 
G^\delta(x)h=\bar{G}^\delta h, \quad  
{H^\delta}(x)(h,\ell)=\bar{H}^\delta(h,\ell)
\end{aligned}
\end{equation}
with $f^\delta\in X$, $\bar{G}^\delta\in L(X,\mathbb{R})=X^*$, $\bar{H}^\delta\in L(X^2,\mathbb{R})$, where trivially $Q_k^\delta$ coincides with $J^\delta$.

\bigskip

To provide a convergence analysis, we start with the case of an a priori choice of $\alpha_k$ 
\begin{equation}\label{alpha_apriori}
\alpha_k=\alpha_0\theta^k
\end{equation}
for some $\theta\in(0,1)$ and make, among others, the following assumption.
\begin{assumption}\label{ass1}
For some topology $\calT$ on $X$, 
\begin{itemize}
\item the sublevel set 
$\{x\in\bigcup_{\delta\in(0,\bar{\delta})}\tilde{M}^\delta\, : \, \calR(x)\leq R\}
= \bigcup_{\delta\in(0,\bar{\delta})}\tilde{M}_R^\delta$ 
is $\calT$ compact, with $R=(1+\tfrac{a-b}{\tau(a-b)-c} \tfrac{\tau b+c}{(a\theta-b)})\calR^\dagger + \tfrac{a-b}{\tau(a-b)-c} \tfrac{\tau b+c}{\alpha_0} J(x_0)$
with $\tau$ as in \eqref{discrprinc}, $\alpha_0$ as in \eqref{alpha_apriori}, and $a,b,c$ as in \eqref{abc1};
\item $M$ is $\calT$ closed with respect to the family of sets $(\tilde{M}^\delta)_{\delta\in(0,\bar{\delta})}$ in the following sense: 
\begin{eqnarray*}
&&\mbox{For any sequence }(z_n)_{n\in\mathbb{N}}\subseteq X\,, \ (\delta_n)_{n\in\mathbb{N}}\in(0,\bar{\delta}]\,, \ \delta_n\to0\mbox{ as }n\to\infty\nonumber\\
&&\Bigl(\forall n\in\mathbb{N}: \ z_n\in \tilde{M}^{\delta_n} \mbox{ and } z_n \stackrel{\calT}{\longrightarrow} z \Bigr) \Rightarrow z\in M
\end{eqnarray*}
\item
$\lim_{\delta\to0} \sup_{x\in M^\delta}(J(x)-J^\delta(x))\leq0$; 
\item
$J$ is $\calT$ lower semicontinuos.
\end{itemize}
\end{assumption}

Comparably to the tangential cone condition in the context of nonlinear Landweber iteration \cite{HNS95} and more recently also the IRGNM \cite{IRGNMIvanov} we impose a restriction on the nonlinearity / nonconvexity of $J$
\begin{equation}\label{abc1}
\begin{aligned}
G^\delta(x)(x_+-x^+)+\tfrac12 H^\delta(x)\Bigl((x_+-x)^2-(x-x^+)^2\Bigr)\geq a J^\delta(x_+)-bJ^\delta(x)-cJ^\delta(x^+)\\
\mbox{ for all } x,x_+\in \tilde{M}^\delta , \ x^+=x^\dagger \,, \quad 
\delta\in(0,\bar{\delta})\,,
\end{aligned}
\end{equation}
with $a>b\geq0$, $c\geq0$.

\begin{theorem}\label{th:apriori}
Let conditions \eqref{etadelta}, \eqref{normalize_delta}, \eqref{GHR}, \eqref{abc1}, and Assumptions \ref{ass0}, \ref{ass1} hold, assume that $\alpha_k$ is chosen a priori according to \eqref{alpha_apriori}, and $k_*$ is chosen according to the discrepancy principle \eqref{discrprinc}, with the following constraints on the constants
\[
1>\theta>\frac{b}{a}, \quad 
\tau>\frac{c}{a-b}.
\]
Then 
\begin{itemize}
\item For any $\delta\in(0,\bar{\delta})$, and any $x_0\in \bigcap_{\delta\in(0,\bar{\delta})}\tilde{M}^\delta\,\cap M$, 
\begin{itemize}
\item the iterates $x_k$ are well-defined for all $k\leq k_*(\delta)$ and $k_*(\delta)$ is finite;
\item for all $k\in\{1,\ldots, k_*(\delta)\}$ we have
\[ 
J^\delta(x_k)\leq \tfrac{b}{a} J^\delta(x_{k-1}) + \tfrac{1}{a} \alpha_k\calR^\dagger + \tfrac{c}{a}\eta; 
\]
\item for all $k\in\{1,\ldots, k_*(\delta)\}$ we have
\[
\calR(x_k)\leq R
\] 
\end{itemize}
\item As $\delta\to0$, the final iterates $x_{k_*(\delta)}$ tend to a solution of the inverse problem \eqref{minJ} $\calT$-subsequentially, i.e., every sequence $x_{k_*(\delta_j)}$ with $\delta_j\to0$ as $j\to\infty$ has a $\calT$ convergent subsequence and the limit of every $\calT$ convergent subsequence solves \eqref{minJ}. 
\end{itemize}
\end{theorem}
{\em Proof.}\
For any  $k\leq k_*-1$, existence of a minimizer follows from Assumption \ref{ass0} by the direct method of calculus of variations. To this end, note that by $x^\dagger\in \tilde{M}^\delta$, implying 
\[
\min_{x\in\tilde{M}^\delta} Q_k(x)+\alpha\calR(x)\leq Q_k(x^\dagger)+\alpha\calR(x^\dagger),
\] and the lower bound 
\[
Q_k(x)\geq Q_k(x^\dagger)+ a J^\delta(x)-bJ^\delta(x_k)-c J^\delta(x^\dagger)\geq Q_k(x^\dagger)-bJ^\delta(x_k)-c J^\delta(x^\dagger),
\] 
which yields 
\[
\tilde{M}^\delta_r\supseteq\{x\in \tilde{M}^\delta \, : \, Q_k(x)+\alpha\calR(x)\leq Q_k(x^\dagger)+\alpha\calR(x^\dagger)\}
\] 
for $r=\calR(x^\dagger)+\frac{1}{\alpha}(bJ^\delta(x_k)+cJ^\delta(x^\dagger))$, 
 it suffices to restrict the search for a minimizer to the set $\tilde{M}^\delta_r$ as defined in Assumption \ref{ass0}.

For a hence existing minimizer $x_{k+1}$, its minimality together with feasibility of $x^\dagger$ for \eqref{Newton} yields 
\begin{equation}\label{minimality}
\begin{aligned}
G^\delta(x_k)(x_{k+1}-x_k)+\tfrac12 H^\delta(x_k)(x_{k+1}-x_k)^2 + \alpha_k\calR(x_{k+1})\\
\leq 
G^\delta(x_k)(x^\dagger-x_k)+\tfrac12 H^\delta(x_k)(x^\dagger-x_k)^2 + \alpha_k\calR(x^\dagger),
\end{aligned}
\end{equation}
which with \eqref{abc1} implies
\begin{equation}\label{minimality_abc}
a J^\delta(x_{k+1}) + \alpha_k\calR(x_{k+1}) \leq b J^\delta(x_k) + cJ^\delta(x^\dagger) +\alpha_k\calR(x^\dagger)
\end{equation}
thus, with the a priori choice \eqref{alpha_apriori}, and \eqref{etadelta}, abbreviating $J_k=J^\delta(x_k)$, $\calR_k=\calR(x_k)$, $\calR^\dagger=\calR(x^\dagger)$
\begin{equation}\label{minimality_abc_apriori}
J_{k+1} + \tfrac{\alpha_0}{a}\theta^k\calR_{k+1} \leq \tfrac{b}{a} J_k + \tfrac{\alpha_0}{a} \theta^k\calR^\dagger + \tfrac{c}{a}\eta \,.
\end{equation}
Inductively, with $\calR\geq0$, we conclude that for all $k\leq k_*$
\begin{equation}\label{Jk}
\begin{aligned}
J_k&\leq (\tfrac{b}{a})^k J_0 + \tfrac{\alpha_0}{a}\calR^\dagger\sum_{j=0}^{k-1} (\tfrac{b}{a})^j\theta^{k-1-j}+ \tfrac{c}{a}\eta \sum_{j=0}^{k-1} (\tfrac{b}{a})^j\\
&\leq (\tfrac{b}{a})^k J_0 + \tfrac{\alpha_0}{a\theta-b}\calR^\dagger\theta^k+ \tfrac{c}{a-b}\eta\,. 
\end{aligned}
\end{equation}
Using the minimality of $k_*$ according to \eqref{discrprinc}, we get, for all $k\leq k_*-1$, that 
$\eta\leq\frac{J_k}{\tau}$ and therefore, together with \eqref{Jk} 
\[
(1-\tfrac{c}{\tau(a-b)}) J_k \leq (\tfrac{b}{a})^k J_0 + \tfrac{\alpha_0}{a\theta-b}\calR^\dagger\theta^k\,.
\]
Inserting this back into \eqref{minimality_abc_apriori} with $J^\delta\geq0$, after multiplication by $\frac{a}{\alpha_k}$ and again using \eqref{discrprinc} yields
\begin{equation}\label{estR}
\begin{aligned}
\calR_{k+1}&\leq \tfrac{b}{\alpha_k} J_k + \calR^\dagger +\tfrac{c}{\alpha_k} \tfrac{J_k}{\tau}
\leq \calR^\dagger +\tfrac{\tau b+c}{\tau \alpha_0} \theta^{-k} J_k\\
&\leq \calR^\dagger + \tfrac{a-b}{\tau(a-b)-c} \tfrac{\tau b+c}{\alpha_0} \Bigl(
(\tfrac{b}{a\theta})^k J_0 + \tfrac{\alpha_0}{a\theta-b}\calR^\dagger\Bigr)= R
\end{aligned}
\end{equation}
for all $k\leq k_*-1$.

From \eqref{Jk}, which holds for all $k\leq k_*$ and $\tau>\tfrac{c}{a-b}$, as well as $\tfrac{b}{a}<\theta$, we conclude that the stopping index according to \eqref{discrprinc} is reached after finitely many, namely at most $\frac{\log((\tau-\tfrac{c}{a-b})\eta)-\log(J_0+\tfrac{\alpha_0}{a\theta-b}\calR^\dagger)}{\log\theta}$ steps.

Setting $k=k_*-1$ in \eqref{estR} yields $\calR(x_{k_*(\delta)})\leq R$, 
which implies $\calT$ convergence of a subsequence $x^j$ of $x_{k_*(\delta)}$ to some $\bar{x}$, which by Assumption \ref{ass1} lies in $M$.

By definition of $k_*$ and \eqref{etadelta} we have $J(x_{k_*(\delta)})\leq\tau\eta(\delta)+J(x_{k_*(\delta)})-J^\delta(x_{k_*(\delta)})\to0$ as $\delta\to0$; $\calT$ lower semicontinuity therefore yields $J(\bar{x})=0$.

\begin{flushright}
$\diamondsuit$
\end{flushright}

\bigskip
We now consider convergence with an a posteriori according choice of $\alpha_k$ according to the discrepancy principle type rule (which can also be interpreted as an inexact Newton condition) 
\begin{equation}\label{alpha_aposteriori}
\ul{\sigma}\leq\sigma_k(\alpha_k):=\frac{Q_k^\delta(X_{k+1}(\alpha_k))}{J^\delta(x_k)}\leq\ol{\sigma}
\end{equation}
with $0<\ul{\sigma}<\ol{\sigma}<1$; note that in \eqref{alpha_aposteriori}, the denominator of $\sigma_k(\alpha_k)$ will be positive and bounded away from zero by $\tau\eta(\delta)$ for all $k\leq k_*(\delta)-1$ by \eqref{discrprinc}.
In order to obtain well-definedness of $\sigma_k(\alpha)$ as a function of $\alpha$, we will assume that the mapping 
\[
\alpha\mapsto Q_k^\delta(X_{k+1}(\alpha)) \mbox{ with }X_{k+1}(\alpha)=\mbox{argmin}_{x\in\tilde{M}^\delta}(Q_k^\delta(x)+\alpha\calR(x))
\] 
is single valued, which is, e.g., the case if the minimizer of $Q_k^\delta(x)+\alpha\calR(x)$ over $\tilde{M}^\delta$ is unique. The latter can be achieved, e.g., by assuming convexity of $Q_k^\delta$ -- choosing $H^\delta$ as a positive semidefinite approximation of the (not necessarily positive semidefinite) true Hessian ${J^\delta}''$ -- and strict convexity of $\calR$.

For the a posteriori choice \eqref{alpha_aposteriori} we have to slightly modify the setting to guarantee existence of $\alpha_k$ such that \eqref{alpha_aposteriori} holds. The latter is possible if for some appropriate point $x^*$, 
the quotient $\frac{Q^\delta_k(x^*)}{J^\delta(x_k)}$ is large enough 
\begin{equation}\label{casea}
\ul{\sigma}<\frac{Q^\delta_k(x^*)}{J^\delta(x_k)}
\end{equation}
as we will show below.
This leads us to the following case distinction for updating the iterates
\begin{eqnarray*}
&&\mbox{If \eqref{casea} holds, choose $\alpha_k$ according to \eqref{alpha_aposteriori} and $x_{k+1}$ as in \eqref{Newton}} \\
&&\mbox{otherwise set $x_{k+1}=x^*$.} 
\end{eqnarray*}
Here $x^*\in \bigcap_{\delta\in(0,\bar{\delta})} \tilde{M}^\delta \cap M$ is a point of attraction of $\calR$ in the sense of the following assumption.
\begin{assumption}\label{ass2}
For some topology $\calT_1$ on $X$, 
\begin{itemize}
\item 
$\calR(x^*)=0$ and  for any sequence $(x_j)_{j\in\mathbb{N}}\subseteq X$ 
\begin{equation}\label{coercivity}
\calR(x_j)\to0 \ \Rightarrow \ x_j \stackrel{\calT_1}{\longrightarrow} x^*
\end{equation}
\item sublevel sets of $\calR$ are $\calT_1$ compact;
\item $\calR$ is $\calT_1$ lower semicontinuous;
\item the mapping $x\mapsto G^\delta(x_k)(x-x_k)+\tfrac12 H^\delta(x_k)(x-x_k)^2$ is $\calT_1$ continuous
\item $\tilde{M}^\delta$ is $\calT_1$ closed.
\end{itemize}
\end{assumption}
A simple example of a functional $\calR$ satisfying this assumption is some power of the norm distance from the a priori guess $x^*$, $\calR(x)=\|x-x^*\|^p$, along with the weak or weak* topology $\calT_1$, provided $X$ is reflexive or the dual of a separable space.

\begin{lemma}\label{lem:cont}
The mappings $\alpha\mapsto\calR(x_{k+1}(\alpha))$ and $\alpha\mapsto -Q_k(x_{k+1}(\alpha))$, where $x_{k+1}(\alpha)\in X_{k+1}(\alpha)$ (cf. \eqref{Newton}) are monotonically decreasing.

If additionally Assumption \ref{ass2} and \eqref{abc1} hold, and the mapping $\alpha\mapsto Q_k(X_{k+1}(\alpha))$ is single valued, then the mapping $\alpha\mapsto \sigma_k(\alpha)$ is well-defined and continuous on $(0,\infty)$.
\end{lemma}
{\em Proof.}\
For two values $\alpha$, $\tilde{\alpha}$, minimality implies
\[
\begin{aligned}
&Q_k(x_{k+1}(\alpha))+\alpha\calR(x_{k+1}(\alpha))
\leq Q_k(x_{k+1}(\tilde{\alpha}))+\alpha\calR(x_{k+1}(\tilde{\alpha}))\\
&= Q_k(x_{k+1}(\tilde{\alpha}))+\tilde{\alpha}\calR(x_{k+1}(\tilde{\alpha}))
+(\alpha-\tilde{\alpha})\calR(x_{k+1}(\tilde{\alpha}))\\
&\leq Q_k(x_{k+1}(\alpha))+\tilde{\alpha}\calR(x_{k+1}(\alpha))
+(\alpha-\tilde{\alpha})\calR(x_{k+1}(\tilde{\alpha}))
\end{aligned}
\]
which implies
\[
0\geq (\alpha-\tilde{\alpha})(\calR(x_{k+1}(\alpha))-\calR(x_{k+1}(\tilde{\alpha}))).
\]
Hence, $\alpha\mapsto\calR(x_{k+1}(\alpha))$ is monotonically decreasing and 
\[
Q_k(x_{k+1}(\alpha))-Q_k(x_{k+1}(\tilde{\alpha}))
\leq\alpha(\calR(x_{k+1}(\tilde{\alpha}))-\calR(x_{k+1}(\alpha)))\leq0 \mbox{ for }\alpha\leq\tilde{\alpha}
\]
that is, $\alpha\mapsto Q_k(x_{k+1}(\alpha))$ is monotonically increasing.

To prove continuity of the mapping $\alpha\mapsto Q_k(x_{k+1}(\alpha))$ under the assumption that this mapping is single valued, consider $\bar{\alpha}>0$ and a sequence $(\alpha_\ell)_{\ell\in\mathbb{N}}$ converging to $\bar{\alpha}>0$.
Minimality and \eqref{abc1} yield
\[
\begin{aligned}
&\alpha_\ell \calR(x_{k+1}(\alpha_\ell))
\leq Q_k(x^\dagger) +\alpha_\ell\calR(x^\dagger) - Q_k(x_{k+1}(\alpha_\ell))\\
&\leq -a J^\delta(x_{k+1}(\alpha_\ell))+b J^\delta(x_k)+c J^\delta(x^\dagger)+\alpha_\ell\calR(x^\dagger)
\leq b J^\delta(x_k)+c J^\delta(x^\dagger)+\alpha_\ell\calR(x^\dagger),
\end{aligned}
\]
which by strict positivity of $\bar{\alpha}$ implies boundedness of $(\calR(x_{k+1}(\alpha_\ell)))_{\ell\in\mathbb{N}}$. By Assumption \ref{ass2} there exists a $\calT_1$ convergent subsequence $(x_{k+1}(\alpha_{\ell_j}))_{j\in\mathbb{N}}$ whose limit $\bar{x}$ lies in $\tilde{M}^\delta$ and even in $X_{k+1}(\alpha)$, due to the fact that $Q_k(x_{k+1}(\alpha_{\ell_j}))\to Q_k(\bar{x})$ and the estimate
\[
\begin{aligned}
&Q_k(\bar{x}) +\bar{\alpha}\calR(\bar{x}) 
\leq \liminf_{j\to\infty} \Bigl(Q_k(x_{k+1}(\alpha_{\ell_j})) +\alpha_{\ell_j} \calR(x_{k+1}(\alpha_{\ell_j}))\Bigr)\\
&\leq \liminf_{j\to\infty} \Bigl(Q_k(x_{k+1}(\bar{\alpha})) +\alpha_{\ell_j} \calR(x_{k+1}(\bar{\alpha}))\Bigr)
= Q_k(x_{k+1}(\bar{\alpha})) +\bar{\alpha} \calR(x_{k+1}(\bar{\alpha}))\,.
\end{aligned}
\]
A subsequence-subsequence argument together with $\calT_1$ continuity of $Q_k$ and the assumed single valuedness of the mapping $\alpha\mapsto Q_k(X_{k+1}(\alpha))$ implies convergence 
$Q_k(x_{k+1}(\alpha_\ell))\to Q_k(x_{k+1}(\bar{\alpha}))$, hence, after division by $J^\delta(x_k)$, convergence $\sigma_k(x_{k+1}(\alpha_\ell))\to \sigma_k(x_{k+1}(\bar{\alpha}))$.
\begin{flushright}
$\diamondsuit$
\end{flushright}

To prove convergence of the iterates, we need a slightly stronger condition than \eqref{abc1}, namely
\begin{equation}\label{abc2}
\begin{aligned}
\ul{a}J^\delta(x_+)-\ul{b}J^\delta(x)\leq G^\delta(x)(x_+-x)+\tfrac12 H^\delta(x)(x_+-x)^2
\leq \ol{a}J^\delta(x_+)-\ol{b}J^\delta(x)\\
\mbox{ for all } x,x_+\in \tilde{M}^\delta \,, \quad 
\delta\in(0,\bar{\delta})\,,
\end{aligned}
\end{equation}
with $\ul{a},\ul{b},\ol{a},\ol{b}\geq0$. 
Note that \eqref{abc2} implies $(\ol{a}-\ul{a})J^\delta(x_+)+(\ul{b}-\ol{b})J^\delta(x)\geq0$, hence by nonnegativity of $J^\delta$ and the fact that $J^\delta(x^\dagger)\leq\eta$ can get arbitrarily close to zero, $\ol{a}\geq\ul{a}$ and $\ul{b}\geq\ol{b}$.
In fact, \eqref{abc2} implies \eqref{abc1} with $a=\ul{a}$, $b=\ul{b}-\ol{b}$, $c=\ol{a}$.

\begin{theorem}\label{th:aposteriori}
Let conditions \eqref{etadelta}, \eqref{normalize_delta}, \eqref{GHR}, \eqref{abc2}, and Assumptions \ref{ass0}, \ref{ass1}, \ref{ass2} hold, assume that $\alpha_k$ is chosen a posteriori according to \eqref{alpha_aposteriori} if \eqref{casea} holds (otherwise set $x_{k+1}:=x^*$), and $k_*$ is chosen according to the discrepancy principle \eqref{discrprinc}, with the following constraints on the constants
\begin{equation}\label{const_thaposteriori}
1+\frac{\bar{a}}{\tau}<\ul{\sigma}+\ol{b}\,, \quad
\ol{\sigma}+\ul{b}<1+\ul{a}
\,.
\end{equation}
Then 
\begin{itemize}
\item For any $\delta\in(0,\bar{\delta})$, and any $x_0\in \bigcap_{\delta\in(0,\bar{\delta})}\cap M$, 
\begin{itemize}
\item the iterates $x_k$ are well-defined for all $k\leq k_*(\delta)$ and $k_*(\delta)$ is finite;
\item for all $k\in\{1,\ldots, k_*(\delta)\}$ and $q=\frac{\ol{\sigma} -1+\ul{b}}{\ul{a}}<1$ we have
\[ 
J^\delta(x_k)\leq q J^\delta(x_{k-1});
\]
\item for all $k\in\{1,\ldots, k_*(\delta)\}$ and $x^\dagger$ satisfying \eqref{etadelta} we have
\[
\calR(x_k)\leq \calR(x^\dagger) \mbox{ and }x^\dagger\mbox{ solves \eqref{minJ}.}
\] 
\end{itemize}
\item As $\delta\to0$, the final iterates $x_{k_*(\delta)}$ tend to a solution of the inverse problem \eqref{minJ} $\calT$-subsequentially, i.e., every sequence $x_{k_*(\delta_j)}$ with $\delta_j\to0$ as $j\to\infty$ has a $\calT$ convergent subsequence and the limit of every $\calT$ convergent subsequence solves \eqref{minJ}. 
\end{itemize}
\end{theorem}
Note that the conditions \eqref{const_thaposteriori} on the constants can be satisfied by choosing $\tau$ sufficiently large and $\ul{\sigma}<\ol{\sigma}$ in an appropriate way, provided the constants in \eqref{abc2} satisfy 
\[\ul{b} <\ul{a}+\ol{b}\,,\]
since then we can choose $\ul{\sigma},\ol{\sigma}$ to satisfy $1-\ol{b}<\ul{\sigma}<\ol{\sigma}<1+\ul{a}-\ul{b}$, so that 
\eqref{const_thaposteriori} can be achieved by making $\tau$ large enough.

{\em Proof.}\
Existence of minimizers $x_{k+1}(\alpha)$ of \eqref{Newton} with $\alpha>0$ in place of $\alpha_k$ follows like in the a priori setting of Theorem \ref{th:apriori}, using the fact that \eqref{abc2} implies \eqref{abc1}.

To prove that $\alpha_k$ satisfying \eqref{alpha_aposteriori} exists under condition \eqref{casea}, we first of all verify the upper bound with $\alpha=0$ (which actually does not require \eqref{casea}). To this end, we make use of minimality \eqref{minimality} and the upper bound in \eqref{abc2} to conclude
\[
\begin{aligned}
\sigma_k(\alpha)\leq
\frac{J^\delta(x_k)+G^\delta(x_k)(x^\dagger-x_k)+\tfrac12 H^\delta(x_k)(x^\dagger-x_k)^2 + \alpha(\calR(x^\dagger)-\calR(x_{k+1}(\alpha)))}{J^\delta(x_k)}\\
\leq 1-\ol{b}+\bar{a}\frac{J^\delta(x^\dagger)}{J^\delta(x_k)}+\alpha\frac{\calR(x^\dagger)-\calR(x_{k+1}(\alpha))}{J^\delta(x_k)}\,,
\end{aligned}
\]
so that by \eqref{discrprinc}, for any $k\in\{1,\ldots,k_*-1\}$
\[
\lim_{\alpha\searrow0}\sigma_k(\alpha)\leq 1-\ol{b}+\frac{\bar{a}}{\tau}<\ul{\sigma}\,.
\]
On the other hand, minimality and the fact that $x^*\in\tilde{M}^\delta$ together with the lower bound in \eqref{abc2} and $\calR(x^*)=0$ yield
\[
\begin{aligned}
&\ul{a}J^\delta(x_{k+1}(\alpha))-\ul{b}J^\delta(x_k) + \alpha\calR(x_{k+1}(\alpha))\\
&\leq G^\delta(x_k)(x_{k+1}(\alpha)-x_k)+\tfrac12 H^\delta(x_k)(x_{k+1}(\alpha)-x_k)^2 + \alpha\calR(x_{k+1}(\alpha))\\
&\leq G^\delta(x_k)(x^*-x_k)+\tfrac12 H^\delta(x_k)(x^*-x_k)^2
\end{aligned}
\]
which by nonnegativity of $\ul{a}J^\delta(x_{k+1}(\alpha))$ yields
\[
\calR(x_{k+1}(\alpha))\leq \frac{1}{\alpha}\Bigl(\ul{b}J^\delta(x_k)+G^\delta(x_k)(x^*-x_k)+\tfrac12 H^\delta(x_k)(x^*-x_k)^2\Bigr)\to0\mbox{ as }\alpha\to\infty
\]
which by Assumption \ref{ass2} implies $\calT_1$ convergence of $x_{k+1}(\alpha)$ to $x^*$, thus, by \eqref{casea} $\lim_{\alpha\to\infty} \sigma_k(\alpha)\geq \ul{\sigma}$.
The Intermediate Value Theorem together with continuity of the mapping $\alpha\mapsto \sigma_k(\alpha)$ according to Lemma \ref{lem:cont} implies existence of an $\alpha\in(0,\infty)$ such that $\ul{\sigma}\leq\sigma_k(\alpha)\leq\ol{\sigma}$.

In both cases we get geometric decay of the cost function values: If \eqref{casea} is satisfied, this follows from the lower bound in \eqref{abc2} and the upper  bound in \eqref{alpha_aposteriori}
\[
J^\delta(x_{k+1})\leq\tfrac{1}{\ul{a}}\Bigl( \ul{b}J^\delta(x_k)+G^\delta(x_k)(x_{k+1}-x_k)+\tfrac12 H^\delta(x_k)(x_{k+1}-x_k)^2\Bigr)\leq \frac{\ol{\sigma} -1+\ul{b}}{\ul{a}} J^\delta(x_k)\,.
\]
Otherwise, negation of \eqref{casea} and the fact that in that case we set $x_{k+1}=x^*$, together with the lower bound in \eqref{abc2} directly yields
\[
J^\delta(x_{k+1})=J^\delta(x^*)\leq\tfrac{1}{\ul{a}}\Bigl( \ul{b}J^\delta(x_k)+G^\delta(x_k)(x^*-x_k)+\tfrac12 H^\delta(x_k)(x^*-x_k)^2\Bigr)\leq \frac{\ul{\sigma} -1+\ul{b}}{\ul{a}} J^\delta(x_k)
\]

This implies that $k_*$ is finite, more precisely $k_*\leq \frac{\log(\tau\eta)-\log(J^\delta(x_0))}{\log(q)}$.

To establish the bound on $\calR(x_{k+1})$, we again employ minimality \eqref{minimality} together with \eqref{abc2}, which in case \eqref{casea} with \eqref{alpha_aposteriori} yields
\[
\begin{aligned}
&\ul{\sigma}J^\delta(x_k) + \alpha_k\calR(x_{k+1})\\
&\leq 
J^\delta(x_k)+G^\delta(x_k)(x_{k+1}-x_k)+\tfrac12 H^\delta(x_k)(x_{k+1}-x_k)^2 + \alpha_k\calR(x_{k+1})\\
&\leq 
J^\delta(x_k)+G^\delta(x_k)(x^\dagger-x_k)+\tfrac12 H^\delta(x_k)(x^\dagger-x_k)^2 + \alpha_k\calR(x^\dagger)\\
&\leq 
\ol{a}J^\delta(x^\dagger) +(1-\ol{b})J^\delta(x_k)+\alpha_k\calR(x^\dagger)\,,
\end{aligned}
\]
hence, due to \eqref{discrprinc}, $\tau(\ol{b}+\ul{\sigma}-1)\geq\ol{a}$, 
\[
\calR(x_{k+1})\leq \calR(x^\dagger)+\frac{1}{\alpha_k} \Bigl(\ol{a}J^\delta(x^\dagger)-(\ol{b}+\ul{\sigma}-1) J^\delta(x_k)\Bigr)\leq \calR(x^\dagger)\,.
\]
If \eqref{casea} fails to hold then we set $x_{k+1}=x^*$, hence get $\calR(x_{k+1})=0$.

The rest of the proof is the same as for Theorem \ref{th:apriori}.
\begin{flushright}
$\diamondsuit$
\end{flushright}

\begin{remark} \label{rem:conv-nonl}
Condition \eqref{abc2} is motivated by the fact that 
\[
G^\delta(x)(x_+-x)+\tfrac12 H^\delta(x)(x_+-x)^2\approx J^\delta(x_+)-J^\delta(x)\,,
\]
with equality in case of a quadratic functional $J^\delta$ \eqref{Jquad} 
from which (again using nonnegativity of $J^\delta$) we expect values $\ul{a}\leq1$, $\ul{b}\geq1$, $\ol{a}\geq1$, $\ol{b}\leq1$ where these constants can be chosen the closer to one the closer $J^\delta$ is to a quadratic functional. Also note that \eqref{abc2} holds with $\ul{a}=\ul{b}=\ol{a}=\ol{b}$ in the quadratic case \eqref{Jquad} independently of the definiteness of the Hessian, so does not necessarily relate to convexity of $J^\delta$. Indeed, while nonnegativity of the Hessian would be enforced by assuming $J^\delta\geq0$ on all of $X$, we only assume this to hold on $\tilde{M}^\delta$ cf. \eqref{normalize_delta}.

A sufficient condition for \eqref{abc2} (with $\ul{a}=1-\tilde{c}$, $\ul{b}=1+\tilde{c}$, $\ol{a}=1+\tilde{c}$, $\ol{b}=1-\tilde{c}$) is
\begin{equation}\label{tcc}
\begin{aligned}
\tilde{c}(J^\delta(x_+)+J^\delta(x))
\geq |J^\delta(x_+)-J^\delta(x)-G^\delta(x)(x_+-x)-\tfrac12 H^\delta(x)(x_+-x)^2|
\\
\mbox{ for all } x,x_+\in \tilde{M}^\delta \,, \quad 
\delta\in(0,\bar{\delta})\,,
\end{aligned}
\end{equation}
which, in its turn is implied by the weak tangential cone condition in the Hilbert space least squares setting \eqref{JFxy}
\begin{equation}\label{tccFxy}
\begin{aligned}
|\langle F(x_+)-F(x)-F'(x)(x_+-x),F(x)-y^\delta\rangle| 
\leq c_{tc}\|F(x_+)-F(x)\|\, \|F(x)-y^\delta\|\\
\mbox{ for all } x,x_+\in \tilde{M}^\delta \,, \quad 
\delta\in(0,\bar{\delta})\,,
\end{aligned}
\end{equation}
with $\tilde{c}=(1+\sqrt{2})c_{tc}$; cf. \eqref{tccFxy_grad}. 
This can be seen by using the fact that the left hand side in \eqref{tccFxy} just equals the left hand side in \eqref{tcc} with \eqref{JFxy}, and by estimating the right hand side with $\alpha:=\|F(x_+)-y^\delta\|$, $\beta:=\|F(x)-y^\delta\|$ as follows
\[
\|F(x_+)-F(x)\|\, \|F(x)-y^\delta\|\leq (\alpha+\beta)\beta\leq\frac{1+\sqrt{2}}{2}(\alpha^2+\beta^2)=(1+\sqrt{2})(J^\delta(x_+)+J^\delta(x))\,.
\]
Condition \eqref{tccFxy} with $x_+=x^\dagger$ is also sufficient for condition \eqref{weaktcc_grad} from the previous section with $\gamma= 1-c_{tc}-\kappa$ provided  $(1+c_{tc})\|F(x)-y^\delta\|\leq2\sqrt{\kappa\eta(\delta)}$ and $\|F'(x)\|\leq 1$ as the estimate
\[
\begin{aligned}
&\langle F'(x)(x-x^\dagger),F(x)-y^\delta\rangle 
\geq \langle F(x)-F(x^\dagger),F(x)-y^\delta\rangle - c_{tc}\|F(x)-F(x^\dagger)\|\, \|F(x)-y^\delta\|\\
&= \|F(x)-y^\delta\|^2 - \langle F(x^\dagger)-y^\delta,F(x)-y^\delta\rangle - c_{tc}(\|F(x)-y^\delta-(F(x^\dagger)-y^\delta)\|\, \|F(x)-y^\delta\|\\
&\geq (1-c_{tc}) \|F(x)-y^\delta\|^2 -(1+c_{tc}) \|F(x^\dagger)-y^\delta\|\,\|F(x)-y^\delta\|\\
&\geq (1-c_{tc}-\kappa) \|F(x)-y^\delta\|^2 -\frac{(1+c_{tc})^2}{4\kappa} \|F(x^\dagger)-y^\delta\|^2
\end{aligned}
\]
following from \eqref{tccFxy} with the triangle inequality and Young's inequality shows.

\medskip

On order to further relate the assumptions \eqref{abc1}, \eqref{abc2} made for Newton's method with those \eqref{convex1}, \eqref{convex2} for the projected gradient method, we will now point out that actually also the sufficient condition \eqref{tcc} involves some convexity.

For this purpose we consider the noise free case $\delta=0$ for simplicity of exposition and use the fact that for $n\in\mathbb{N}_0$, a functional $J\in C^n(X)$ and elements $x,\tilde{x},h\in X$ the identity
\[
(J^{(n-1)}(\tilde{x})-J^{(n-1)}(x))[h^{n-1}] = \int_0^1 J^{(n)}[x+\theta(\tilde{x}-x))[\tilde{x}-x,h^{n-1}]\, d\theta
\]
holds. Thus we can rewrite the left hand sides of the nonlinearity conditions \eqref{convex1}, \eqref{tcc} as 
\begin{equation}\label{Jppgrad}
\begin{aligned}
&\langle \nabla J(x)- \nabla J(x^\dagger),x-x^\dagger\rangle
=J'(x)- J'(x^\dagger)[x-x^\dagger]\\
&=\int_0^1 J''(x^\dagger+\theta(x-x^\dagger))[(x-x^\dagger)]^2\, d\theta \,, 
\end{aligned}
\end{equation}
and, with $J(x^\dagger)=0$ 
\[
J(x_+)+J(x) = (J(x_+)-J(x^\dagger))+(J(x)-J(x^\dagger))
\]
where, assuming $J'(x^\dagger)=0$ (as is the case in the examples from section  \ref{sec:appl})
\begin{equation}\label{JppNewton}
\begin{aligned}
&J(x)-J(x^\dagger)=
\int_0^1 \Bigl(J'(x^\dagger+\theta(x-x^\dagger))-J'(x^\dagger)\Bigr)[x-x^\dagger]\, d\theta\\
&=\int_0^1 \int_0^1 \theta J''(x^\dagger+\theta\sigma(x-x^\dagger))[(x-x^\dagger)^2]\,d\sigma\, d\theta
\end{aligned}
\end{equation}
and likewise for $x$ replaced by $x_+$. 
Similarly, using the identities $\int_0^1\, d\theta=1$, $\int_0^1\int_0^1\theta \,d\sigma\, d\theta=\frac12$, 
one sees that for the right hand side in \eqref{tcc} with $G:=J'$, $H:=J''$, the identity 
\[
\begin{aligned}
&J(x_+)-J(x)-J'(x)(x_+-x)-\tfrac12 J''(x)(x_+-x)^2 \\
&= \int_0^1\int_0^1\int_0^1\theta^2\sigma J'''(x^\dagger+\theta\sigma\rho(x-x^\dagger))[(x-x^\dagger)^3]\, d\rho\, d\sigma\, d\theta
\end{aligned}
\]
holds.

Since the left hand sides in both \eqref{convex1}, \eqref{tcc} both have to be nonnegative (in some uniform sense) we see from \eqref{Jppgrad} and \eqref{JppNewton} (and setting $x_+=x^\dagger$ in \eqref{tcc} to see necessity) that $J''$ needs to be positive definite (in some uniform sense) in order for \eqref{convex1}, \eqref{tcc} to hold.
This amounts to a convexity condition on $J$.
\end{remark}

\begin{remark}
Alternatively to \eqref{Newton} one could consider the projected versions (based on unconstrained minimization)
\begin{equation}\label{Newton1}
\tilde{x}_{k+1}\in \mbox{argmin}_{x\in X} J^\delta(x_k)+G(x_k)(x-x_k)+\tfrac12 H(x_k)(x-x_k)^2 + \alpha_k\calR(x)
\quad
x_{k+1}=\mbox{Proj}_{\tilde{M}^\delta}(\tilde{x}_{k+1})
\end{equation}
see \cite{muge12}
which, however, analogously to the projected Landweber iteration from \cite[Section 3.2]{Eicke92} only converges under a sufficiently strong source condition. 
\end{remark}

\section{Application in diffusion/impedance identification}\label{sec:appl}
\def\EE{\mathbf{E}}
\def\JJ{\mathbf{J}}
\def\vv{\mathbf{v}}
\def\ww{\mathbf{w}}
Following the seminal idea from \cite{KohnVogelius87} we consider variational formulations of the problem of identifying the spatially varying parameter $\sigma$ in the elliptic PDE
\begin{equation}\label{ellPDE}
\nabla\cdot(\sigma\nabla\phi)=0\mbox{ in }\Omega 
\end{equation}
from observations of $\phi$.
Depending on what kind of observations we consider, this problem arises in several applications that we will consider here, namely
\begin{itemize}
\item[(a)]
in classical electrical impedance tomography EIT, where it is known as Calderon's problem and $\sigma$ plays the role of an electrical conductivity, 
\item[(b)]
in impedance acoustic tomography IAT, a novel hybrid imaging method, again for reconstructing $\sigma$ as a conductivity;
\item[(c)] 
but also as a simplified version of the inverse groundwater filtration problem GWF of recovering the diffusion coefficient $\sigma$ in an aquifer. 
\end{itemize}
Although we will finally be only able to verify the crucial conditions \eqref{convex1}, \eqref{tcc} for GWF, we stick to the electromagnetic context notation wise, since in our numerical experiments we will focus on a version of EIT that is known as impedance acoustic tomography IAT, see, e.g., \cite{WidlakScherzer2012}.
In Section \ref{sec:numexp_EIT-IAT} we will also allow for experiments with several excitations (and corresponding measurements), hence consider
\[
\nabla\cdot(\sigma\nabla\phi_i)=0\mbox{ in }\Omega \,, \quad i\in\{1,\ldots,I\} \,.
\]
However for simplicity of notation, we will focus on the case $I=1$, i.e., \eqref{ellPDE}, in the current section.
The observations are, depending on the application 
\[
\begin{aligned}
&v = \phi\vert_{\partial\Omega} \mbox{ (the voltage at the boundary in EIT),}\\
&\mathcal{H}=\sigma |\nabla \phi|^2 \mbox{ (the power density in IET),} \\
&p= \phi\mbox{ or } g = \nabla\phi \mbox{ (the piezometric head or its gradient in GWF),} 
\end{aligned}
\]
where for EIT and IAT we will consider the more realistic complete electrode model in Section \ref{sec:numexp_EIT-IAT}.
Concerning GWF, measurements are actually done on the piezometric head itself, however this allows to recover an approximation of its its gradient by means of regularized numerical differentiation, see, e.g. \cite{HankeScherzer2001} and the refernces therein.

Considering a smooth and simply connected bounded domain $\Omega\subseteq\mathbb{R}^2$ and using the vector fields $\EE$ (the electric field), $\JJ$ (the current density), where $\nabla =\left(\begin{array}{c}\partial_1\\ \partial_2\end{array}\right)$, $\nabla^\bot=\left(\begin{array}{c}-\partial_2\\ \partial_1\end{array}\right)$ we can equivalently rephrase \eqref{ellPDE} as
\[
\sigma \EE = \JJ\,, \quad \EE=\nabla \phi\,, \quad \JJ=\nabla^\bot\psi\,,
\]
for some potential $\psi$ (note that we are using the opposite sign convention as compared to the usual engineering notation).
The cost function part pertaining to this model is, analogously to \cite{KohnVogelius87}, therefore often called the Kohn-Vogelius functional
\begin{equation}\label{Jmod_KV}
J_{mod}^{KV}(\sigma,\EE,\JJ)=\tfrac12\int_\Omega \left|\sqrt{\sigma}\EE-\tfrac{1}{\sqrt{\sigma}}\JJ\right|^2\,d\Omega,
\end{equation}
where we denote the infinitesimal area element by $d\Omega$ to avoid confusion with the abbreviation $x_k$ for the iterates in the first three sections of this paper.
Alternatively, we will consider the output least squares type cost function term  
\begin{equation}\label{Jmod_LS}
J_{mod}^{LS}(\sigma,\EE,\JJ)=\tfrac12\int_\Omega \left|\sigma\EE-\JJ\right|^2\,d\Omega\,.
\end{equation}
Note that \eqref{Jmod_LS} is quadratic with respect to $J$, thus quadratic with respect to $\psi$.

Excitation is imposed via the current $j$ through the boundary, i.e., as Dirichlet boundary condition on $\psi$.

To incorporate the observations, we will consider the functionals
\begin{equation}\label{Jobs_EIT_IAT_GWF}
\begin{aligned}
&J_{obs}^{EIT}(\phi;v) = \tfrac12\int_{\partial\Omega} (\phi-v)^2\, d\Omega \mbox{ for EIT,}\\
&J_{obs_1}^{IAT}(\EE,\JJ;\mathcal{H}) = \tfrac12\int_\Omega (\JJ\cdot\EE-\mathcal{H})^2\, d\Omega \mbox{ or }
J_{obs_2}^{IAT}(\sigma,\EE;\mathcal{H}) = \tfrac12\int_\Omega (\sigma|\EE|^2-\mathcal{H})^2\, d\Omega \mbox{ for IAT,}\\
&J_{obs_1}^{GWF}(\phi;p) = \tfrac12\|\phi-p\|_{H^s(\Omega)}^2\mbox{ or }J_{obs_2}^{GWF}(\EE;g) = \|\EE-g\|_{L^2(\Omega)}^2 \mbox{ for GWF,}
\end{aligned}
\end{equation}
where again for GWF the use of the $H^s(\Omega)$ norm or flux data can be justified by some pre-smoothing procedure applied to the given measurements.

\medskip

Using these functionals as building blocks and incorporating the excitation via injection of the current $j$ through the boundary we can write the above parameter identification problems in several minimization based formulations. We will now list a few of them, where $j$ sometimes appears explicitely, sometimes in tangentially integrated form, meaning that for a parametrization $\Gamma$ of the boundary $\partial\Omega$ (normalized to $\|\dot{\Gamma}\|=1$) we define $\alpha(\Gamma(s))=\int_0^s j(\Gamma(r))\, d r$ so that $\mathbf{J}\cdot\nu=\nabla^\bot\psi\cdot \nu=\frac{d\alpha}{ds}=j$. 
Moreover we will sometimes work with smooth extensions $\phi_0$, $\psi_0$ of $v$, $\alpha$ to the interior of $\Omega$. 
While, as already mentioned, the observation functional will depend on the application, we always have both $J_{mod}^{KV}$ and $J_{mod}^{LS}$ at our disposal to incorporate the model, thus will only write $J_{mod}$ below. 
There will also be versions based on an elimination of $\sigma$ by writing, for fixed $\phi,\psi$, the minimizer of $J_{mod}$ with respect to $\sigma$ under the constraint $\underline{\sigma}\leq\sigma\leq\overline{\sigma}$ as 
\[
\sigma(\EE,\JJ)=\max\{\underline{\sigma},\min\{\overline{\sigma},\tfrac{|\JJ|}{|\EE|}\}\} \mbox{ pointwise in }\Omega\,.
\]
Alternatively it is also possible to eliminate $\phi,\psi$ by writing them as $\phi(\sigma)$, $\psi(\sigma)$ mimimizing $J_{mod}$ with respect to $\phi,\psi$. This together with the integrated current $\alpha$ leads to boundary value problems for the elliptic PDE \eqref{ellPDE} and a similar PDE for $\psi$ 
\[
\begin{aligned}
&\phi(\sigma)\mbox{ solves } \left\{\begin{array}{rcll}\nabla\cdot(\sigma\nabla\phi)&=&0&\mbox{ in }\Omega\\
\phi&=&v&\mbox{ on }\partial\Omega\end{array}\right. \\
&\phi_N(\sigma)\mbox{ solves } \left\{\begin{array}{rcll}\nabla\cdot(\sigma\nabla\phi)&=&0&\mbox{ in }\Omega\\
\nabla\phi\cdot\nu&=&j&\mbox{ on }\partial\Omega\, \quad\int_\Omega \phi\, d\Omega =0\end{array}\right. \\
&\psi(\sigma)\mbox{ solves } \left\{\begin{array}{rcll}\nabla^\bot\cdot(\frac{1}{\sigma}\nabla^\bot\psi)&=&0&\mbox{ in }\Omega\\
\psi&=&\alpha&\mbox{ on }\partial\Omega\end{array}\right. \\ 
&\EE(\sigma) = \tfrac{\nabla^\bot\psi(\sigma)}{\sigma}\mbox{ pointwise in }\Omega
\end{aligned}
\]
(the latter two lines imply that $\nabla^\bot\cdot\EE(\sigma)=0$ so existence of $\phi$ such that $\EE(\sigma)=\nabla\phi$) 
and corresponds to the classical reduced formulation of the inverse problem.
Note that $\phi(\sigma)$ is only defined in case of $v$ being observed, i.e., for EIT.

\paragraph{EIT:}
\[
\begin{aligned}
(i)&\min_{\sigma,\phi,\psi}\{J_{mod}(\sigma,\nabla\phi,\nabla^\bot\psi)+\beta J_{obs}^{EIT}(\phi;v) \, : \, \sigma\in L^2_{[\underline{\sigma},\overline{\sigma}]}(\Omega)\,,\ \phi\in H_\diamondsuit^1(\Omega), \psi\in H_0^1(\Omega)+\psi_0\}\\
(ii)&\min_{\sigma,\phi,\psi}\{J_{mod}(\sigma,\nabla\phi,\nabla^\bot\psi)\, : \, \sigma\in L^2_{[\underline{\sigma},\overline{\sigma}]}(\Omega)\,,\  \phi\in H_0^1(\Omega)+\phi_0\,,\ \psi\in H_0^1(\Omega)+\psi_0\}\\
(iii)&\min_{\phi,\psi}\{J_{mod}(\sigma(\nabla\phi,\nabla^\bot\psi),\nabla\phi,\nabla^\bot\psi)+\beta J_{obs}^{EIT}(\phi;v)\, : \, \phi\in H_\diamondsuit^1(\Omega), \psi\in H_0^1(\Omega)+\psi_0\}\\
(iv)&\min_{\phi,\psi}\{J_{mod}(\sigma(\nabla\phi,\nabla^\bot\psi),\nabla\phi,\nabla^\bot\psi)\, : \, \phi\in H_0^1(\Omega)+\phi_0, \psi\in H_0^1(\Omega)+\psi_0\}\\
(v)&\min_{\sigma}\{J_{mod}(\sigma,\nabla\phi(\sigma),\nabla^\bot\psi(\sigma))\, : \, \sigma \in L^2_{[\underline{\sigma},\overline{\sigma}]}(\Omega)\}\\
(vi)&\min_{\sigma}\{J_{obs}^{EIT}(\phi_N(\sigma);v)\, : \, \sigma \in L^2_{[\underline{\sigma},\overline{\sigma}]}(\Omega)\}
\end{aligned}
\]
\paragraph{IAT:}
\[
\begin{aligned}
(i)&\min_{\sigma,\phi,\psi}\{J_{mod}(\sigma,\nabla\phi,\nabla^\bot\psi)
+\beta \left\{\begin{array}{l}J_{obs_1}^{IAT}(\nabla\phi,\nabla^\bot\psi;\mathcal{H})\\J_{obs_2}^{IAT}(\sigma,\nabla\phi;\mathcal{H})\end{array}\right. 
\, : \, \sigma\in L^2_{[\underline{\sigma},\overline{\sigma}]}(\Omega)\,, \ \phi\in H_\diamondsuit^1(\Omega), \psi\in H_0^1(\Omega)+\psi_0\}\\
(ii)&\min_{\phi,\psi}\{J_{mod}(\sigma(\nabla\phi,\nabla^\bot\psi),\nabla\phi,\nabla^\bot\psi)
+\beta \left\{\begin{array}{l}J_{obs_1}^{IAT}(\nabla\phi,\nabla^\bot\psi;\mathcal{H})\\J_{obs_2}^{IAT}(\sigma,\nabla\phi;\mathcal{H})\end{array}\right. 
\, : \, \phi\in H_\diamondsuit^1(\Omega), \psi\in H_0^1(\Omega)+\psi_0\}\\
(iii)&\min_{\sigma}\{\left\{\begin{array}{l}J_{obs_1}^{IAT}(\EE(\sigma),\nabla^\bot\psi(\sigma);\mathcal{H})\\
J_{obs_2}^{IAT}(\sigma,\EE(\sigma);\mathcal{H})\end{array}\right.\, : \, \sigma \in L^2_{[\underline{\sigma},\overline{\sigma}]}(\Omega)\}
\end{aligned}
\]
\paragraph{GWF:}
\[
\begin{aligned}
(i)&\min_{\sigma,\phi,\psi}\{J_{mod}(\sigma,\nabla\phi,\nabla^\bot\psi)
+\beta \left\{\begin{array}{l}J_{obs_1}^{GWF}(\phi;p)\\J_{obs_2}^{GWF}(\nabla\phi;g)\end{array}\right. 
\, : \, \sigma\in L^2_{[\underline{\sigma},\overline{\sigma}]}(\Omega)\,, \ \phi\in H_\diamondsuit^1(\Omega), \psi\in H_0^1(\Omega)+\psi_0\}\\
(ii)&\min_{\phi,\psi}\{J_{mod}(\sigma(\nabla\phi,\nabla^\bot\psi),\nabla\phi,\nabla^\bot\psi)
+\beta \left\{\begin{array}{l}J_{obs_1}^{GWF}(\phi;p)\\J_{obs_2}^{GWF}(\nabla\phi;g)\end{array}\right. 
\, : \, \phi\in H_\diamondsuit^1(\Omega), \psi\in H_0^1(\Omega)+\psi_0\}\\
(iii)&\min_{\sigma}\{J_{obs_2}^{GWF}(\EE(\sigma);g)\, : \, \sigma \in L^2_{[\underline{\sigma},\overline{\sigma}]}(\Omega)\}
\end{aligned}
\]
where
\[
L^2_{[\underline{\sigma},\overline{\sigma}]}(\Omega)=\{\sigma\in L^2(\Omega)\, : \, \underline{\sigma}\leq\sigma\leq\overline{\sigma}\}
\,, \quad 
H_\diamondsuit^1(\Omega) =\{\phi\in H^1(\Omega)\, : \, \int_\Omega \phi\, d\Omega =0\}\,,
\]
and $\beta>0$ is a fixed parameter; we will simply set it to one in our computations.
Note that $J_{mod}(\sigma,\EE(\sigma),\nabla^\bot\psi(\sigma))=0$, therefore, the model term does not appear in the last instances of IAT and GWF, respectively. However, due to the bound constraints incorporated into the definition of $\sigma(\phi,\psi)$, a nonzero value of $J_{mod}(\sigma(\phi,\psi),\nabla\phi,\nabla^\bot\psi)$ is possible, which is why it appears in the third and fourth instances of EIT.
The sixth instance of EIT is just the classical reduced formulation. 

\medskip

As far as convexity is concerned, the Hessians of the functionals in \eqref{Jmod_KV}, \eqref{Jmod_LS}, \eqref{Jobs_EIT_IAT_GWF} compute as 
\[
\begin{aligned}
&J_{mod}^{KV''}(\sigma,\EE,\JJ)[(h,\vv,\ww)^2]=\int_\Omega \Bigl\{\left|\tfrac{h}{\sqrt{\sigma}^3}\JJ+\sqrt{\sigma}\vv-\tfrac{1}{\sqrt{\sigma}}\ww\right|^2+2\tfrac{h}{\sigma}(\sqrt{\sigma}\EE-\tfrac{1}{\sqrt{\sigma}}\JJ)\cdot\vv\Bigr\}\,d\Omega\\
&J_{mod}^{LS''}(\sigma,\EE,\JJ)[(h,\vv,\ww)^2]=\int_\Omega \Bigl\{\left|h\EE+\sigma\vv-\ww\right|^2+h(\sigma\EE-\JJ)\cdot\vv\Bigr\}\,d\Omega\\
&J_{obs}^{EIT''}(\phi;v)[u^2] = \int_{\partial\Omega} u^2\, d\Omega \\
&J_{obs_1}^{IAT''}(\EE,\JJ;\mathcal{H})[(\vv,\ww)^2] = \int_\Omega |\ww\cdot\EE+\JJ\cdot\vv|^2 + (\JJ\cdot\EE-\mathcal{H})\, \ww\cdot\vv\, d\Omega \\
&J_{obs_2}^{IAT''}(\sigma,\EE;\mathcal{H})[(h,\vv)^2] = \int_\Omega (|h|\EE+2\sigma\EE\cdot\vv)^2+(\sigma|\EE|^2-\mathcal{H})2(h\vv\cdot\EE+\sigma|v|^2)\, d\Omega \\
&J_{obs_1}^{GWF''}(\phi;p)[u^2] = \|u\|_{H^s(\Omega)}^2\\
&J_{obs_2}^{GWF''}(\EE;g)[\vv^2] = \|\vv\|_{L^2(\Omega)}^2\,.
\end{aligned}
\]
Thus, the Hessians of $J_{mod}^{KV}$, $J_{mod}^{LS}$, $J_{obs_1}^{IAT}$, $J_{obs_2}^{IAT}$ can only be guaranteed to be positive at their minimal points, whereas those of $J_{obs}^{EIT}$, $J_{obs_1}^{GWF}$, $J_{obs_2}^{GWF}$ are always positive.
Since $J_{obs}^{EIT}$ only acts on the boundary, its additive combination with $J_{mod}^{KV}$ or $J_{mod}^{LS}$ cannot be expected to yield a globally convex functional. Likewise, combinations of $J_{obs_1}^{IAT}$ or $J_{obs_2}^{IAT}$  with $J_{mod}^{KV}$ or $J_{mod}^{LS}$ cannot be expected to be overall convex. This corresponds to the known fact that also for other formulations of EIT and IAT, the usual nonlinearity/convexity conditions fail to hold.

A combination satisfying the nonlinearity assumption \eqref{tcc} and therefore also \eqref{abc2}, \eqref{convex1} is GWF with 
\[
J^\delta(\sigma,\phi,\psi)=J_{mod}^{LS}(\sigma,\nabla\phi,\nabla^\bot\psi)+\beta J_{obs}^{GWF}(\nabla\phi;g^\delta)
\] 
To verify this, we show that \eqref{weaktcc_grad}, \eqref{tccFxy} is satisfied for 
$F(\sigma,\phi,\psi)=\left(\begin{array}{c}\sigma\nabla\phi-\nabla^\bot\psi\\ \nabla\phi\end{array}\right)$ by estimating (with the abbreviations $\EE=\nabla\psi$, $\JJ=\nabla^\bot\psi$)
\[
\begin{aligned}
&|\langle F(x_+)-F(x)-F'(x)(x_+-x),F(x)-y^\delta\rangle| \\
&=\int_\Omega (\sigma_+-\sigma)(\EE_+-\EE)(\sigma\EE-\JJ)\, d\Omega\\
&\leq \|\sigma_+-\sigma\|_{L^\infty(\Omega)} \|\EE_+-\EE\|_{L^2(\Omega)} \|\sigma\EE-\JJ\|_{L^2(\Omega)}\\
&\leq (\overline{\sigma}-\underline{\sigma}) \sqrt{\|\sigma_+\EE_+-J_+-\sigma\EE+J\|_{L^2(\Omega)}^2+\|\EE_+-\EE\|_{L^2(\Omega)}^2} \\
&\qquad\qquad\ \cdot\sqrt{\|\sigma\EE-\JJ\|_{L^2(\Omega)}^2+\|\EE-g^\delta\|_{L^2(\Omega)}^2}\\
&= (\overline{\sigma}-\underline{\sigma}) \|F(x_+)-F(x)\|\, \|F(x)-y^\delta\|\\
\end{aligned}
\]
which directly implies \eqref{tccFxy} with $c_{tc}=\frac{\overline{\sigma}-\underline{\sigma}}{\sup_{x\in \tilde{M}^\delta}\|F'(x)\|}$ and hence \eqref{weaktcc_grad} with $\gamma= 1-c_{tc}-\kappa$ provided  $(1+c_{tc})\|F(x)-y^\delta\|\leq2\sqrt{\kappa\eta(\delta)}$. In order to obtain a finite value of 
\[
\sup_{x\in \tilde{M}^\delta}\|F'(x)\|=\sup_{(\sigma,\phi,\psi)\in\tilde{M}^\delta}\sup_{(h,\vv,\ww)\in L^2(\Omega)^5\setminus\{0\}}\frac{\int_\Omega (h\nabla\phi+\sigma\vv-\ww)\cdot(\sigma\nabla\phi-\nabla^\bot\psi)\, d\Omega}{\|(h,\vv,\ww)\|_{L^2(\Omega)^5}}
\]
we choose $\tilde{M}^\delta$ to be a bounded subset of $L^\infty(\Omega)\times W^{1,\infty}(\Omega)\times H^1(\Omega)$ with an apriori bound satisfied by the exact solution of the inverse problem.


\section{Numerical results for IAT and EIT}\label{sec:numexp_EIT-IAT}
In this section, we will provide some numerical results for the problem of identifying the conductivity $\sigma$ in \eqref{ellPDE}. As already mentioned, we will work with the more realistic complete electrode model (CEM) instead of idealized continuous boundary excitation and observations. Moreover, we will focus on the hybrid tomographic application IAT and we will only show one set of reconstructions of EIT. More extensive numerical tests for IAT but also for GWF and EIT can be found in the PhD thesis \cite{ThesisKha2021}.

\subsection{The complete electrode model and setting for the cost functions}
\begin{figure}[!h]
\begin{minipage}{0.45\textwidth}
\begin{center}
\begin{tikzpicture} 
\draw (0,0) circle (3 and 2);
\coordinate (P) at ($(0, 0) + (30:3cm and 2cm)$);
\draw[black, line width=1pt] (P) arc (30:70:3cm and 2cm) node[pos=0.4, above]{$g_1$};
\coordinate (P) at ($(0, 0) + (110:3cm and 2cm)$);
\draw[black, line width=1pt] (P) arc (110:150:3cm and 2cm) node[pos=0.6, above]{$g_2$};
\coordinate (P) at ($(0, 0) + (200:3cm and 2cm)$);
\draw[black, line width=1pt] (P) arc (200:270:3cm and 2cm) node[pos=0.5, below]{$g_3$};
\coordinate (P) at ($(0, 0) + (305:3cm and 2cm)$);
\draw[black, line width=1pt] (P) arc (305:360:3cm and 2cm) node[pos=0.5, right]{$g_4$};
\draw[red, line width=2pt] 
(3,0) circle(1pt) node[blue, right] {$e_1^a$} arc (0:30:3cm and 2cm) circle(1pt) node[blue, right] {$e_1^b$} node[pos=0.5, right]{$e_1$};
\coordinate (P) at ($(0, 0) + (70:3cm and 2cm)$);
\draw[red, line width=2pt]
(P) circle (1pt) node[above, blue] {$e_2^a$} arc (70:110:3cm and 2cm) circle (1pt) node[above, blue] {$e_2^b$} node[pos=0.5, above]{$e_2$};
\coordinate (P) at ($(0, 0) + (150:3cm and 2cm)$);
\draw[red, line width=2pt]
(P) circle (1pt) node[left, blue] {$e_3^a$} arc (150:200:3cm and 2cm) circle (1pt) node[left, blue] {$e_3^b$} node[pos=0.5, left]{$e_3$};
\coordinate (P) at ($(0, 0) + (270:3cm and 2cm)$);
\draw[red, line width=2pt]
(P) circle (1pt) node[below, blue] {$e_4^a$} arc (270:305:3cm and 2cm) circle (1pt) node[below, blue] {$e_4^b$} node[pos=0.5, below]{$e_4$};
\end{tikzpicture}
\end{center}
\end{minipage}
\begin{minipage}{0.52\textwidth}
\includegraphics[width=\textwidth]{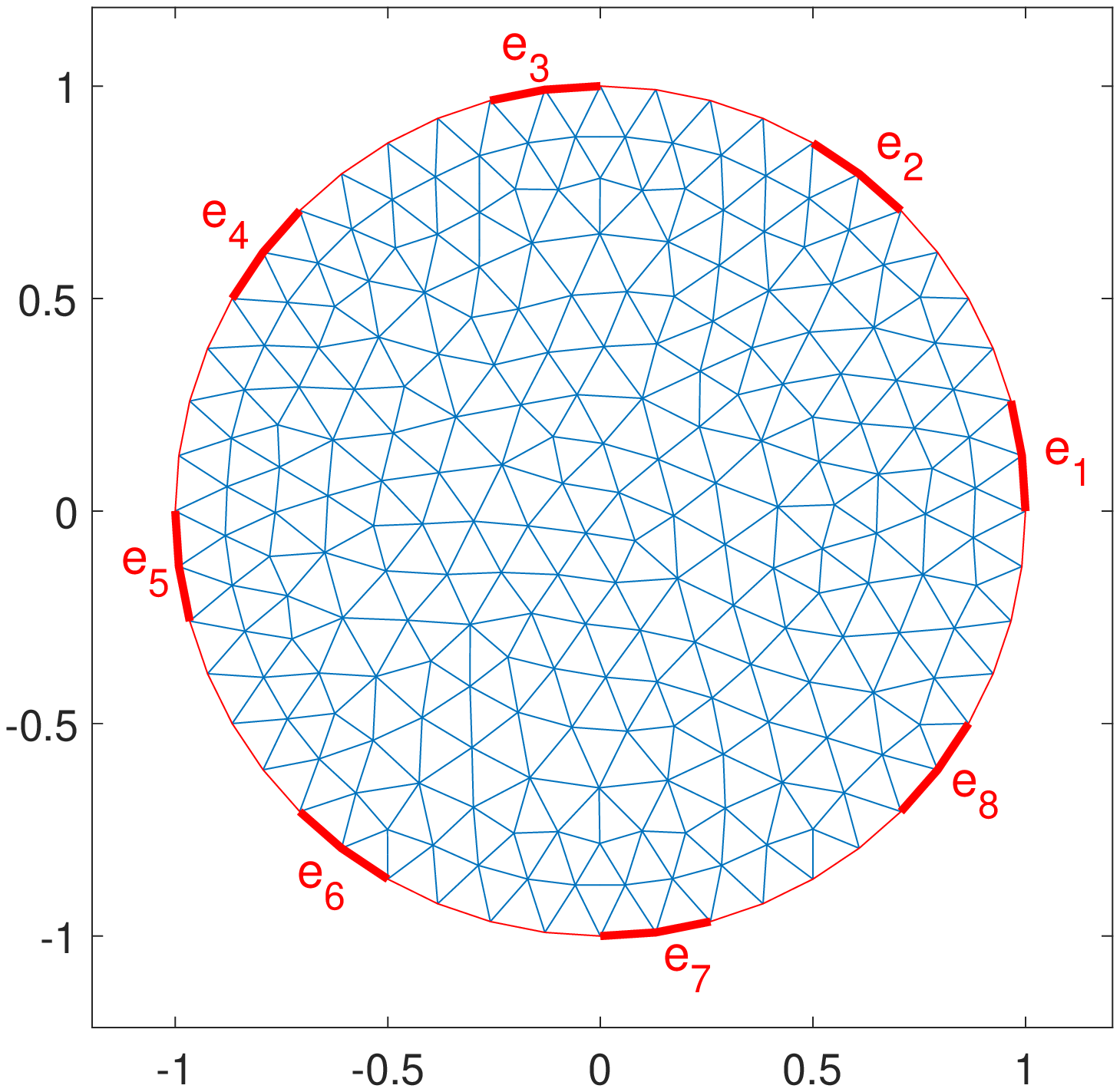} 
\end{minipage}
\caption{left: Electrodes (in red) on the boundary with $L=4$; right: finite element discretization with $8$ electrodes
\label{pic:CEM}}
\end{figure}

In the complete electrode model (CEM) current is fed in through a finite number of electrodes, $e_1, \dots, e_L$, see Figure \ref{pic:CEM}. In case of boundary measurements, as relevant for EIT, they are also taken at these electrodes. 
 
Let $J_i$, $E_i$, $i = 1,2, \dots, I$, be the current density and the electric field in the $i$th measurement, and let $\phi_i$, $\psi_i$ be the potentials for $J_i$, $E_i$; then they must satisfy
\begin{subequations} \label{eq:EIT_CEM}
\begin{alignat}{2}
& \sqrt{\sigma}\nabla \phi_i - \frac{1}{\sqrt{\sigma}} \nabla^\bot \psi_i = 0 &&\quad \text{in \;} \Omega, \label{eq:EIT_CEM_eq} \\
& \phi_i + z_\ell \nabla^\bot \psi_i \cdot \nu = v_{\ell,i} &&\quad \text{on \;} e_\ell, \ell=1,2,\dots,L, \label{eq:EIT_CEM_constrain_el1} \\
& \int_{e_\ell} \nabla^\bot \psi_i \cdot \nu \dd s = j_{\ell,i} &&\quad \text{for \;} \ell=1,2,\dots,L, \label{eq:EIT_CEM_constrain_el2} \\
& \nabla^\bot \psi_i \cdot \nu = 0 &&\quad \text{on \;} \partial \Omega \backslash \cup_{\ell=1}^L e_\ell, \label{eq:EIT_CEM_constrain_gl}
\end{alignat}
\end{subequations}
$\forall i= 1, 2, \dots, I$, where 
\begin{itemize}
\item[] $j_{\ell,i}$, $v_{\ell,i}$ are the applied current and measured voltage on $e_\ell$ at the $i$th measurement,
\item[] $\{z_\ell\}_{\ell=1}^{L}$ is the set of positive contact impedances. 
\end{itemize}
By assuming $\psi(e_i^a) = 0$ and using \eqref{eq:EIT_CEM_constrain_el2}, \eqref{eq:EIT_CEM_constrain_gl} \eqref{eq:EIT_CEM_constrain_el1}, we get
\[
\begin{aligned}
& \psi_i|_{g_{\ell}} = \bar j_{\ell,i}, && \quad \forall \ell \in \{1,\dots,L\}, \\
& \int_{e_\ell^a}^x \phi_i \dd s - z_\ell \psi_i(x)  = \bar v_{\ell,i}(x), && \quad \forall x \in e_\ell, \forall \ell \in \{1,\dots,L\},
\end{aligned}
\]
where $\bar j_{\ell,i} = -\sum_{k=1}^{\ell} j_{k,i}$, $\bar v_{\ell,i}(x) = -z_\ell \Big( -\sum_{k=1}^{\ell-1} j_{k,i} \Big)+ v_{\ell,i} d_{e_\ell} (x)$ and $d_{e_\ell} (x)$ is the length of $e_\ell$ from $e_\ell^a$ to $x$.

In the case of EIT, the data $(j_{\ell,i}, v_{\ell,i})_{\ell,i} \in \mathbb{R}^{2LI}$ can be considered as $(\bar j_{\ell,i}, \bar v_{\ell,i})_{\ell,i} \in \prod_{i=1}^I \prod_{\ell=1}^L (L^2(g_\ell) \times L^2(e_\ell))$ and the cost function part corresponding to observations is chosen as $J^{EIT}_{obs}: H^1(\Omega)^{2I} \to \mathbb{R}$,
\begin{equation} \label{eq:costfunJobs_EIT}
J^{EIT}_{obs}(\Phi, \Psi; \bar j, \bar v) = \frac{1}{2} \sum_{i=1}^I \sum_{\ell=1}^L \Big( \int_{g_\ell} |\psi_i - \bar j_{\ell,i}|^2 \dd s + \int_{e_\ell} \Big|\int \phi_i|_{e_\ell} - z_{\ell} \psi_i - \bar v_{\ell,i} \Big|^2 \dd s \Big),
\end{equation}
where $\Phi = (\phi_i)_i$, $\Psi = (\psi_i)_i$, $\bar j = (\bar j_{\ell,i})_{\ell,i}$, $\bar v = (\bar v_{\ell,i})_{\ell,i}$ and $\int \phi_i|_{e_\ell}: e_\ell \to \mathbb{R}, x \mapsto \int_{e_\ell^a}^x \phi_i \dd s$. 

In the case of IAT, instead of $(j_{\ell,i}, v_{\ell,i})_{\ell,i}$, we observe $\mathcal{H} = (\mathcal{H}_i)_i = (\sigma |\nabla \phi_i|^2)_i$ and the cost function part corresponding to these observations is $J^{IAT}_{obs}: L^2(\Omega) \times H^1(\Omega)^I \to \mathbb{R}$,
\begin{equation} \label{eq:costfunJobs_IAT}
J^{IAT}_{obs}(\sigma, \Phi; \mathcal{H}) = \frac{1}{2} \sum_{i=1}^I \int_\Omega \Big| \sigma |\nabla \phi_i|^2 - \mathcal{H}_i \Big|^2 \dd \Omega.
\end{equation}

In both cases of EIT and IAT, we choose the cost function part corresponding to the model as $J^{KV}_{mod}: L^2(\Omega) \times H^1(\Omega)^{2I} \to \mathbb{R}$,
\begin{equation} \label{eq:costfunJmod_KV}
J^{KV}_{mod}(\sigma, \Phi, \Psi) = \frac{1}{2} \sum_{i=1}^I \int_\Omega \Big| \sqrt{\sigma} \nabla \phi_i - \frac{1}{\sqrt{\sigma}} \nabla^\bot \psi_i \Big|^2 \dd \Omega
\end{equation}
and combine it with the model part to 
\begin{equation} \label{eq:costfunJ_IAT_aao}
J^{IAT}(\sigma, \Phi, \Psi) =  J^{KV}_{mod}(\sigma, \Phi, \Psi) + J^{IAT}_{obs}(\sigma, \Phi; \mathcal{H})
\end{equation}
\begin{equation} \label{eq:costfunJ_EIT_aao}
J^{EIT}(\sigma, \Phi, \Psi) =  J^{KV}_{mod}(\sigma, \Phi, \Psi) + J^{EIT}_{obs}(\Phi, \Psi; \bar j, \bar v)
\end{equation}
on the admissible sets
\[
M^{IAT}_{\ad} = M^{EIT}_{\ad} = L^2_{[\underline{\sigma}, \overline{\sigma}]} (\Omega) \times H_\diamondsuit^1(\Omega) \times H_0^1(\Omega)+\Phi_0
\]

As in the previous section, besides the resulting all-at-once versions (cf. EIT (i), (ii) and IAT (i)) we also consider some of the reduced versions of the cost function. 

The first version involves eliminating $\sigma$ from the cost function (cf. EIT (iii), (iv) and IAT (ii)), by defining, for given $\Phi,\Psi$, the corresponding $\sigma$ by
\begin{equation} \label{eq:sigma}
\begin{split}
\sigma(\Phi,\Psi) &= \argmin_{\sigma} \Big\{ \frac{1}{2} \sum_{i=1}^{I} \int_\Omega \Big| \sqrt{\sigma} \nabla \phi_i - \frac{1}{\sqrt{\sigma}} \nabla^\bot \psi_i \Big|^2 \dd \Omega: \sigma \in L^2_{[\underline{\sigma}, \overline{\sigma}]} (\Omega) \Big\} \\
&= \argmin_{\sigma} \Big\{ \sum_{i=1}^{I} \int_\Omega \Big( \sigma |\nabla \phi_i|^2 + \frac{1}{\sigma} |\nabla^\bot \psi_i|^2 \Big) \dd \Omega: \sigma \in L^2_{[\underline{\sigma}, \overline{\sigma}]} (\Omega) \Big\}
\end{split}
\end{equation}
or explicitly
\begin{equation*}
\begin{split}
\sigma(\Phi,\Psi) &=\min \Bigg\{ \overline{\sigma}, \max\Bigg\{ \underline{\sigma}, \sqrt{\frac{\sum_{i=1}^I |\nabla^\bot \psi_i|^2}{\sum_{i=1}^I |\nabla \phi_i|^2}} \Bigg\} \Bigg\}.
\end{split}
\end{equation*}
For the case of IAT, we set
\begin{equation} \label{eq:costfunJ_IAT_sigma}
\begin{split}
J^{IAT}_\sigma(\Phi, \Psi; \mathcal{H}) & = J^{KV}_{mod} (\sigma(\Phi, \Psi), \Phi, \Psi) + \beta J^{IAT}_{obs} (\sigma(\Phi, \Psi), \Phi; H) 
\end{split}
\end{equation}
for the case of EIT, 
\begin{equation} \label{eq:costfunJ_EIT_sigma}
\begin{split}
J^{EIT}_\sigma(\Phi, \Psi; \bar j, \bar v) & = J^{KV}_{mod} (\sigma(\Phi, \Psi), \Phi, \Psi) + \beta J^{EIT}_{obs} (\Phi, \Psi; \bar j, \bar v) 
\end{split}
\end{equation}
and in both cases
\[
M^{IAT}_{\sigma,\ad} = M^{EIT}_{\sigma,\ad} = H_\diamondsuit^1(\Omega) \times H_0^1(\Omega)+\Phi_0\,.
\]
Note that in spite of the minimizing pre-definition of $\sigma(\Phi, \Psi)$, the model cost function part may be nonzero due to the constraints and therefore still needs to be taken into account.

The second alternative cost function involves eliminating $(\Phi, \Psi)$ from the cost function (cf. EIT (v), (vi) and IAT (iii)) by means of the weak form of the CEM PDE \eqref{eq:EIT_CEM}
\begin{equation} \label{eq:weakform}
\begin{split}
&\int_\Omega \sigma \nabla \phi_i \cdot \nabla p \dd \Omega + \sum_{\ell=1}^L \frac{1}{z_\ell} \int_{e_\ell} (\phi_i - v_{\ell,i}) (p - \xi_{\ell}) \dd s = \sum_{\ell=1}^L j_{\ell,i} \xi_\ell, \forall (p,\xi) \in H^1 (\Omega) \times \mathbb{R}^L,
\end{split}
\end{equation}
and $\nabla^\bot \psi_i = \sigma \nabla \phi_i, \forall i \in \{1,\dots,I\}$, which leads to $J^{KV}_{mod} (\sigma, \Phi(\sigma), \Psi(\sigma)) = 0$. Hence, we have, for the case of IAT,
\begin{equation} \label{eq:costfunJ_IAT_PhiPsi}
\begin{split}
J^{IAT}_{(\Phi,\Psi)}(\sigma; \mathcal{H}) &=  J^{IAT}_{obs} (\sigma, \Phi(\sigma); \mathcal{H}) = \frac{1}{2} \sum_{i=1}^I  \int_\Omega \Big( \sigma |\nabla \phi_i(\sigma)|^2 - \mathcal{H}_i \Big)^2 \dd \Omega;
\end{split}
\end{equation}
for the case of EIT,
\begin{equation} \label{eq:costfunJ_EIT_PhiPsi}
\begin{split}
J^{EIT}_{(\Phi,\Psi)}(\sigma; v) &=  \frac{1}{2} \sum_{i=1}^I \sum_{\ell=1}^L \big| v_{\ell,i}(\sigma) - v_{\ell,i} \big|^2,
\end{split}
\end{equation}
where, $v(\sigma) = (v_{\ell,i}(\sigma))_{i \in \{1,\dots,I\},\ell \in \{1,\dots,L\}}$ is the solution to \eqref{eq:weakform} and
\[
M^{IAT}_{(\Phi, \Psi),\ad} = M^{EIT}_{(\Phi, \Psi),\ad} = L^2_{[\underline{\sigma}, \overline{\sigma}]} (\Omega) 
\]

\subsection{Implementation using the finite method in Matlab}
In order to generate synthetic data by solving the CEM PDE \eqref{eq:EIT_CEM} using the finite element method. In all aour computations, $\Omega$ is the unit circle in $\mathbb{R}^2$ with eight identical electrodes ($L=8$) denoted by $e_1, \dots, e_8$ attached equidistantly on its boundary (see Figure~\ref{pic:CEM}). The domain $\Omega$ is decomposed by a regular finite element mesh defined by nodes $P_k, k\in\{1,\dots, N_\nod = 913\}$ and elements $\Omega_h, h \in \{1,\dots,N_\ele=432\}$ and the ansatz spaces $L^2 (\Omega)$ for $\sigma$ and $H^1 (\Omega)$ for $\phi,\psi$ are approximated by piecewise constant and continuous piecewise quadratic finite elements spaces $\tilde L^2 (\Omega)$ and $\tilde H^1 (\Omega)$, respectvely.

With $L=8$ electrodes, there are $N_\meas = 28$ possible combinations of excitations -- we will use some of them to reconstruct $\sigma$ later. At the $i$th measurement, we impose the injected current $Jsigma_I = (j_{\ell,i})_{\ell=1}^L$ with
 $$\left\{ \begin{array}{ll} j_{\ell_1,i} = 1, j_{\ell_2,i} = -1, &  \text{\, if $\ell_1<\ell_2$ and $\{\ell_1,\ell_2\}$ is the $i$th element of}\\ & \text{\, the family of 2-elements subsets of $\{1,\dots,8\}$}, \\ j_{\ell,i} = 0, & \text{\, otherwise,} \end{array} \right.$$
at the electrodes and then solve the Galerkin discretized weak form (cf. \eqref{eq:weakform})
\begin{equation} \label{eq:weakform_finiteEle}
\begin{split}
\int_\Omega \sigma^\ex \nabla \phi_i \cdot \nabla p \dd x + \sum_{\ell=1}^L \frac{1}{z_\ell} \int_{e_\ell} (\phi_i - v_{\ell,i}) (p - \xi_{\ell}) \dd s = \sum_{\ell=1}^L j_{\ell,i} \xi_\ell, \\
\forall (p,\xi) \in \tilde H^1 (\Omega) \times \mathbb{R}^L, \forall i \in \{1,\dots,N_\meas\}
\end{split}
\end{equation}
to find $(\phi^\ex_i,(v^\ex_{\ell,i}))$ and the corresponding exact data 
\[
\mathcal{H}^\ex = (\sigma^\ex |\nabla \phi^\ex_1|^2, \dots, \sigma^\ex |\nabla \phi^\ex_{N_\meas}|^2)
\mbox{ for IAT;} \quad 
(j_{\ell,i}, v^\ex_{\ell,i})_{\ell,i}
\mbox{ for EIT.}
\]
The synthetic measured data is generated by adding random noise such that
\[
|\mathcal{H}^\delta_i - \mathcal{H}^\ex_i| \le \delta |\mathcal{H}^\ex_i|, \ \forall i
\mbox{ for IAT;} \quad 
|v^\delta_{\ell,i} - v^\ex_{\ell,i}| \le \delta |v^{ex}_{\ell,i}|
\mbox{ for EIT,}
\]
which in an obvious way defines the noisy versions $J^\delta$, $J_\sigma^\delta$, $J_{\Phi,\Psi}^\delta$ and $M^\delta$, $M_{\sigma,\ad}^\delta$, $M_{(\Phi, \Psi),\ad}^\delta$ of the cost functions and admissible sets, respectively.
In our tests we consider three values of $\delta$: $\delta=0$, $\delta=0.01$ and $\delta=0.1$.

To avoid an inverse crime, we used a coarser mesh in our reconstructions.

In \eqref{eq:weakform} we set the value of contact impedances $z_\ell$, $\ell \in \{1,\dots,L\}$, to 0.1. 

The test case considered in all of our computational results is defined by a constant inclusion on a constant background
\[
 \sigma^{\ex} (x) = \left\{ \begin{array}{ll} 5, & \text{\; in \;} \Omega_h \text{\; if \;} \Omega_h \subset B_{0.5}(-0.3,-0.1)\\ 2, & \text{\; otherwise,} \end{array}\right.
\]
where $B_r(p_1,p_2) \subset \mathbb{R}^2$ is the ball centered at $(p_1,p_2)$ with radius $r$.

With each of the three above mentioned cost function combinations (all-at-once, eliminated $\sigma$, eliminated $(\Phi,\Psi)$, the iterates $x_k=(\sigma_k,\Phi_k,\Psi_k)$, or $x_k=(\Phi_k,\Psi_k)$, or $x_k=\sigma_k$, are defined by the projected gradient method \eqref{projgrad} from Section~\ref{sec:grad} 
where $\mu_k$ is found by an Armijo back tracking line search. 
Details on computation of the gradients of the various cost functions can be found in \cite{ThesisKha2021}.
The iteration is stopped by the discrepancy principle \eqref{discrprinc_grad} in the noisy case and as soon as the step size fell below a value $\epsilon_\mu$ (which we set to $10^{10}$ in our tests) in case of exact data.

\subsection{Numerical results for IAT}
We consider four cases of excitations, namely 
\begin{itemize}
\item $I=1$, with  $j_{1,1}=1$, $j_{5,1}=-1$ and $j_{k,1}=0$ otherwise;
\item $I=2$, with $j_{1,1}=j_{3,2}=1$, $j_{5,1}=j_{7,2}=-1$, and $j_{k,i}=0$ otherwise;
\item $I=4$, with
$j_{1,1}=j_{3,2}=j_{2,3}=j_{4,4}=1$, $j_{5,1}=j_{7,2}=j_{6,3}=j_{8,4}=-1$ and $j_{k,i}=0$ otherwise.
\item $I=28$, with all $({8\atop 2})$ combinations of setting $j_{k,i}=1$, $j_{\ell,i}=-1$ for $k\not=\ell\in\{1,\ldots,8\} $
\end{itemize}
The starting value is set to the mean value of the maximal and minimal value for the conductivity $\sigma_0 = \frac{1}{2} (\underline{\sigma} + \overline{\sigma})$ and $\Phi_0, \Psi_0$, if necessary, are gained from the weak form \eqref{eq:weakform} where $\sigma$ is replaced by $\sigma_0$.

The tables below show the data about \emph{the number of iterations}, \emph{the error $\| \sigma_{end} - \sigma^{\ex} \|_{L^2 (\Omega)}$}, \emph{the CPU time} (in seconds) and \emph{the CPU time for each iteration} for various versions of cost functions.
 
\begin{table}
\begin{tabular}{ll|r|r|r|r|}
\hline
\multicolumn{2}{|c|}{all-at-once version} & \multicolumn{1}{|c|}{number} & \multicolumn{1}{|c|}{$L^2$ error} & \multicolumn{1}{|c|}{CPU-time} & \multicolumn{1}{|c|}{CPU-time} \\
\multicolumn{2}{|c|}{IAT} & \multicolumn{1}{|c|}{of iterations} & \multicolumn{1}{|c|}{$\|\sigma_{end} - \sigma^{\ex}\|$} & \multicolumn{1}{|c|}{(in seconds)} & \multicolumn{1}{|c|}{per iteration} \\
\hline

\hline
\multicolumn{1}{|l|}{\multirow{3}{*}{$I=1$}} & \multicolumn{1}{|l|}{$\delta = 0$} & 5 025 130 & 0.658 72 & 943 527 & 0.187 76 \\
\cline{2-6}
\multicolumn{1}{|c|}{}&\multicolumn{1}{|l|}{$\delta = 0.01$} & 4 959 452 & 0.665 37 & 930 984 & 0.187 72  \\
\cline{2-6}
\multicolumn{1}{|c|}{}&\multicolumn{1}{|l|}{$\delta = 0.1$} & 5 178 542 & 0.806 76 & 984 851 & 0.190 18 \\
\cline{2-6}
\hline

\hline
\multicolumn{1}{|l|}{\multirow{3}{*}{$I=2$}} & \multicolumn{1}{|l|}{$\delta = 0$} & 1 109 245 & 0.387 57 & 266 270 & 0.240 05 \\
\cline{2-6}
\multicolumn{1}{|c|}{}&\multicolumn{1}{|l|}{$\delta = 0.01$} & 1 114 829 & 0.387 46 & 264 697 & 0.237 43 \\
\cline{2-6}
\multicolumn{1}{|c|}{}&\multicolumn{1}{|l|}{$\delta = 0.1$} & 1 520 239 & 0.567 00 & 344 894 & 0.226 87 \\
\cline{2-6}
\hline

\hline
\multicolumn{1}{|l|}{\multirow{3}{*}{$I=4$}} & \multicolumn{1}{|l|}{$\delta = 0$} & 301 651 & 0.310 70 & 73 664 & 0.244 20 \\
\cline{2-6}
\multicolumn{1}{|c|}{}&\multicolumn{1}{|l|}{$\delta = 0.01$} & 308 170 & 0.314 96 & 74 474 & 0.241 66 \\
\cline{2-6}
\multicolumn{1}{|c|}{}&\multicolumn{1}{|l|}{$\delta = 0.1$} & 326 561 & 0.412 61 & 79 693 & 0.244 04 \\
\cline{2-6}
\hline

\hline
\multicolumn{1}{|l|}{\multirow{3}{*}{$I=28$}} & \multicolumn{1}{|l|}{$\delta = 0$} & 249 816 & 0.305 35 & 82 681 & 0.330 97 \\
\cline{2-6}
\multicolumn{1}{|c|}{}&\multicolumn{1}{|l|}{$\delta = 0.01$} & 245 306 & 0.306 75 & 81 263 & 0.331 27 \\
\cline{2-6}
\multicolumn{1}{|c|}{}&\multicolumn{1}{|l|}{$\delta = 0.1$} & 292 914 & 0.326 76 & 96 378 & 0.329 03 \\
\cline{2-6}
\hline
\end{tabular}
\caption{IAT, all-at-once version \eqref{eq:costfunJmod_KV}, \eqref{eq:costfunJobs_IAT}. \label{tab:IAT_all-at-onceVer}}
\end{table}

\begin{table}
\begin{tabular}{ll|r|r|r|r|}
\hline
\multicolumn{2}{|c|}{eliminating-$\sigma$} & \multicolumn{1}{|c|}{number} & \multicolumn{1}{|c|}{$L^2$ error} & \multicolumn{1}{|c|}{CPU-time} & \multicolumn{1}{|c|}{CPU-time} \\
\multicolumn{2}{|c|}{version, IAT} & \multicolumn{1}{|c|}{of iterations} & \multicolumn{1}{|c|}{$\|\sigma_{end} - \sigma^{\ex}\|$} & \multicolumn{1}{|c|}{(in seconds)} & \multicolumn{1}{|c|}{per iteration} \\
\hline

\hline
\multicolumn{1}{|l|}{\multirow{3}{*}{$I=1$}} & \multicolumn{1}{|l|}{$\delta = 0$} & 2 224 & 0.582 89 & 583 & 0.262 28 \\
\cline{2-6}
\multicolumn{1}{|c|}{}&\multicolumn{1}{|l|}{$\delta = 0.01$} & 1 230 & 0.797 77 & 319 & 0.259 46  \\
\cline{2-6}
\multicolumn{1}{|c|}{}&\multicolumn{1}{|l|}{$\delta = 0.1$} & 840 & 0.944 67 & 215 & 0.256 36 \\
\cline{2-6}
\hline

\hline
\multicolumn{1}{|l|}{\multirow{3}{*}{$I=2$}} & \multicolumn{1}{|l|}{$\delta = 0$} & 3 488 & 0.360 50 & 1 232 & 0.353 21 \\
\cline{2-6}
\multicolumn{1}{|c|}{}&\multicolumn{1}{|l|}{$\delta = 0.01$} & 2 894 & 0.373 53 & 920 & 0.317 86 \\
\cline{2-6}
\multicolumn{1}{|c|}{}&\multicolumn{1}{|l|}{$\delta = 0.1$} & 2 645 & 0.455 20 & 959 & 0.362 74 \\
\cline{2-6}
\hline

\hline
\multicolumn{1}{|l|}{\multirow{3}{*}{$I=4$}} & \multicolumn{1}{|l|}{$\delta = 0$} & 2 730 & 0.324 27 & 801 & 0.293 49 \\
\cline{2-6}
\multicolumn{1}{|c|}{}&\multicolumn{1}{|l|}{$\delta = 0.01$} & 3 393 & 0.321 52 & 1 254 & 0.369 59 \\
\cline{2-6}
\multicolumn{1}{|c|}{}&\multicolumn{1}{|l|}{$\delta = 0.1$} & 1 836 & 0.433 52 & 562 & 0.305 89 \\
\cline{2-6}
\hline

\hline
\multicolumn{1}{|l|}{\multirow{3}{*}{$I=28$}} & \multicolumn{1}{|l|}{$\delta = 0$} & 3 330 & 0.330 32 & 1 451 & 0.435 61 \\
\cline{2-6}
\multicolumn{1}{|c|}{}&\multicolumn{1}{|l|}{$\delta = 0.01$} & 3 251 & 0.331 05 & 1 427 & 0.439 00 \\
\cline{2-6}
\multicolumn{1}{|c|}{}&\multicolumn{1}{|l|}{$\delta = 0.1$} & 3 300 & 0.344 39 & 1 511 & 0.457 87 \\
\cline{2-6}
\hline
\end{tabular}
\caption{IAT, eliminating-$\sigma$ version \eqref{eq:costfunJ_IAT_sigma}. \label{tab:IAT_eliminatingsigmaVer}}
\end{table}

\begin{table}
\begin{tabular}{ll|r|r|r|r|}
\hline
\multicolumn{2}{|c|}{eliminating-$(\Phi,\Psi)$} & \multicolumn{1}{|c|}{number} & \multicolumn{1}{|c|}{$L^2$ error} & \multicolumn{1}{|c|}{CPU-time} & \multicolumn{1}{|c|}{CPU-time} \\
\multicolumn{2}{|c|}{version, IAT} & \multicolumn{1}{|c|}{of iterations} & \multicolumn{1}{|c|}{$\|\sigma_{end} - \sigma^{\ex}\|$} & \multicolumn{1}{|c|}{(in seconds)} & \multicolumn{1}{|c|}{per iteration} \\
\hline
  
\hline
\multicolumn{1}{|l|}{\multirow{3}{*}{$I=1$}} & \multicolumn{1}{|l|}{$\delta = 0$} & 225 302 &  2.73e-09  &   147 936     &    0.656 61  \\
\cline{2-6}
\multicolumn{1}{|c|}{}&\multicolumn{1}{|l|}{$\delta = 0.01$} & 100 041  &  0.018 45   &   83 302    &     0.832 68  \\
\cline{2-6}
\multicolumn{1}{|c|}{}&\multicolumn{1}{|l|}{$\delta = 0.1$} & 79 377  &  0.224 34   &   38 043   &    0.479 27  \\
\cline{2-6}
\hline
     
\hline
\multicolumn{1}{|l|}{\multirow{3}{*}{$I=2$}} & \multicolumn{1}{|l|}{$\delta = 0$} & 55 066 &  1.68e-10   &   30 161 &    0.547 72 \\
\cline{2-6}
\multicolumn{1}{|c|}{}&\multicolumn{1}{|l|}{$\delta = 0.01$} & 61 162  &  0.016 69  &    61 202   &      1.000 65 \\
\cline{2-6}
\multicolumn{1}{|c|}{}&\multicolumn{1}{|l|}{$\delta = 0.1$} & 38 162  &  0.174 07   &   41 314      &   1.082 60 \\
\cline{2-6}
\hline

\hline
\multicolumn{1}{|l|}{\multirow{3}{*}{$I=4$}} & \multicolumn{1}{|l|}{$\delta = 0$} & 13 782 &  4.38e-11  &    10 469      &   0.759 58  \\
\cline{2-6}
\multicolumn{1}{|c|}{}&\multicolumn{1}{|l|}{$\delta = 0.01$} & 19 889  &  0.014 41   &   27 064    &   1.360 76 \\
\cline{2-6}
\multicolumn{1}{|c|}{}&\multicolumn{1}{|l|}{$\delta = 0.1$} & 15 169  &  0.122 03   &   19 721    &     1.300 06 \\
\cline{2-6}
\hline
 
\hline
\multicolumn{1}{|l|}{\multirow{3}{*}{$I=28$}} & \multicolumn{1}{|l|}{$\delta = 0$} & 13 868  & 4.38e-11   &   49 005   &   3.533 65  \\
\cline{2-6}
\multicolumn{1}{|c|}{}&\multicolumn{1}{|l|}{$\delta = 0.01$} &  23 540  &  0.007 44   &   96 096    &     4.082 26 \\
\cline{2-6}
\multicolumn{1}{|c|}{}&\multicolumn{1}{|l|}{$\delta = 0.1$} & 28 721  &  0.067 32  &   120 978    &     4.212 17 \\
\cline{2-6}
\hline
\end{tabular}
\caption{IAT, eliminating-$(\Phi,\Psi)$ version \eqref{eq:costfunJ_IAT_PhiPsi}. \label{tab:IAT_eliminatingPhiPsiVer}}
\end{table}

\begin{figure}
\includegraphics[width=0.29\textwidth]{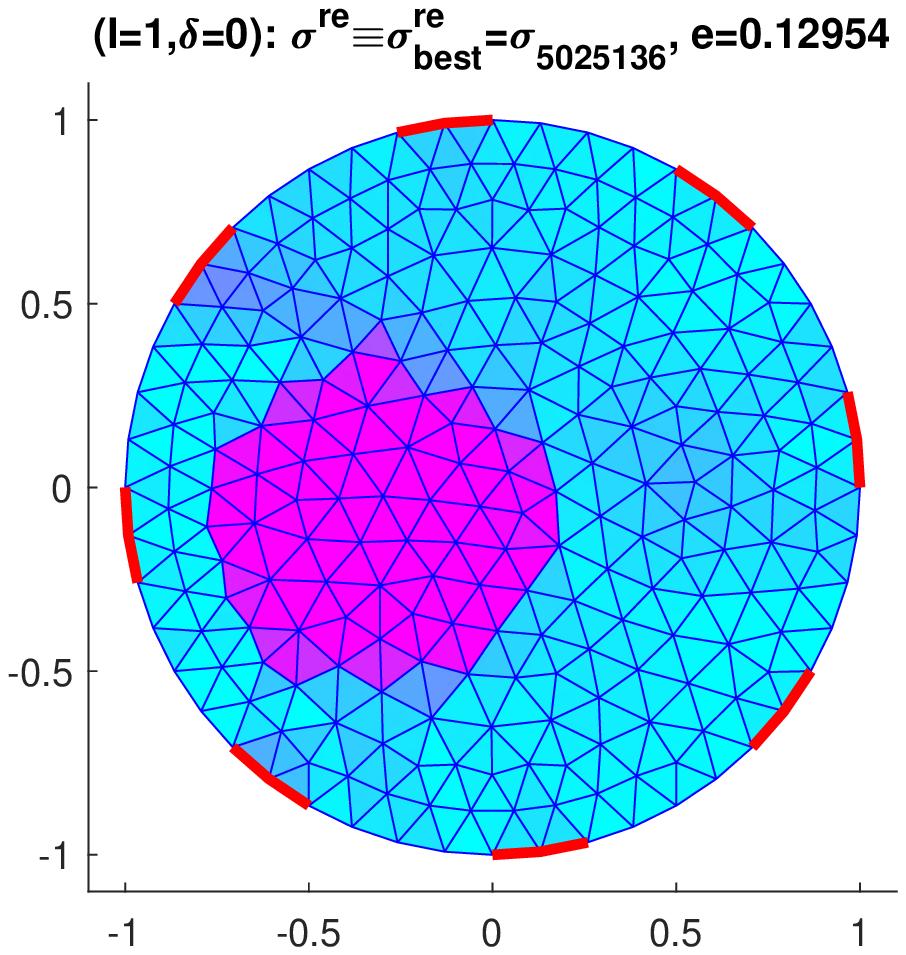}
\includegraphics[width=0.29\textwidth]{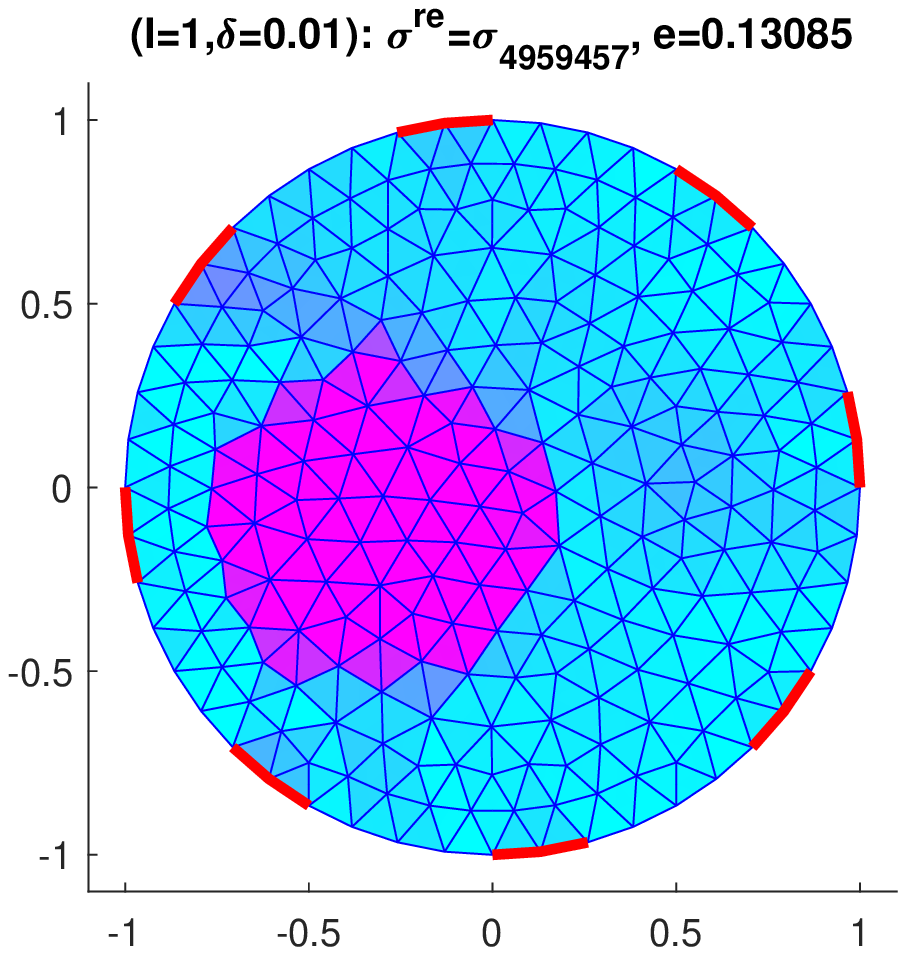}
\includegraphics[width=0.29\textwidth]{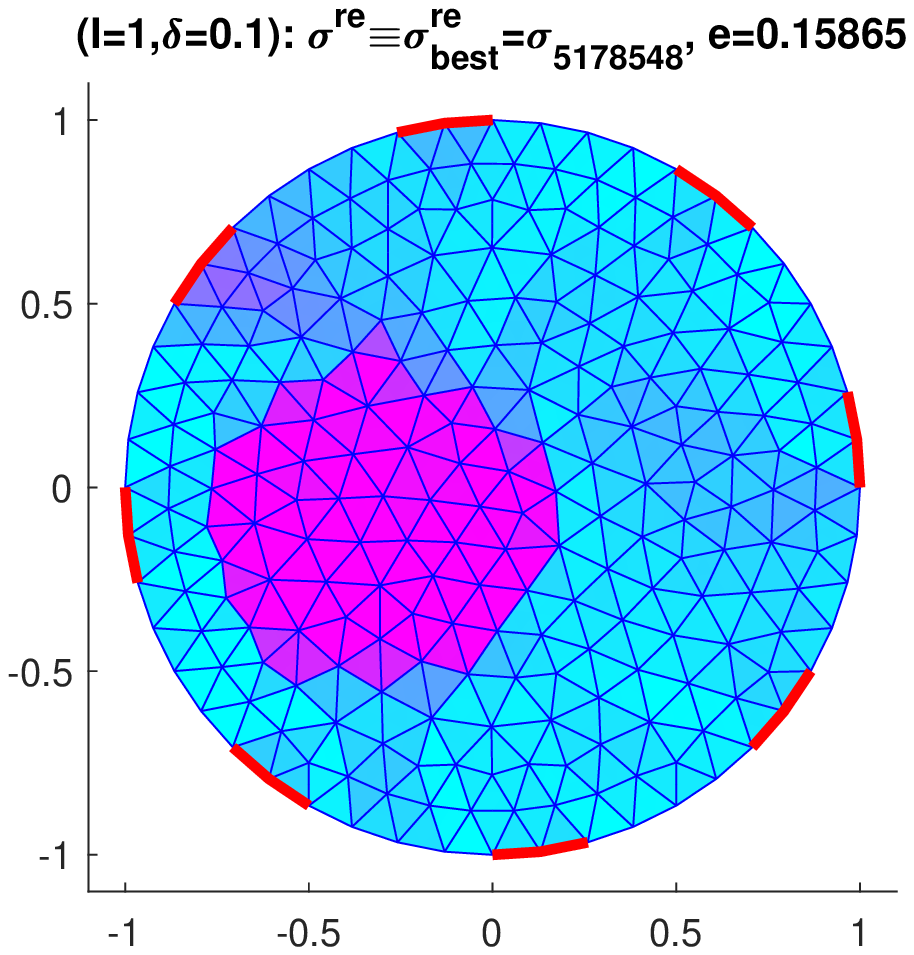}\\
\includegraphics[width=0.29\textwidth]{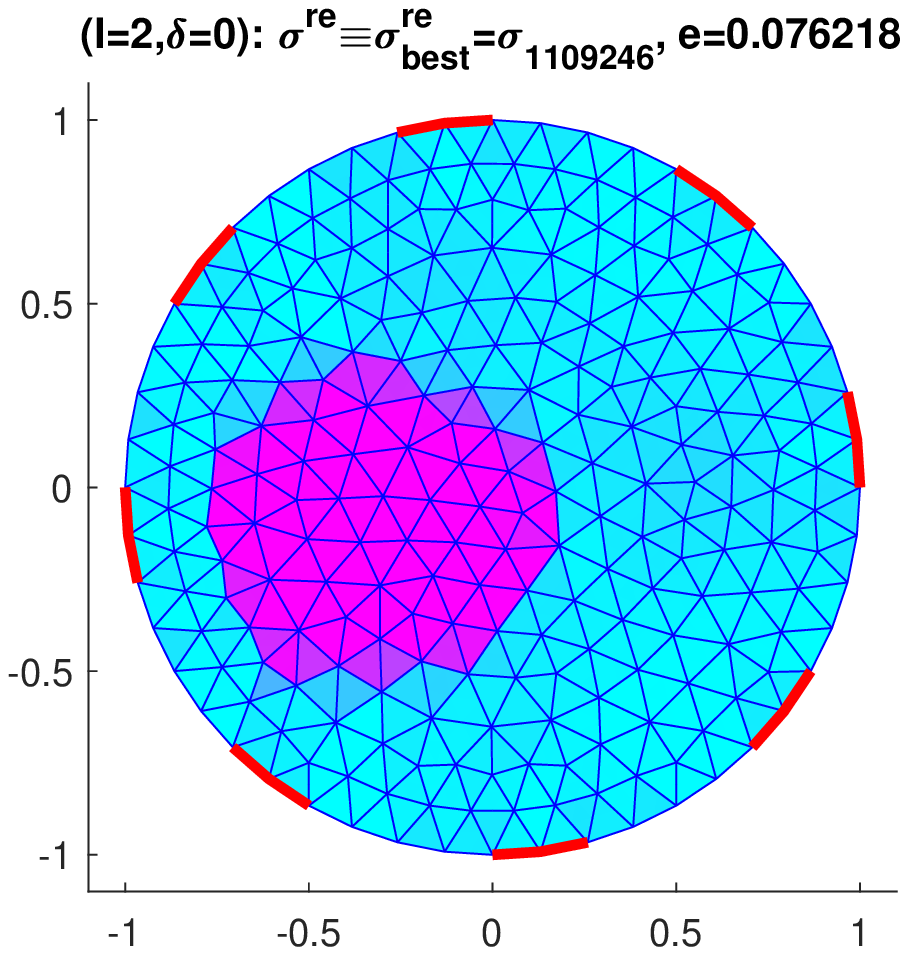}
\includegraphics[width=0.29\textwidth]{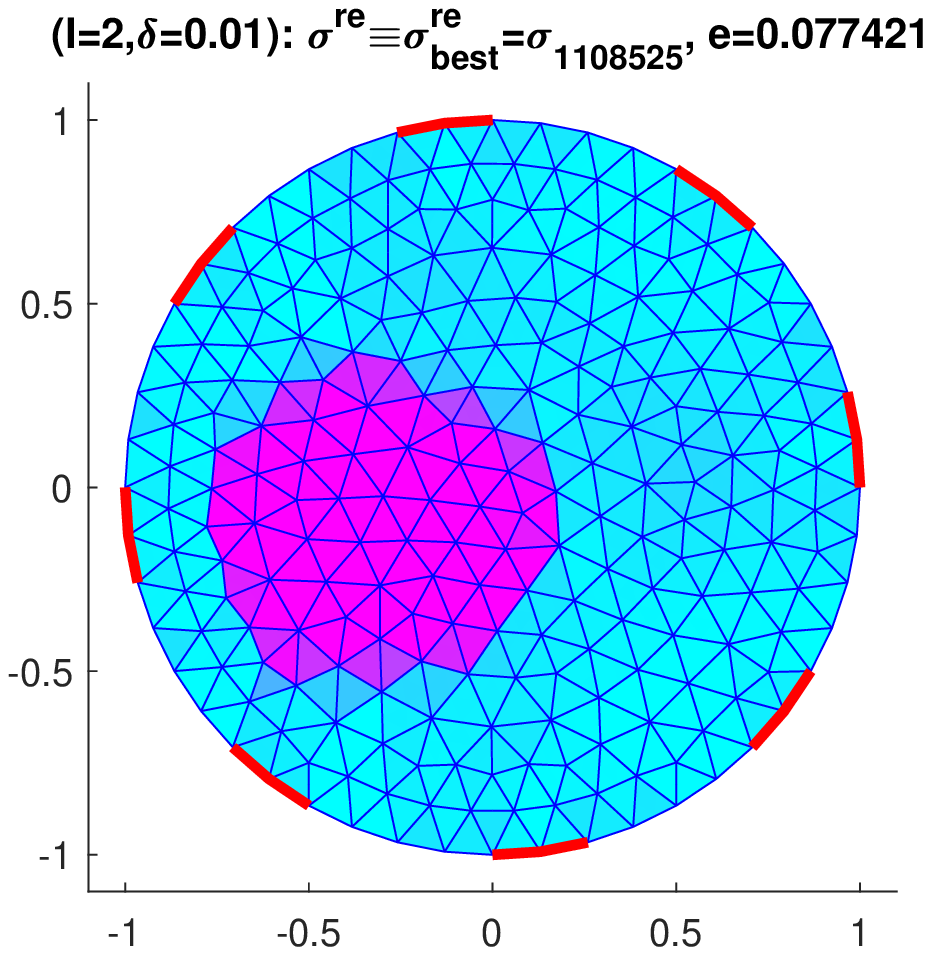}
\includegraphics[width=0.29\textwidth]{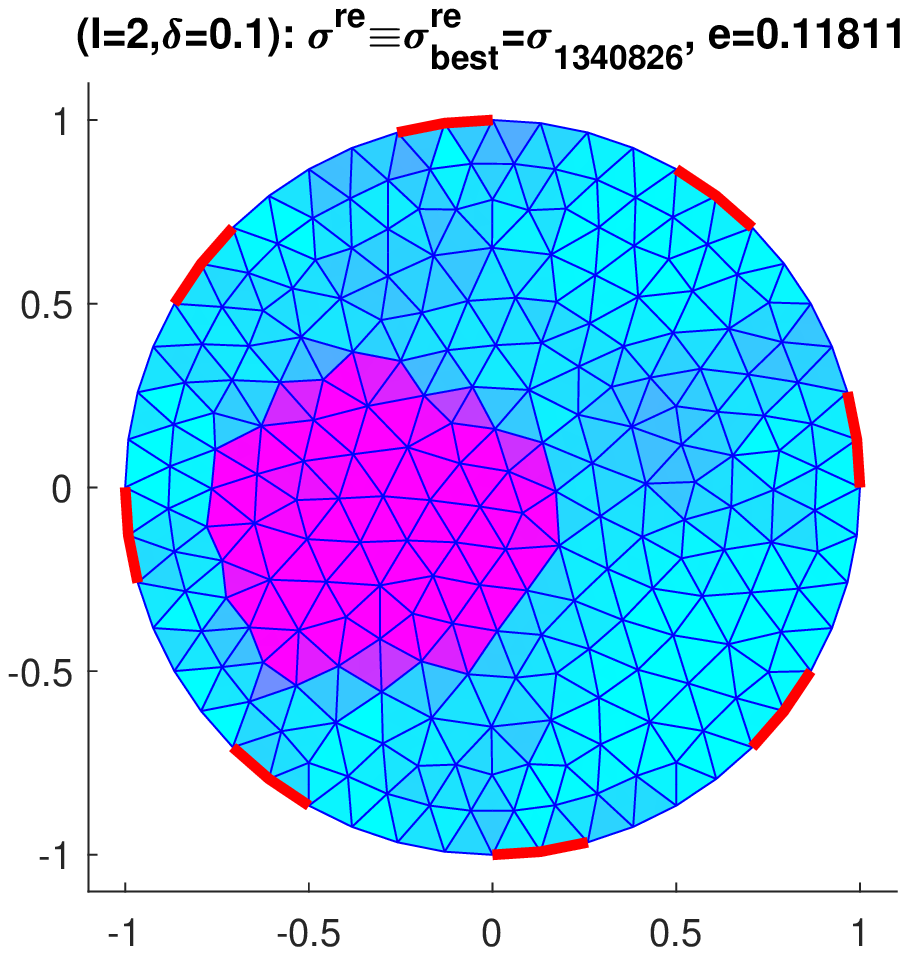}\\
\includegraphics[width=0.29\textwidth]{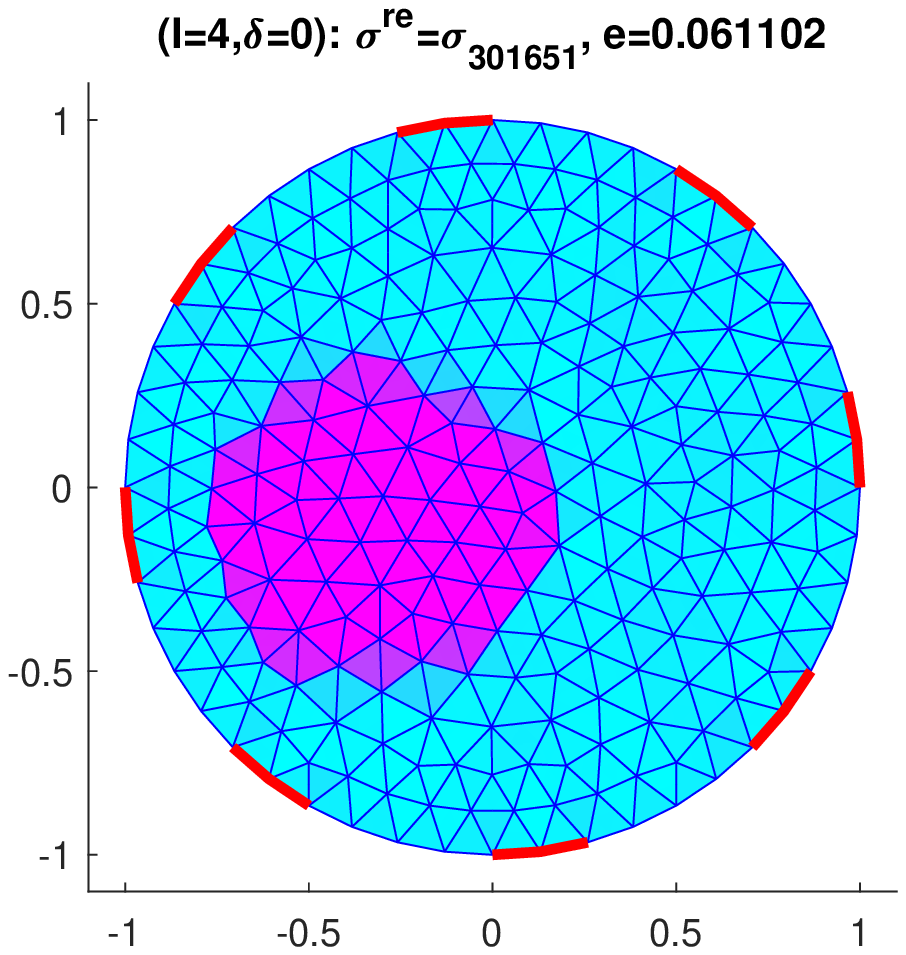}
\includegraphics[width=0.29\textwidth]{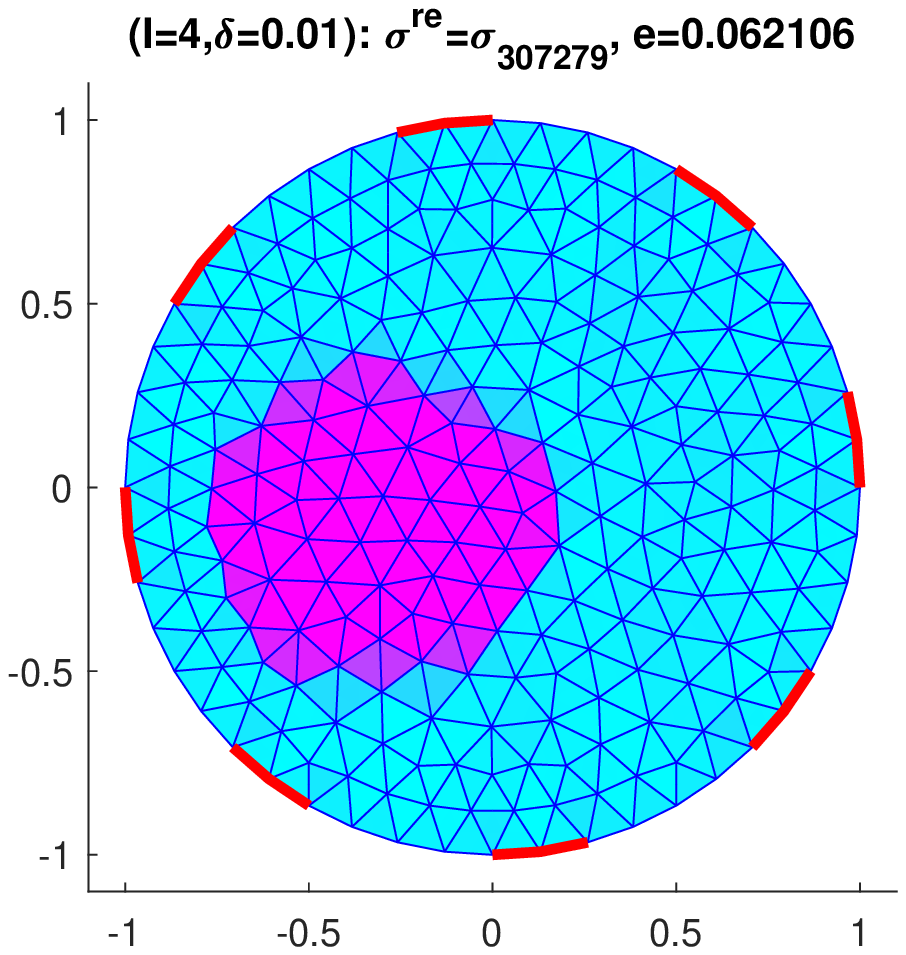}
\includegraphics[width=0.29\textwidth]{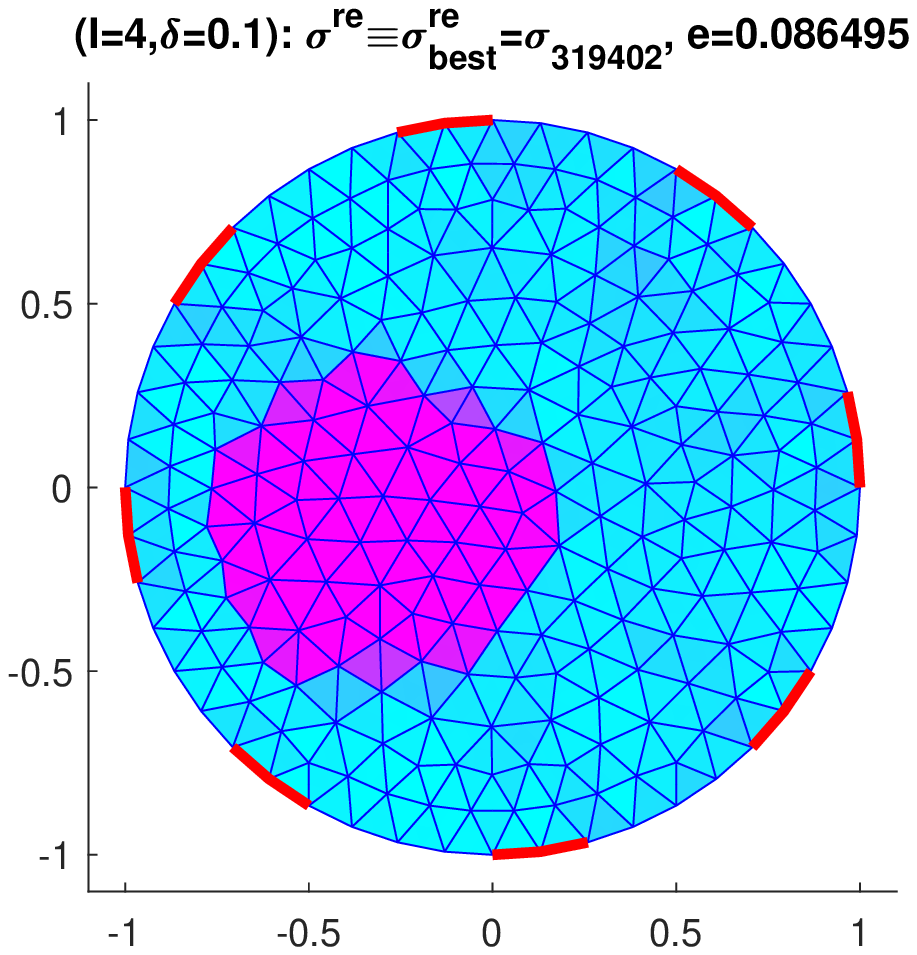}\\
\includegraphics[width=0.29\textwidth]{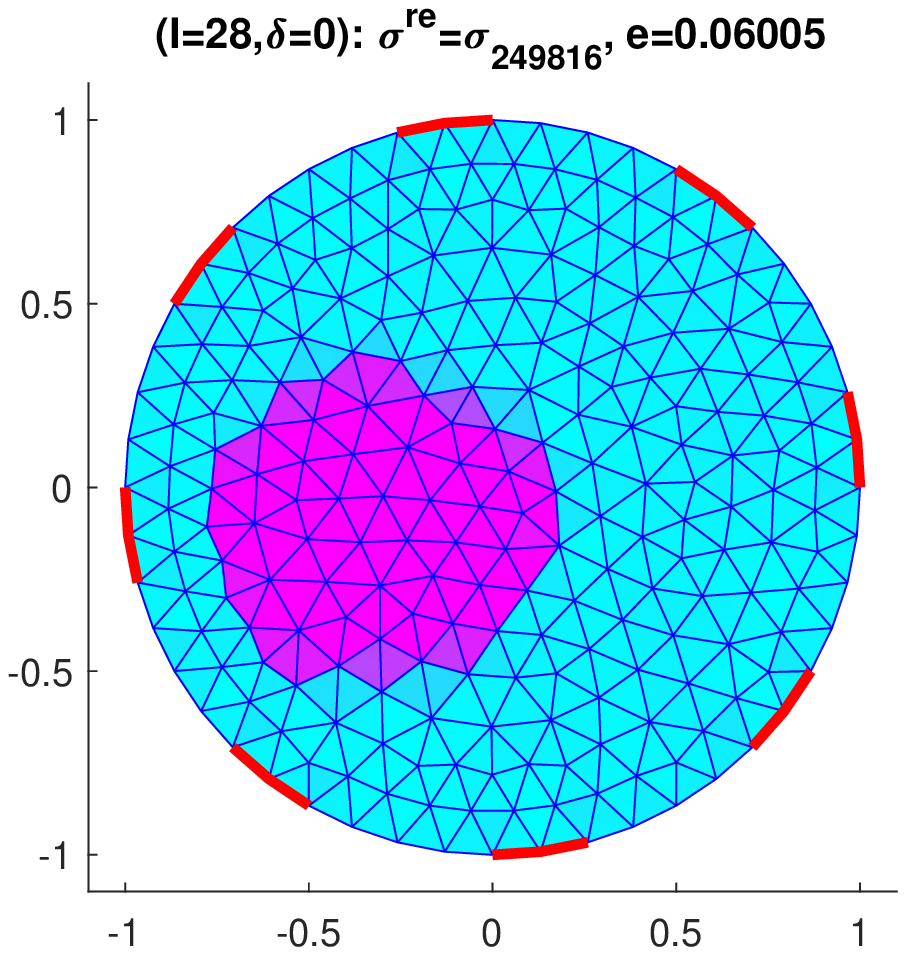}
\includegraphics[width=0.29\textwidth]{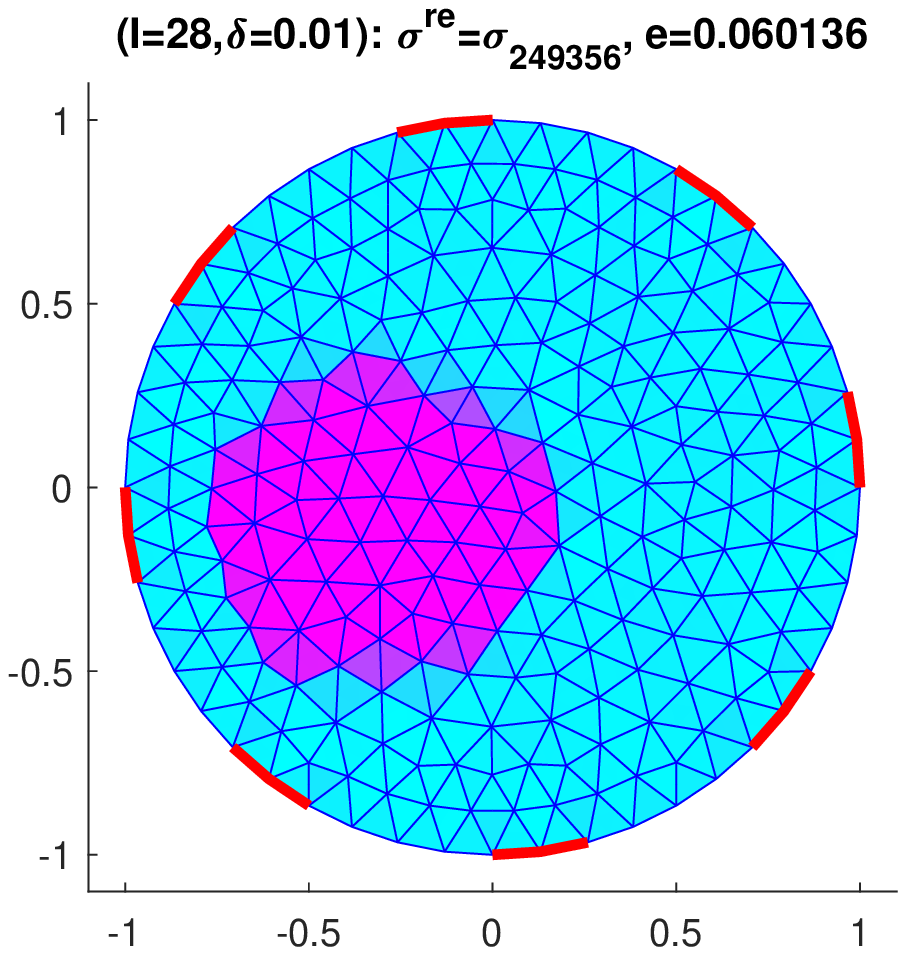}
\includegraphics[width=0.29\textwidth]{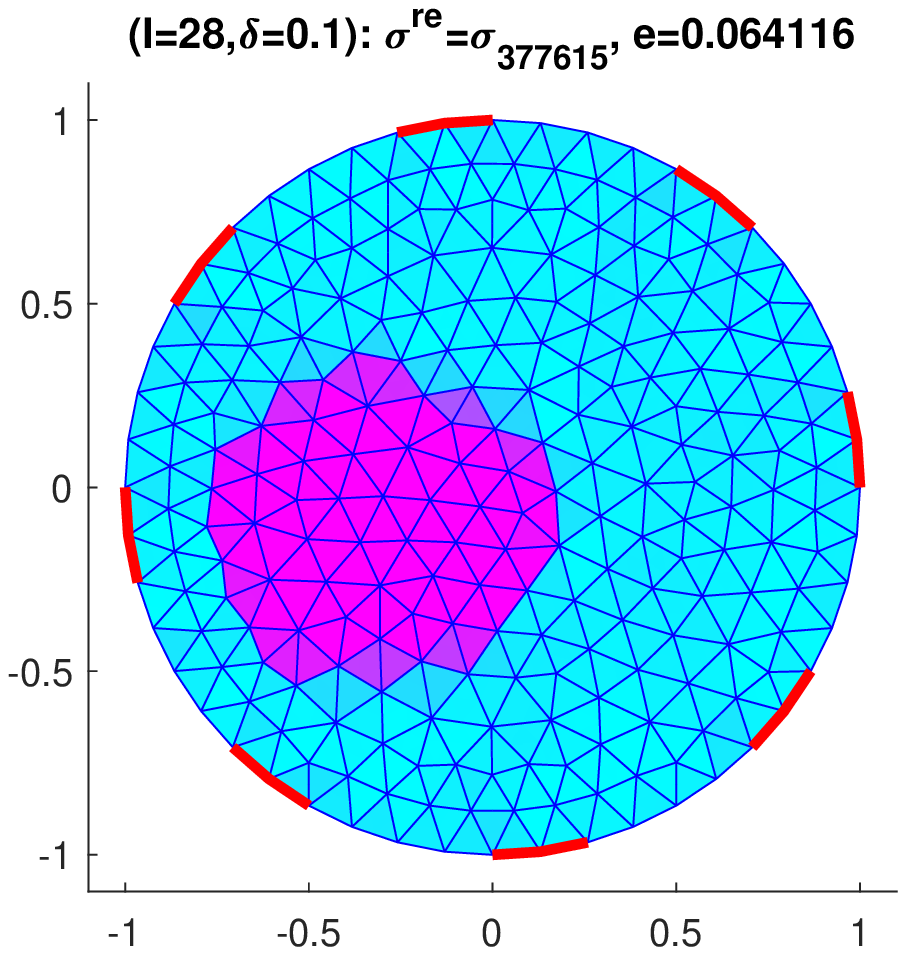}\\
\caption{Reconstructions of $\sigma$ from all-at-once version of cost function IAT \eqref{eq:costfunJmod_KV}, \eqref{eq:costfunJobs_IAT}, in cases $I=1$, $I=2$, $I=4$, $I=28$ (top to bottom) for $\delta=0$, $\delta=0.01$, $\delta=0.1$ (left to right). 
\label{fig:sigmaIAT_all-at-onceVer}}
\end{figure}

\begin{figure}
\includegraphics[width=0.29\textwidth]{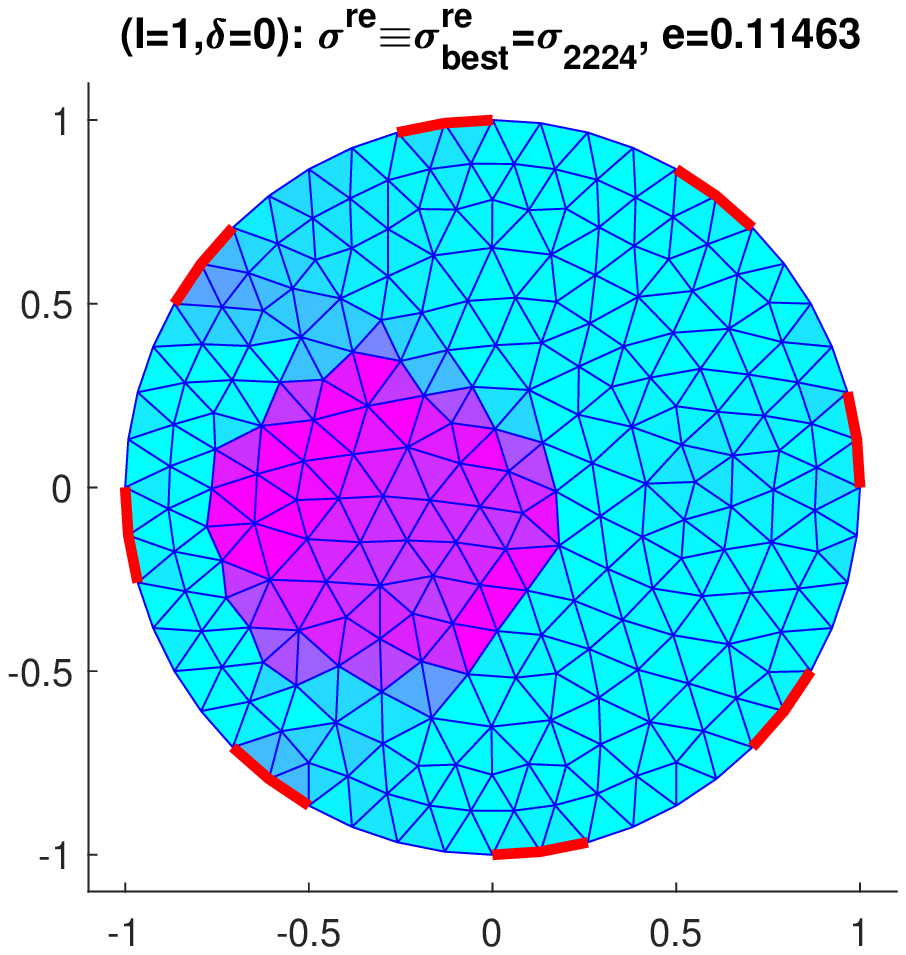}
\includegraphics[width=0.29\textwidth]{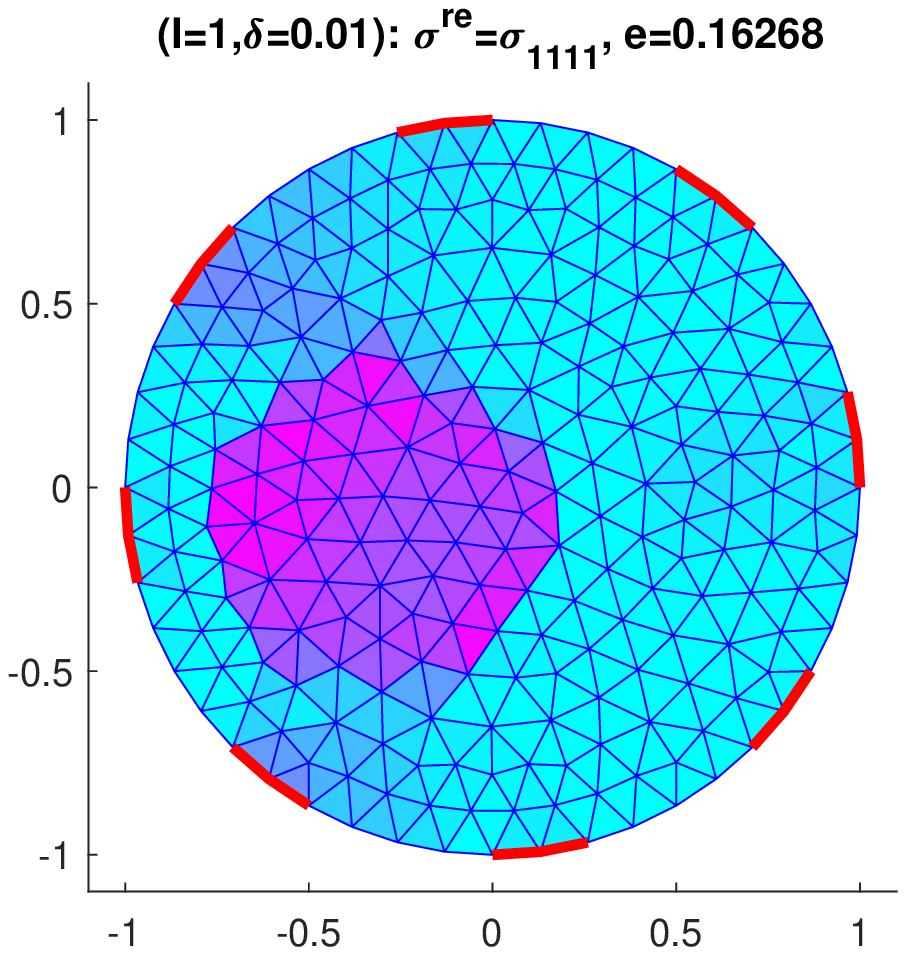}
\includegraphics[width=0.29\textwidth]{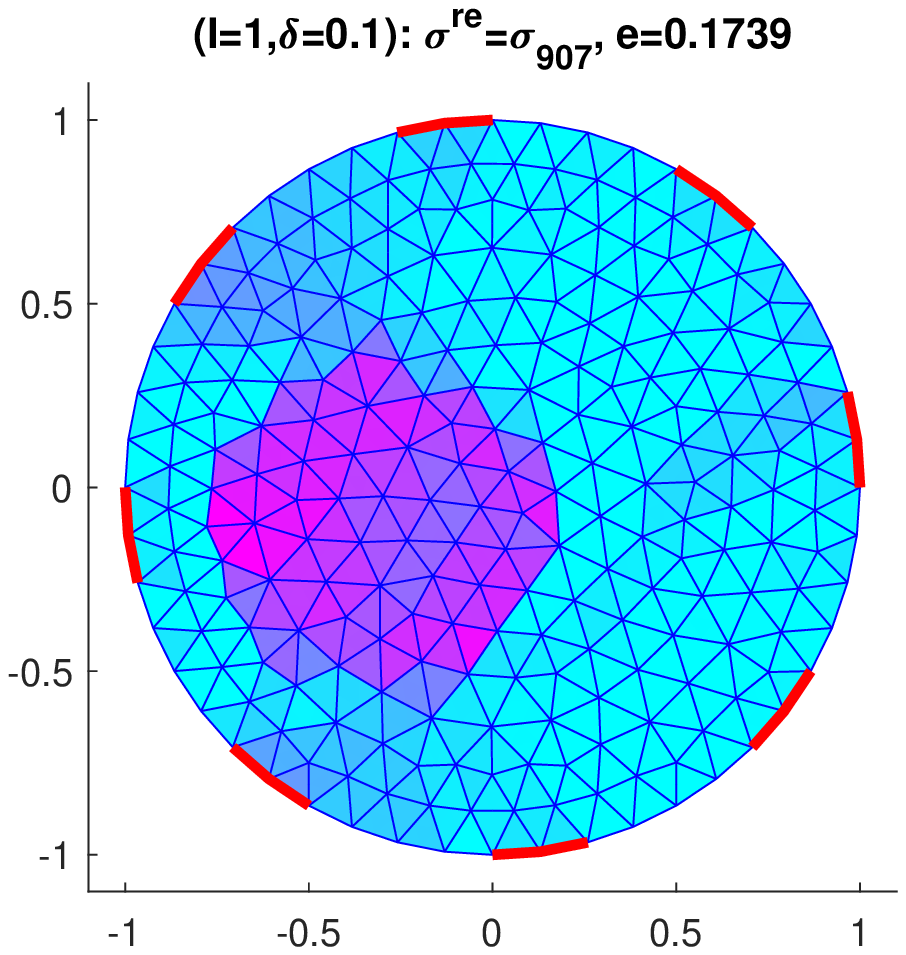}\\
\includegraphics[width=0.29\textwidth]{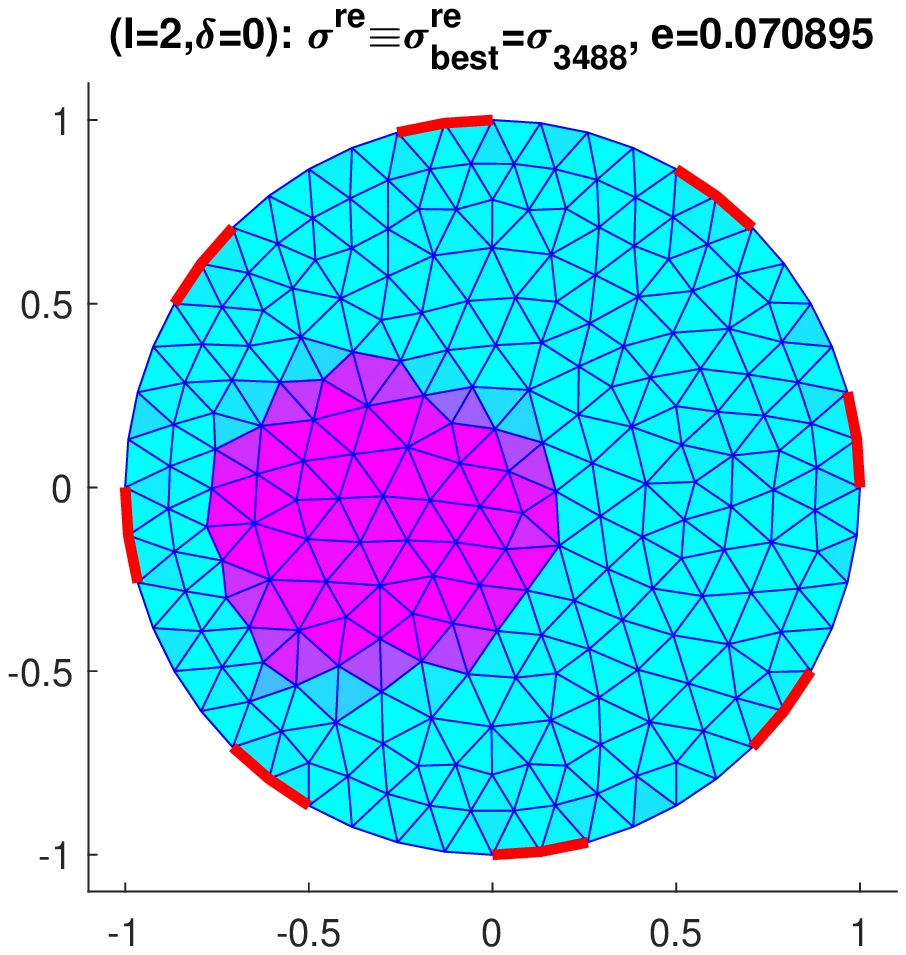}
\includegraphics[width=0.29\textwidth]{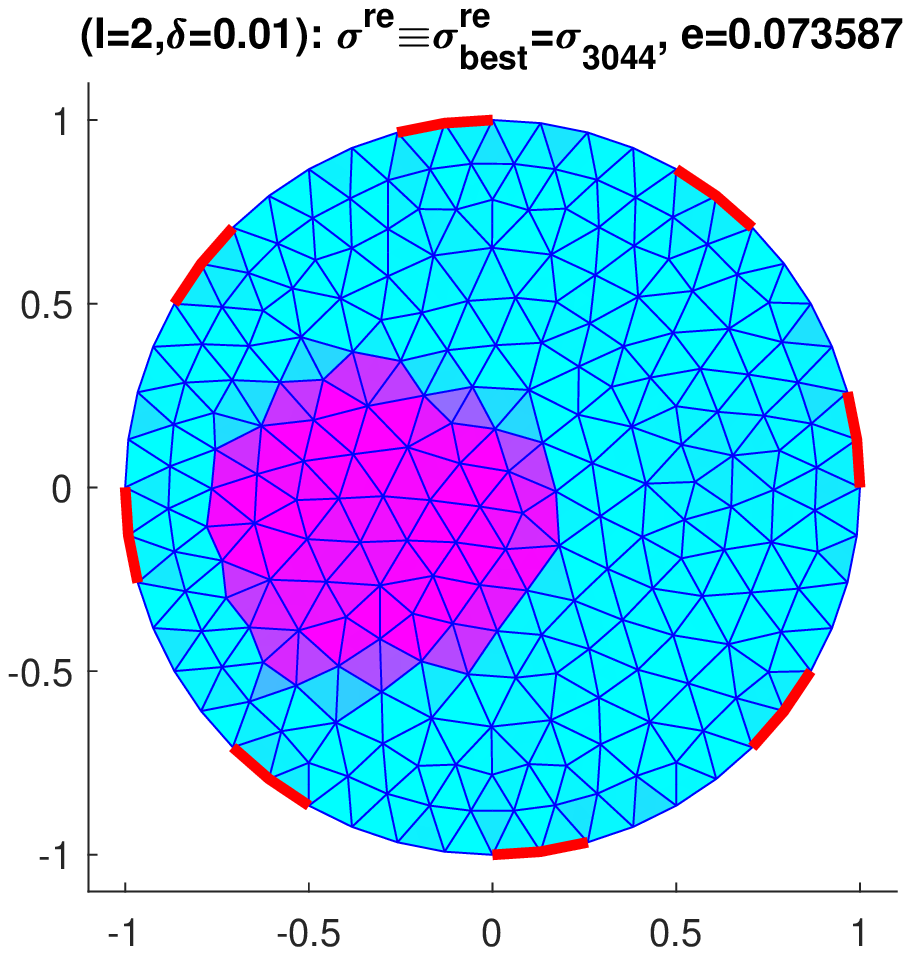}
\includegraphics[width=0.29\textwidth]{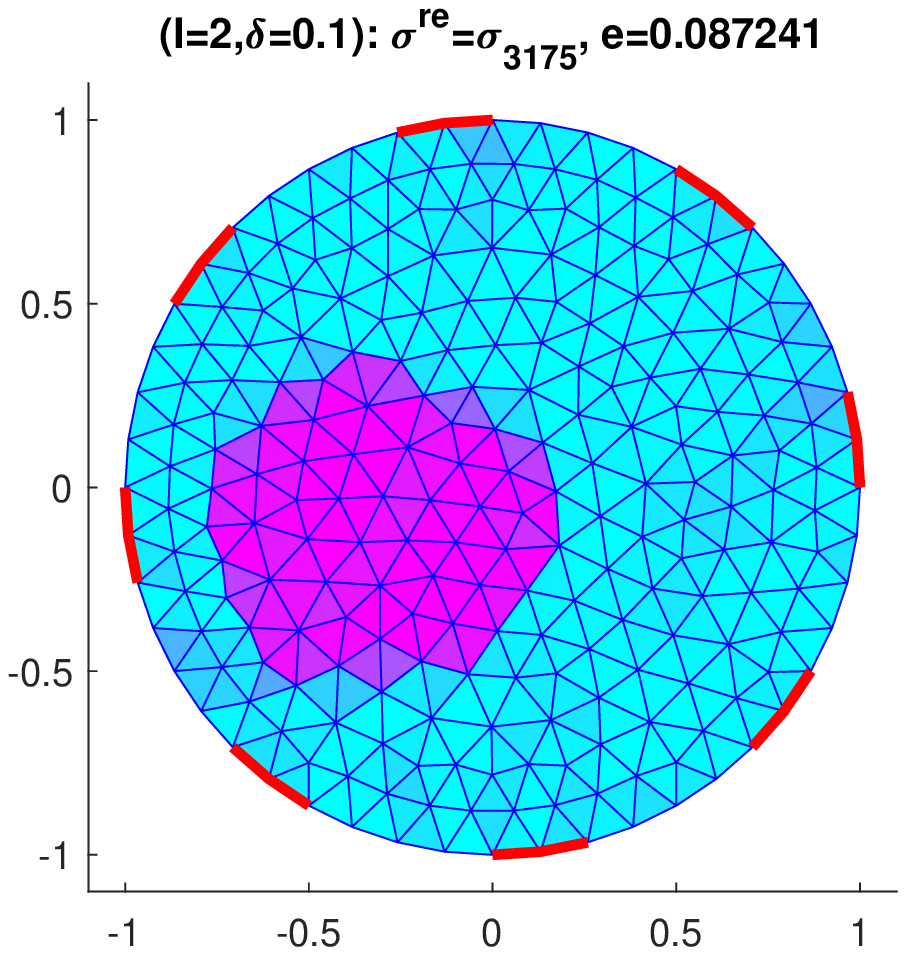}\\
\includegraphics[width=0.29\textwidth]{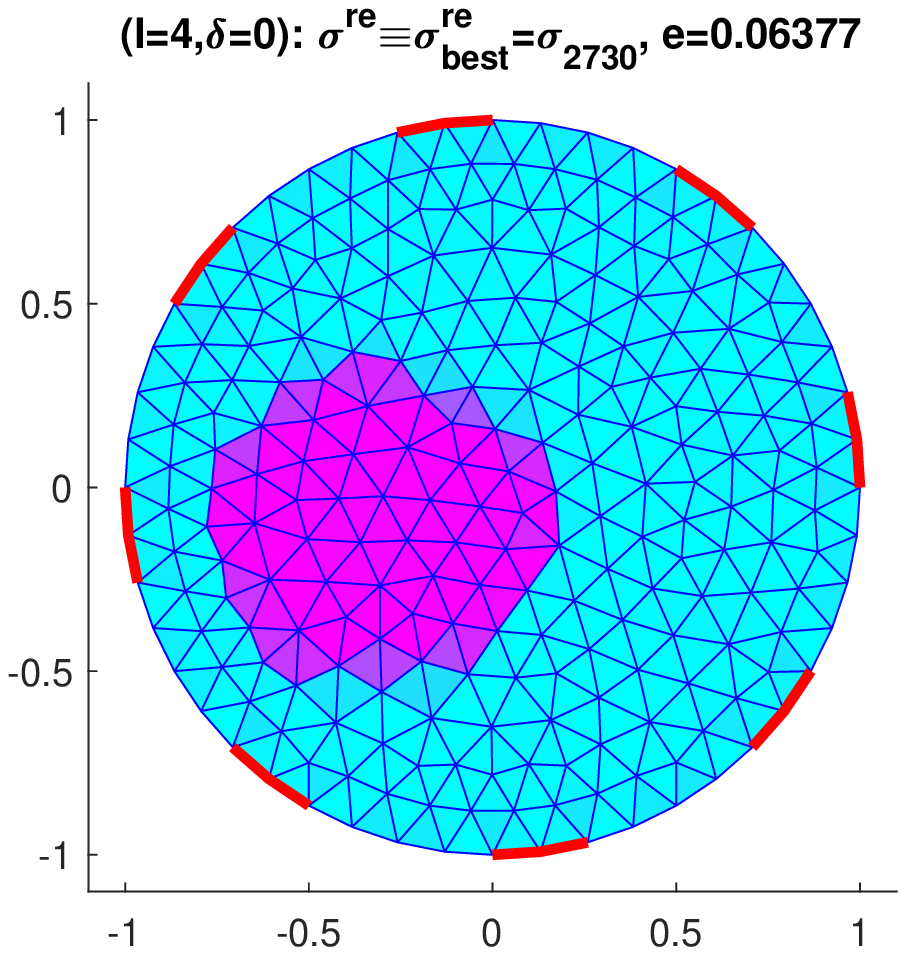}
\includegraphics[width=0.29\textwidth]{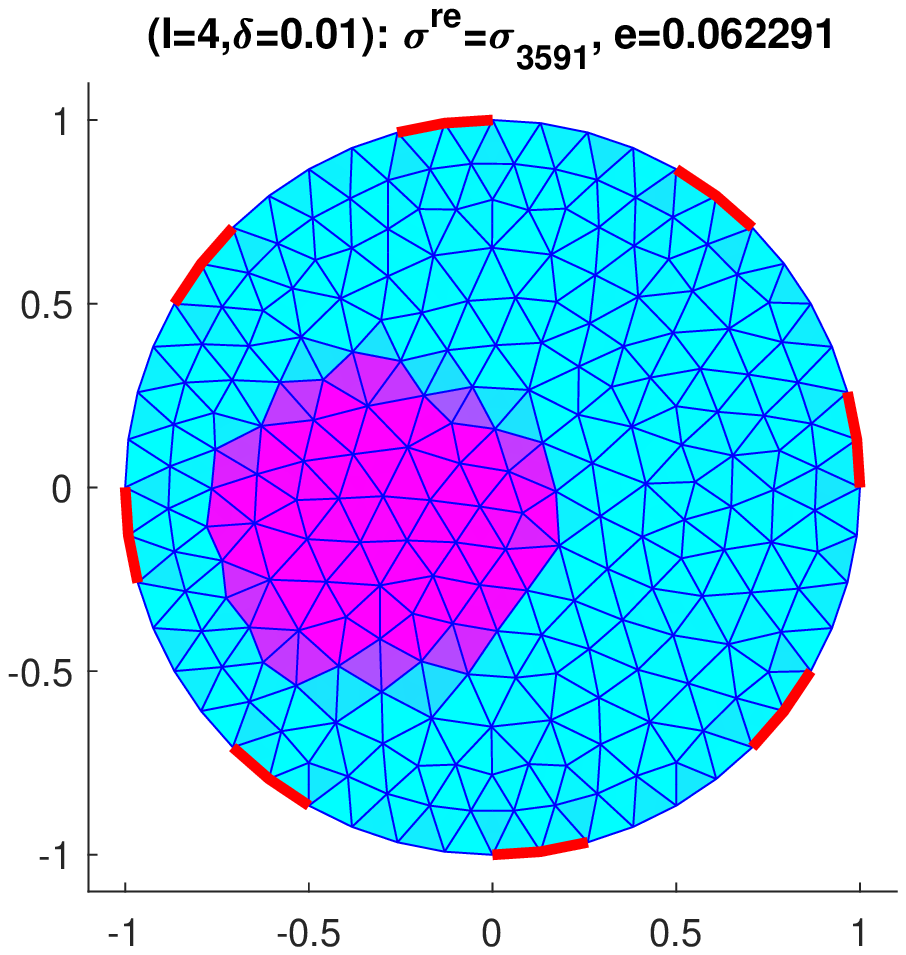}
\includegraphics[width=0.29\textwidth]{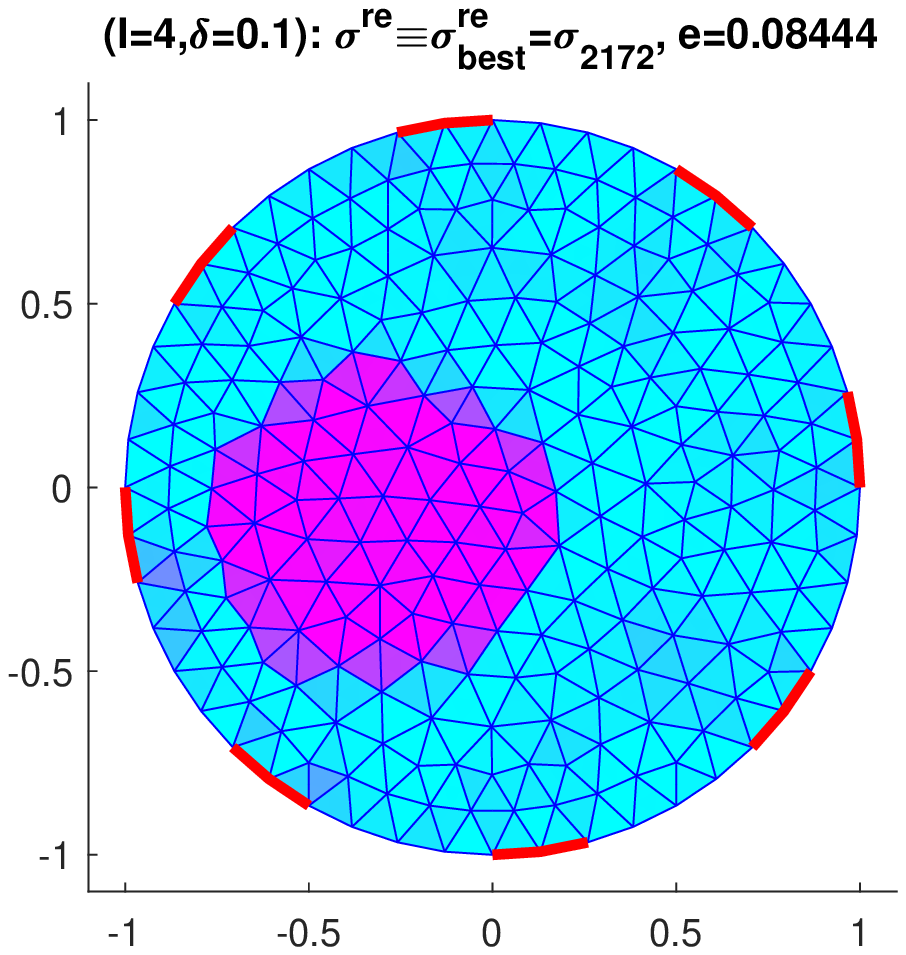}\\
\includegraphics[width=0.29\textwidth]{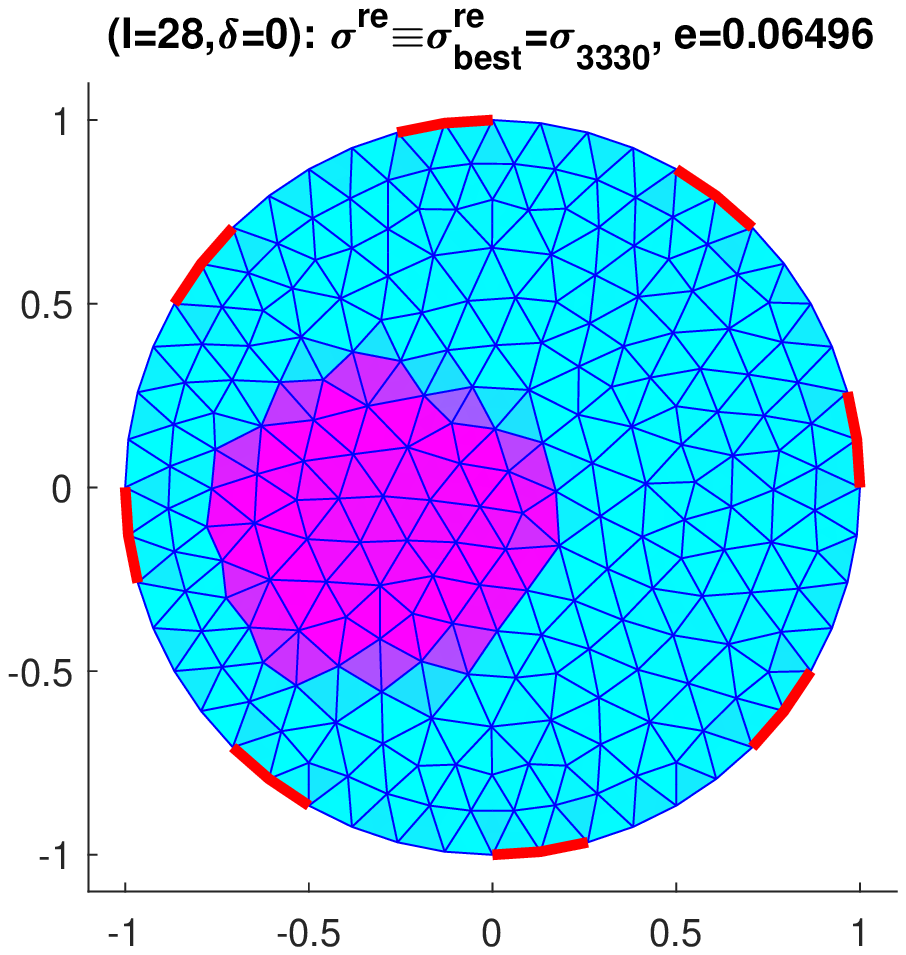}
\includegraphics[width=0.29\textwidth]{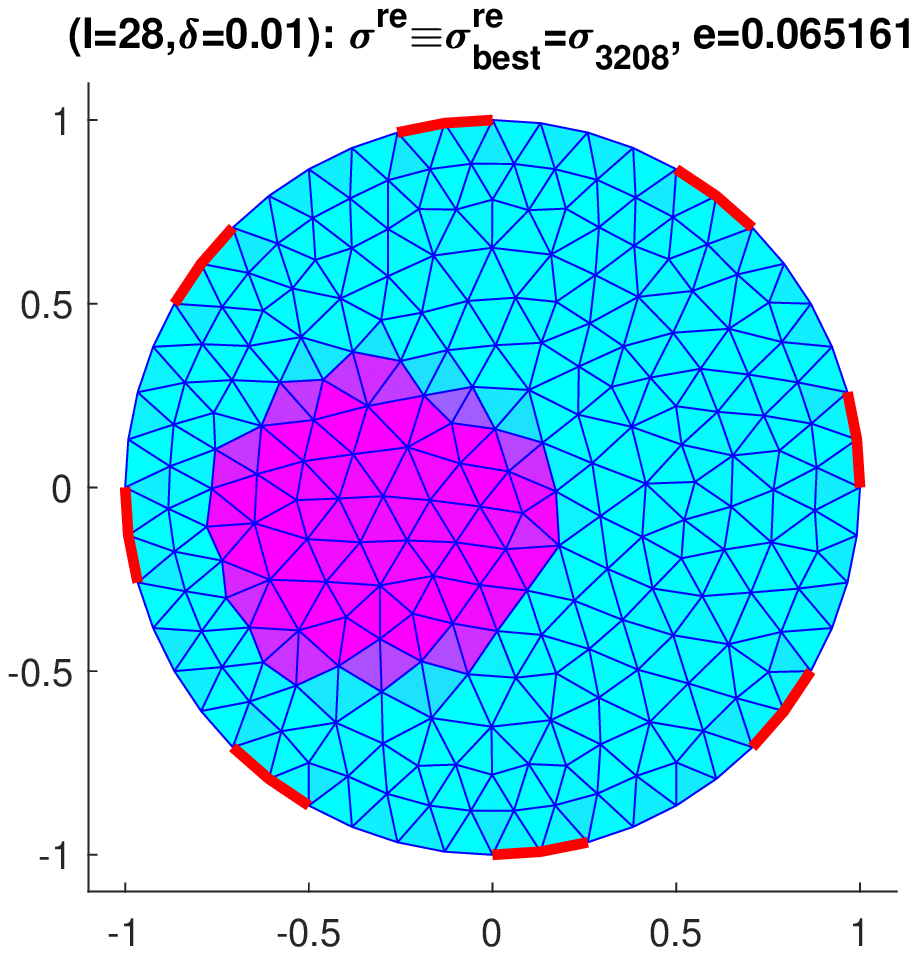}
\includegraphics[width=0.29\textwidth]{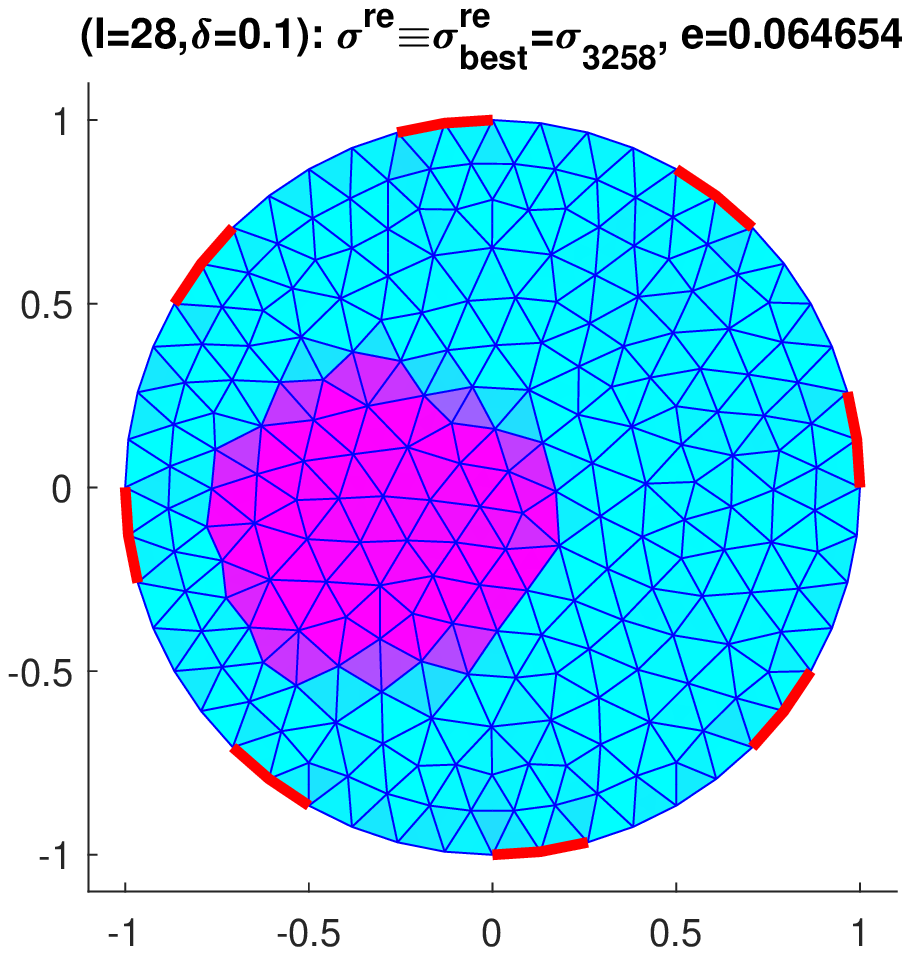}\\
\caption{Reconstructions of $\sigma$ from eliminating-$\sigma$ version of cost function IAT \eqref{eq:costfunJ_IAT_sigma}, in cases $I=1$, $I=2$, $I=4$, $I=28$ (top to bottom) for $\delta=0$, $\delta=0.01$, $\delta=0.1$ (left to right). 
\label{fig:sigma12IAT_eliminatingsigmaVer}}
\end{figure}

\begin{figure}
\includegraphics[width=0.29\textwidth]{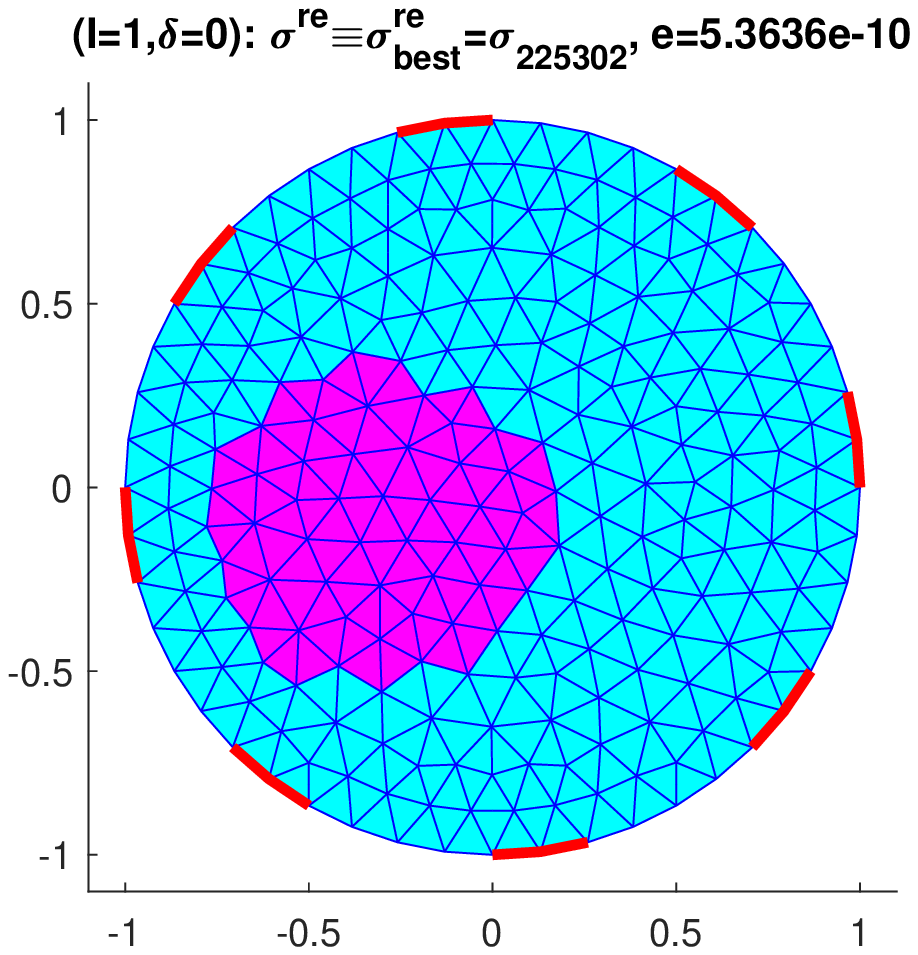}
\includegraphics[width=0.29\textwidth]{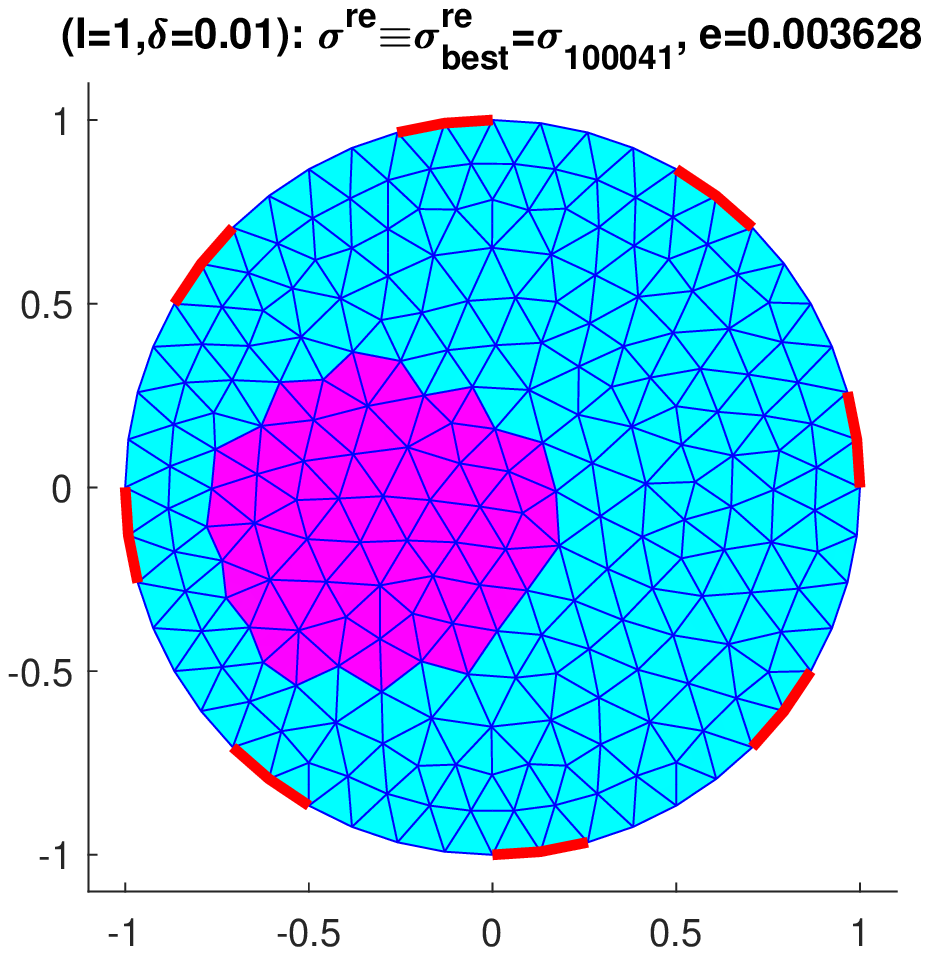}
\includegraphics[width=0.29\textwidth]{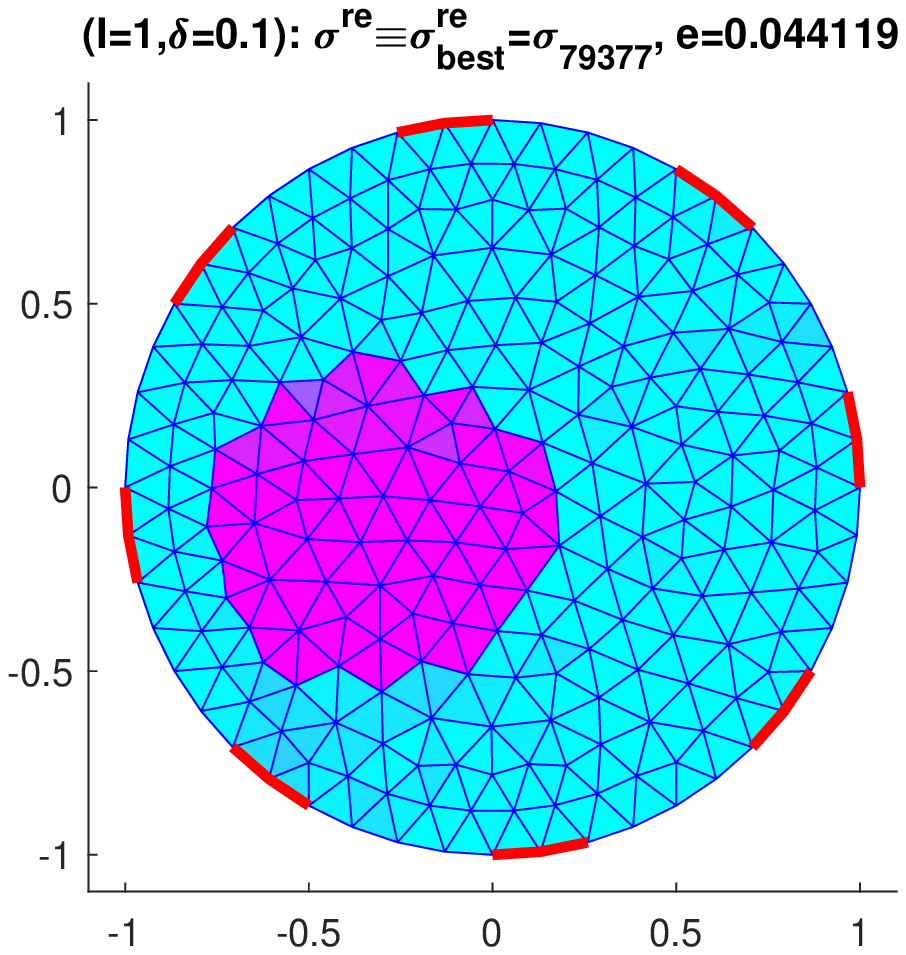}\\
\includegraphics[width=0.29\textwidth]{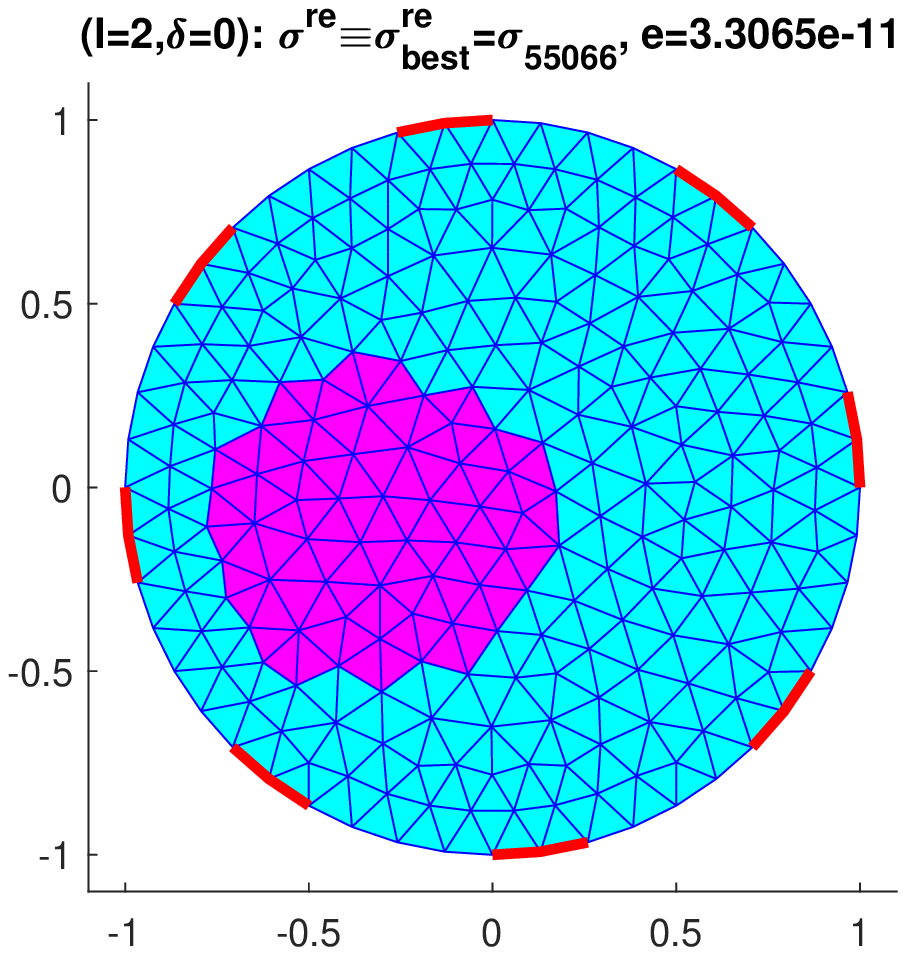}
\includegraphics[width=0.29\textwidth]{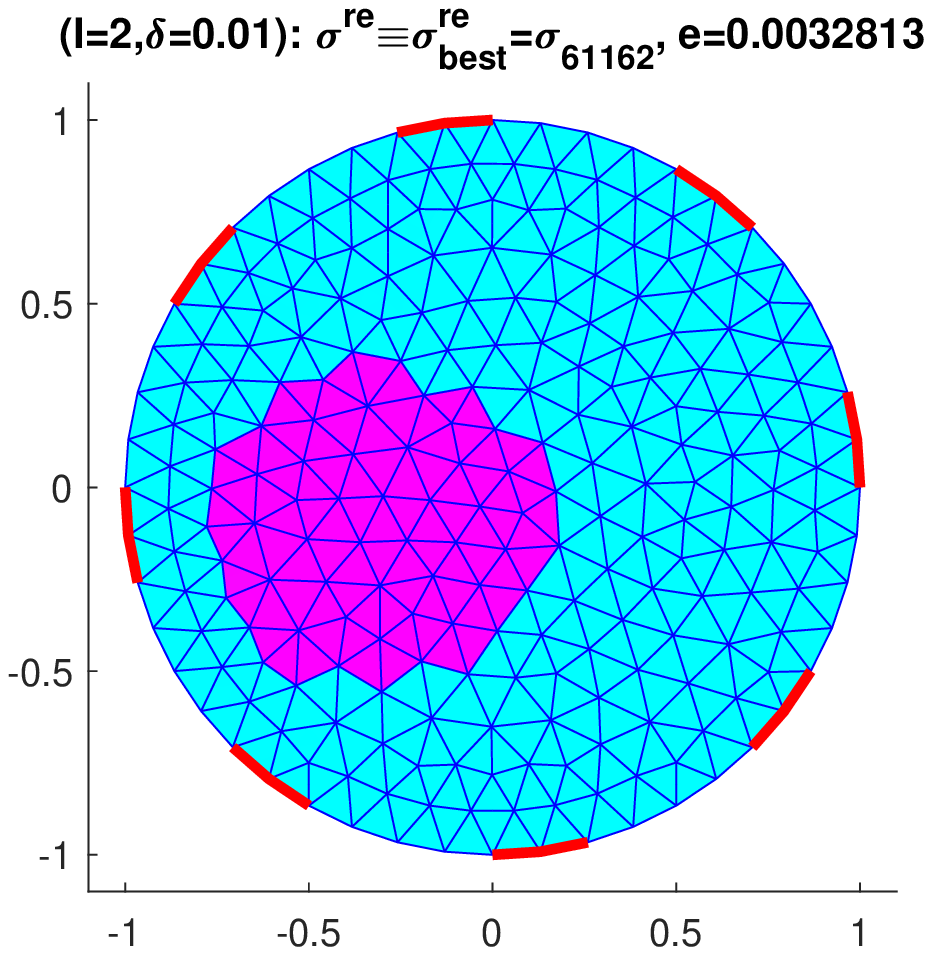}
\includegraphics[width=0.29\textwidth]{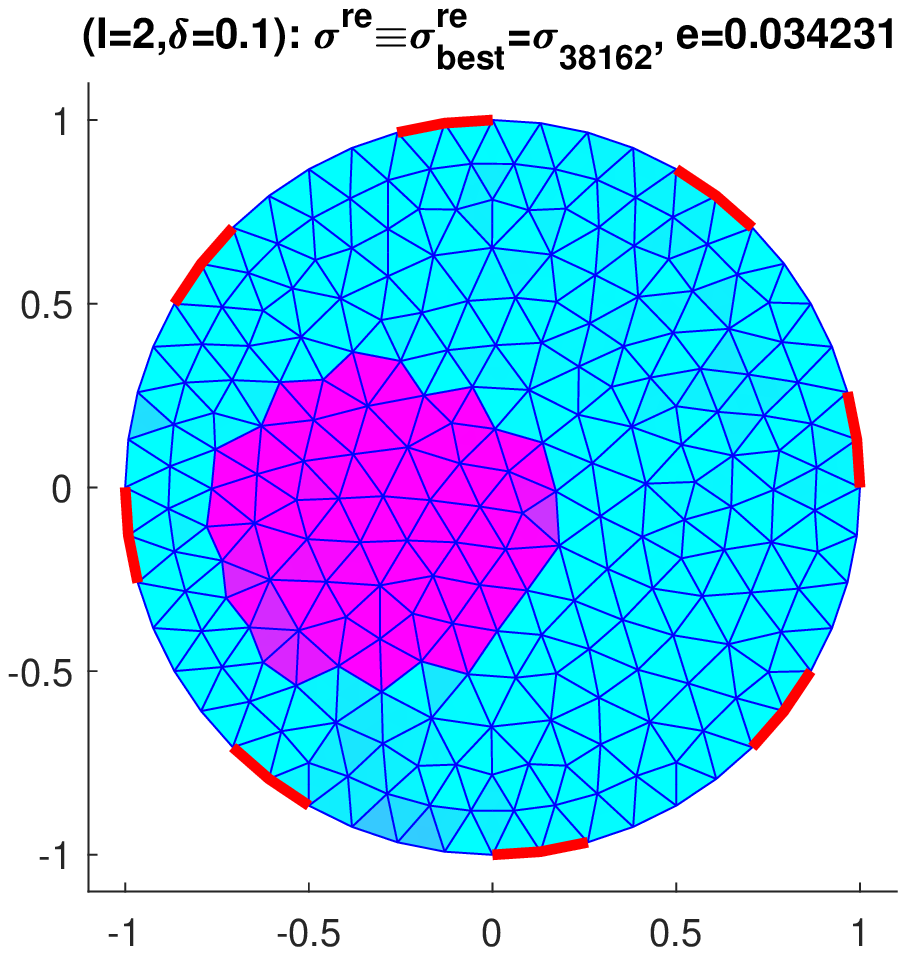}\\
\includegraphics[width=0.29\textwidth]{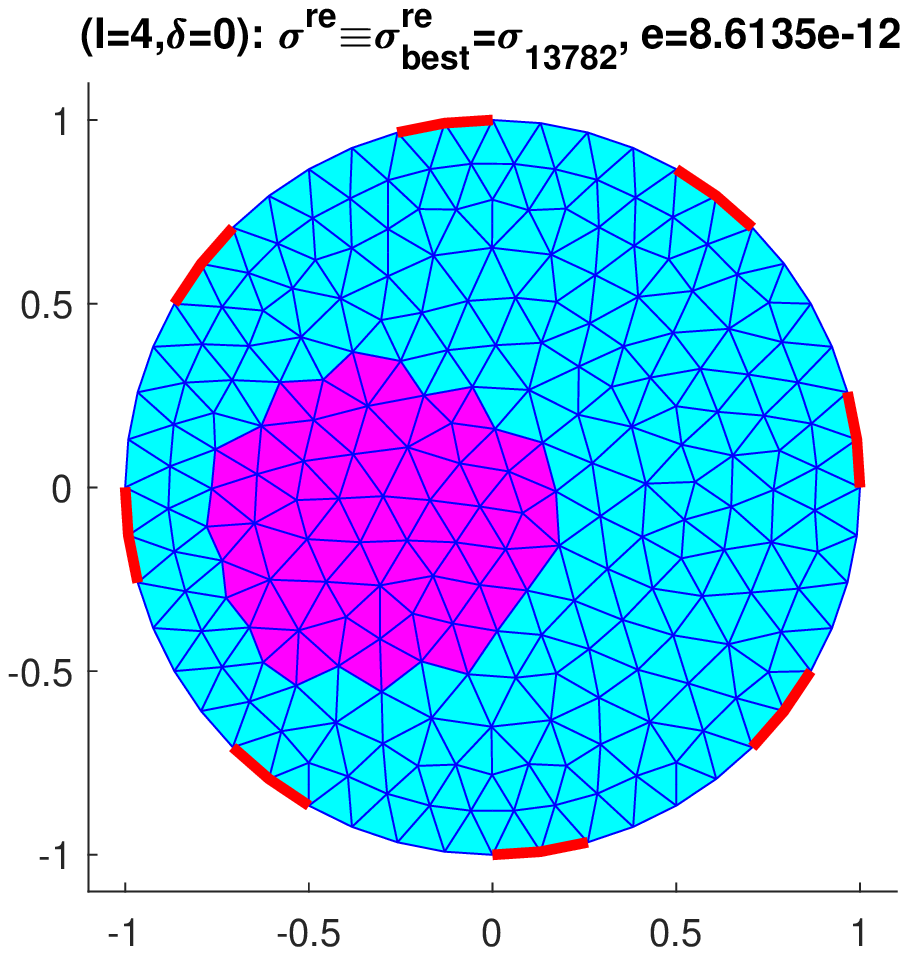}
\includegraphics[width=0.29\textwidth]{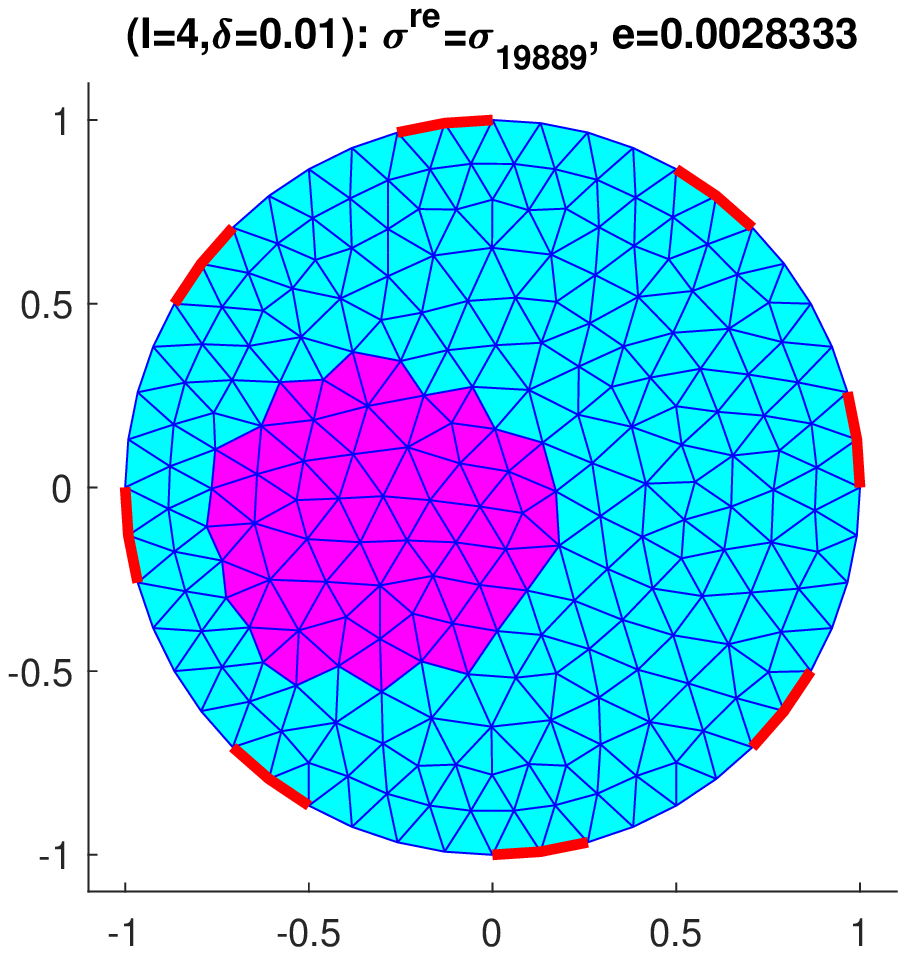}
\includegraphics[width=0.29\textwidth]{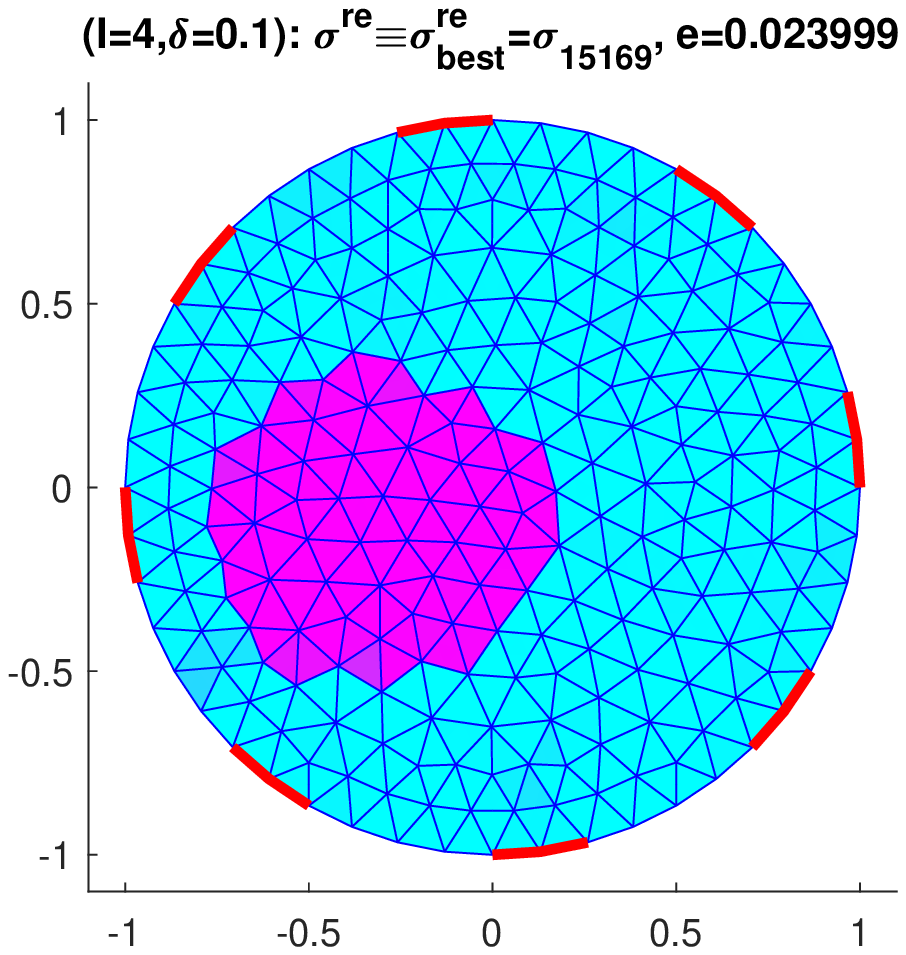}\\
\includegraphics[width=0.29\textwidth]{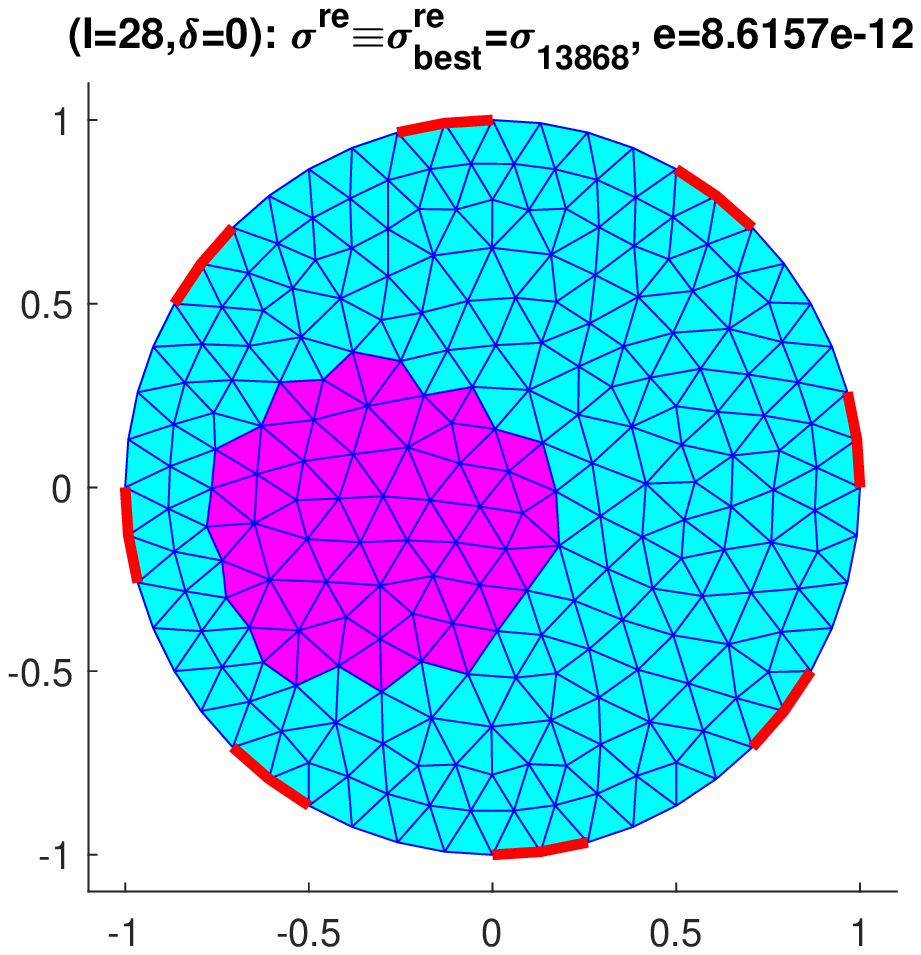}
\includegraphics[width=0.29\textwidth]{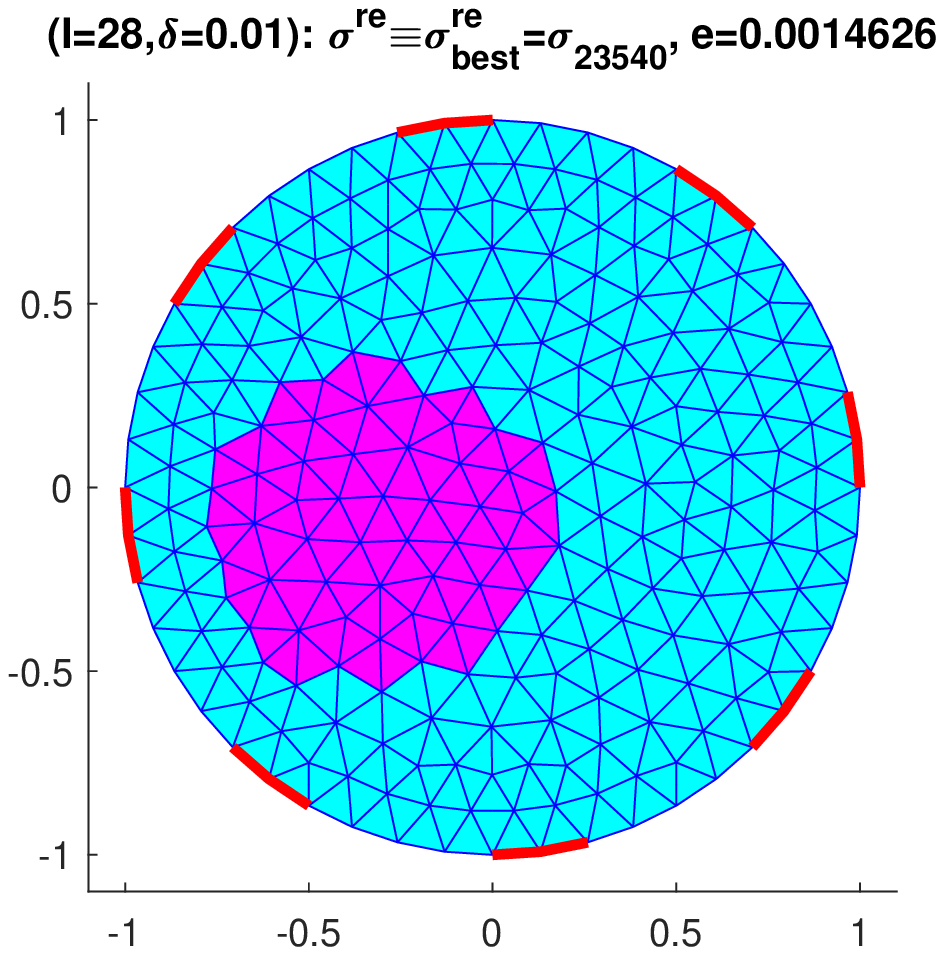}
\includegraphics[width=0.29\textwidth]{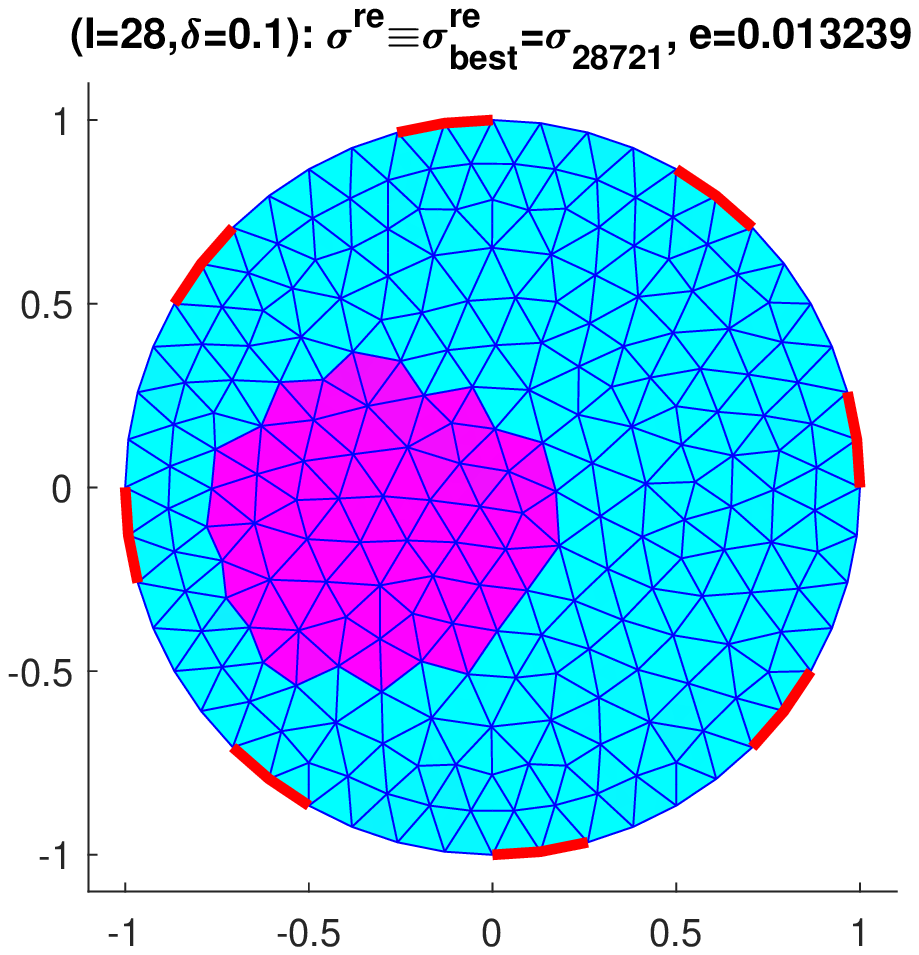}\\
\caption{Reconstructions of $\sigma$ from eliminating-$\Phi-\Psi$ version of cost function IAT \eqref{eq:costfunJ_IAT_PhiPsi}, in cases $I=1$, $I=2$, $I=4$, $I=28$ (top to bottom) for $\delta=0$, $\delta=0.01$, $\delta=0.1$ (left to right). 
\label{fig:sigma12IAT_eliminatingPhiPsiVer}}
\end{figure}

Finally, in Figures~\ref{fig:sigmaIAT_all-at-onceVer}, \ref{fig:sigma12IAT_eliminatingsigmaVer} and \ref{fig:sigma12IAT_eliminatingPhiPsiVer}, we display some pictures of reconstructions.


\subsection{Numerical results for EIT}
\begin{figure}
\begin{center}
\includegraphics[width=0.29\textwidth]{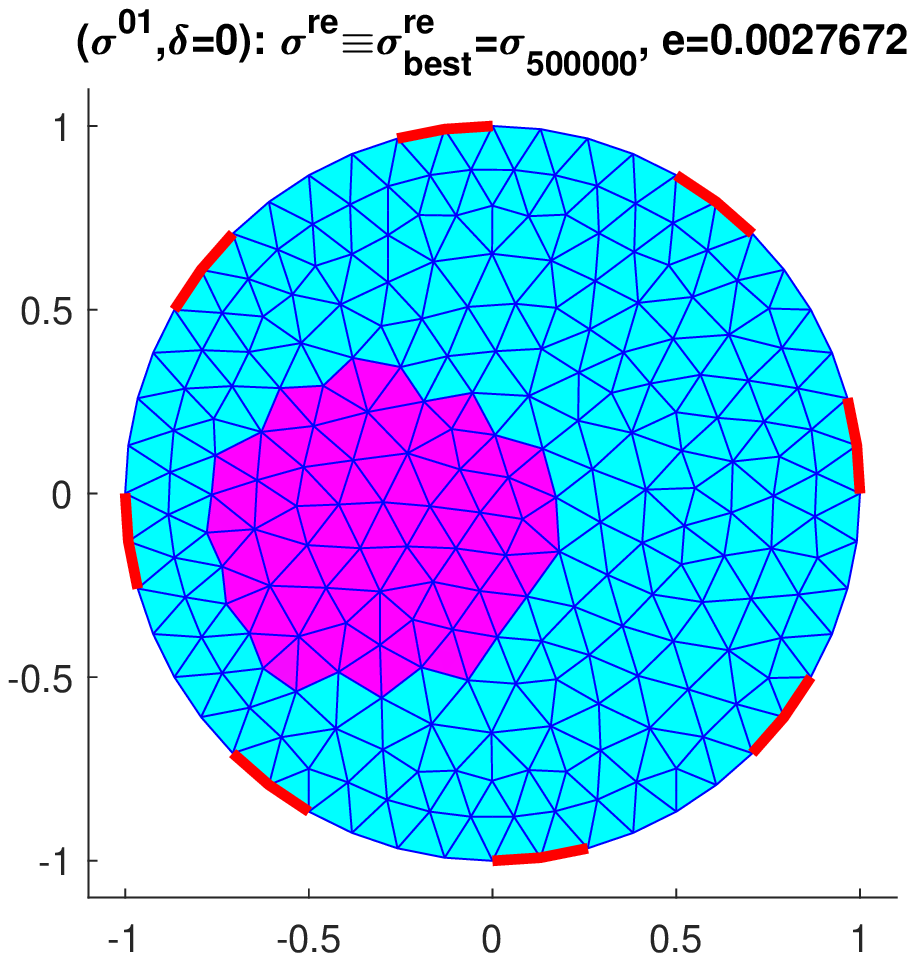}
\includegraphics[width=0.29\textwidth]{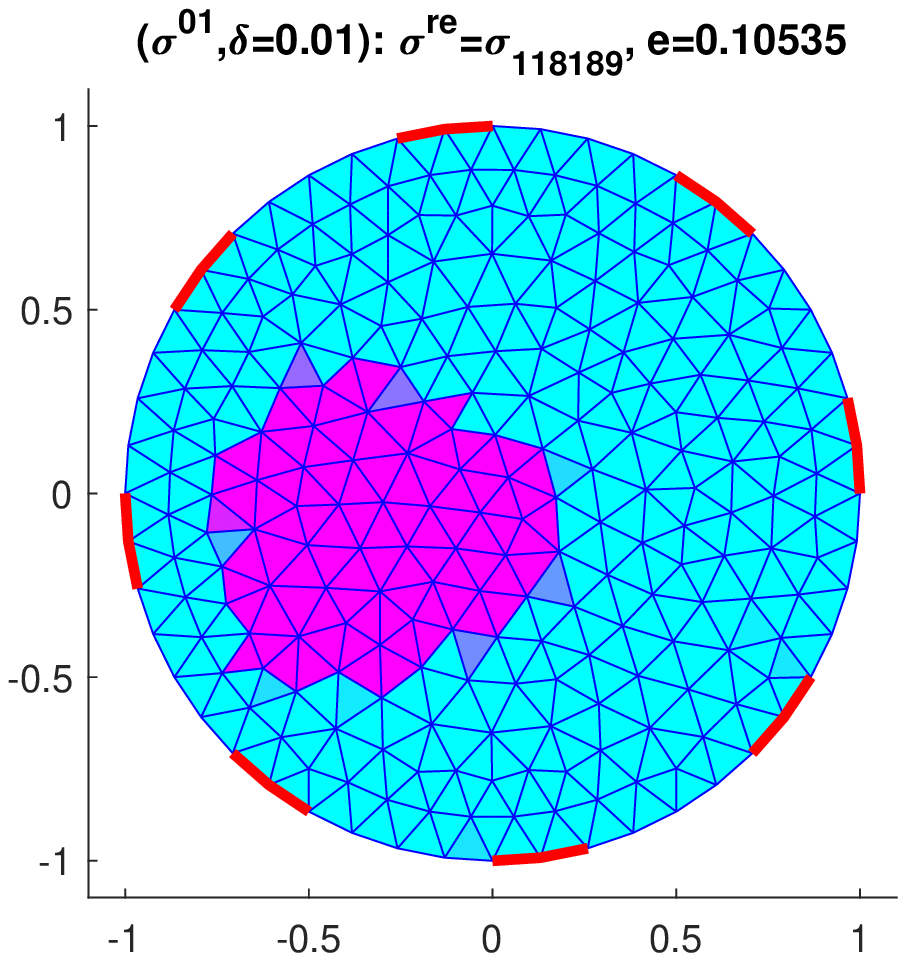}
\includegraphics[width=0.29\textwidth]{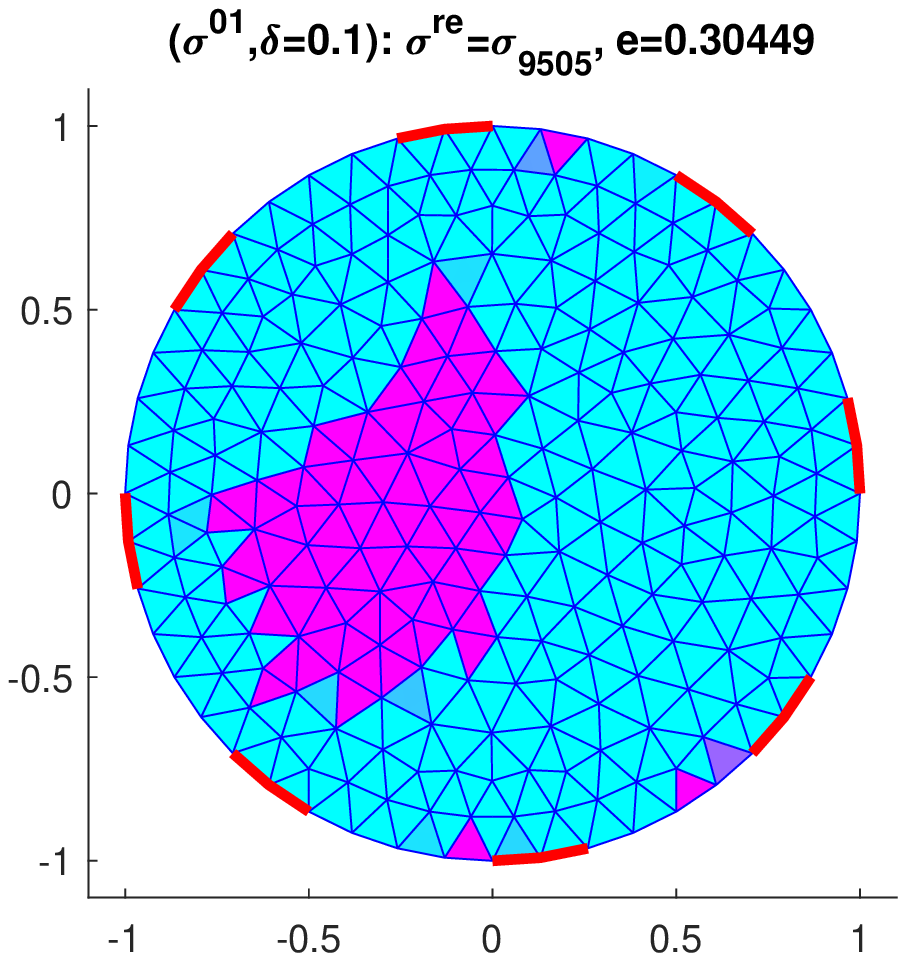}\\
\end{center}
\caption{Reconstructions of $\sigma$ from eliminating-$(\Phi,\Psi)$ version of cost function EIT \eqref{eq:costfunJ_EIT_PhiPsi}, in case $I=28$, 
$\delta=0$, $\delta=0.01$, $\delta=0.1$ (left to right).. 
\label{fig:sigma28EIT_eliminatingPhiPsiVer_sigma01}}
\end{figure}

We close with a few pictures of reconstructions for EIT. Here, for obvious identifiabilty reasons, it is necessary to use all measurements $I=28$. Moreover, the all-at-once and eliminating-$\sigma$ versions failed to converge, so we here only provide results with the classical reduced version of EIT corresponding to \eqref{eq:costfunJ_EIT_PhiPsi}. 
Starting from the constant value $\sigma_{0} = \frac{1}{2} \big( \underline{\sigma} + \overline{\sigma} \big)$ we obtain the reconstructions in Figure~\ref{fig:sigma28EIT_eliminatingPhiPsiVer_sigma01} for noise levels of zero, one and ten per cent. Note that in view of the exponential ill-posedness of this inverse problem, the quality of reconstructions is more than reasonable for this level of data contamination.

More details in particular on numerical tests for EIT can be found in the PhD thesis \cite{ThesisKha2021}.

\section{Conclusions and remarks}
In this paper we have provided convergence results on the iterative solution methods (gradient or Newton type) for minimization based formulations of inverse problems. 
We apply these to the identification of a spatially varying diffusion coefficient in an elliptic PDE from different kinds of measurements, in particular corresponding to the electrical impedance tomography EIT and the impedance acoustic tomography IAT problem, for which we also provide numerical tests.
Future work will, e.g., be concerned with investigations on the convexity conditions: How can an additive combination of functionals and constraints help to satisfy them, e.g., for EIT or IAT?

\section*{Appendix}
\begin{lemma} \label{lem:opialdiscr} (Opial, discrete)
Let $S$ be a non empty subset of a Hilbert space $X$, and $(x_k)_{k\in\mathbb{N}}$ a sequence of elements of $X$. Assume
that
\begin{enumerate}
\item[(i)] for every $z \in S$, $\lim_{k\to\infty}\|x_k - z\|$ exists;
\item[(ii)] every weak sequential limit point of $(x_k)_{k\in\mathbb{N}}$, as $k\to\infty$, belongs to $S$.
\end{enumerate}
Then $x_k$ converges weakly as $k\to\infty$ to a point in $S$.
\end{lemma}
\begin{lemma} \label{lem:opialcont} (Opial, continuous)
Let $S$ be a non empty subset of a Hilbert space $X$, and $x:[0,\infty)\to X$ a map. Assume
that
\begin{enumerate}
\item[(i)] for every $z \in S$, $\lim_{t\to\infty}\|x(t) - z\|$ exists;
\item[(ii)] every weak sequential limit point of $x(t)$, as $t\to\infty$, belongs to $S$.
\end{enumerate}
Then $x(t)$ converges weakly as $t\to\infty$ to a point in $S$.
\end{lemma}

\section*{Acknowledgments}
The work of the first author was supported by the Austrian Science Fund {\sc fwf} under the grants P30054 and DOC 78.

\end{document}